\documentclass[journal]{new-aiaa}

\usepackage{graphicx,subcaption}
\usepackage{amsmath}
\usepackage[mathscr]{euscript}
\usepackage{cases}
\usepackage{dsfont}
\usepackage{breqn}
\usepackage{physics}
\usepackage{algpseudocodex}
\usepackage{algorithm}
\usepackage{bm}
\usepackage{booktabs}
\usepackage{xr}
\usepackage{multirow}
\usepackage{siunitx}
\usepackage{longtable}
\AtBeginDocument{\RenewCommandCopy\qty\SI}
\restylefloat{algorithm}

\newcommand{\R}{\mathbb{R}}

\newcommand{\E}{\mathbb{E}}
\newcommand{\Symmetric}{\mathbb{S}}
\newcommand{\diff}{\mathrm{d}}

\newcommand{\cov}[1]{\mathrm{Cov}[#1 ]}
\let\P\relax
\newcommand{\P}[1]{\mathbb{P}\left[#1\right]}
\let\E\relax
\newcommand{\E}[1]{\mathbb{E}\left[#1\right]}

\newcommand{\suchthat}{\mathrm{s.t.}\quad}
\newcommand{\normal}{\mathcal{N}}

\newcommand{\blkdiag}[1]{\mathrm{blkdiag}(#1)}

\DeclareMathOperator*{\minimize}{minimize}

\usepackage{cleveref}
\crefname{align}{Eq.}{Eqs.}
\crefname{equation}{Eq.}{Eqs.}
\crefname{figure}{Fig.}{Figs.}
\crefname{table}{Table}{Tables}
\crefname{theorem}{Theorem}{Theorems}
\crefname{definition}{Definition}{Definitions}
\crefname{lemma}{Lemma}{Lemmas}
\crefname{remark}{Remark}{Remarks}
\crefname{assumption}{Assumption}{Assumptions}
\crefname{proof}{Proof}{Proofs}
\crefname{algorithm}{Algorithm}{Algorithms}
\crefname{problem}{Problem}{Problems}
\crefname{proposition}{Proposition}{Propositions}
\crefname{corollary}{Corollary}{Corollaries}
\crefname{section}{Section}{Sections}

\newtheorem{problem}{Problem}
\newtheorem{remark}{Remark}

\allowdisplaybreaks[1] 

\newcommand\blfootnote[1]{%
  \begingroup
  \renewcommand\thefootnote{}\footnote{#1}%
  \addtocounter{footnote}{-1}%
  \endgroup
}

\def\showChanges{0} 

\newcommand{\highlight}[1]{%
    \ifnum\showChanges=1
        \textcolor{red}{#1}%
    \else
        #1%
    \fi
}
\newcommand{\mathhl}[1]{
    \ifnum\showChanges=1
        \mathcolor{red}{#1}
    \else
        #1
    \fi
}

\usepackage{marginnote}
\usepackage{tcolorbox}
\newcommand{\margincomment}[1]{%
    \ifnum\showChanges=1
    \textcolor{blue}{\setlength{\baselineskip}{10pt}\footnotesize\raggedright\marginnote{#1}}%
    \else
    \fi
}

\title{Robust Cislunar Low-Thrust Trajectory Optimization under Uncertainties via Sequential Covariance Steering \blfootnote{An earlier version of this paper was presented as paper 24-379 at the 2024 AAS/AIAA Astrodynamics Specialist Conference, Broomfield, CO, August 11-15 2024.}}

\begin{document}

\author{Naoya Kumagai\footnote{Ph.D. Student, School of Aeronautics and Astronautics, Purdue University, West Lafayette, Indiana, 47907; nkumagai@purdue.edu (Corresponding author)}
\ and Kenshiro Oguri\footnote{Assistant Professor, School of Aeronautics and Astronautics, Purdue University, West Lafayette, Indiana, 47907.}
}
\affil{Purdue University, West Lafayette, Indiana, 47907}

\maketitle{}

\begin{abstract}
Spacecraft operations are influenced by uncertainties such as dynamics modeling, navigation, and maneuver execution errors. 
Although mission design has traditionally incorporated heuristic safety margins to mitigate the effect of uncertainties, particularly before/after crucial events, it is yet unclear whether this practice will scale in the cislunar region, which features locally chaotic nonlinear dynamics and involves frequent lunar flybys. 
This paper applies chance-constrained covariance steering and sequential convex programming to simultaneously design an optimal trajectory and trajectory correction policy 
that can probabilistically guarantee safety constraints under the assumed physical/navigational error models. 
The results show that the proposed method can effectively control the state uncertainty in a highly nonlinear environment.
The framework allows faster computation and lossless \highlight{convexification of linear covariance propagation} compared to existing methods, enabling a rapid and accurate comparison of $\Delta V_{99}$ costs for different uncertainty parameters.
We demonstrate the algorithm on several transfers in the Earth-Moon Circular Restricted Three Body Problem. 
\end{abstract}
\section{Introduction}

Mission operations are subject to various uncertainties, such as unmodeled accelerations, orbit determination, maneuver execution errors, missed thrust events, and launch errors. Especially in cislunar space, a challenging environment that features greater dynamics nonlinearities that skew any Gaussian state uncertainties into non-Gaussian, there is an increasing interest in effectively characterizing the uncertainties through orbit determination \cite{qi_analysis_2023} and uncertainty quantification \cite{sharan_accurate_2023}. Also of interest is the design of control strategies that address these uncertainties. 

Compared to traditional mission design which involves heuristic margins for increased safety, especially after crucial events such as planetary gravity assists and launch, operations in cislunar space will encounter far more frequent crucial events, the representative one being lunar flybys.
With a higher frequency of such events, heuristic safety margins may be overly time-consuming to incorporate or even make the transfer designs infeasible. 
These scenarios call for a higher level of automation in the safe trajectory design process that introduces less conservativeness in the safety margins. 
The trajectory and control (i.e. trajectory correction maneuver (TCM)) policy should ideally be designed with the underlying uncertainty models in mind to reduce conservativeness.

A promising and rapidly expanding approach for designing safe controllers that can also provide non-conservative solutions is stochastic optimal control (SOC), which directly considers the state probability distribution in the optimization process. 
Several recent works have addressed SOC in the context of deep-space trajectory design, with varying uncertainty quantification (UQ) and optimization methods. 
Differential dynamic programming (DDP) and unscented transform (UT) are combined in \cite{ozaki_stochastic_2018,ozaki_tube_2020} to design closed-loop feedback policies for low-thrust trajectory design. While an explicit control policy is unavailable from these works and fitting from sigma points is required to retrieve the policy, \cite{yuan_uncertainty-resilient_2024} provides an explicit linear feedback policy. 
Greco et al. \cite{greco_robust_2022} includes orbit determination errors and uses nonlinear programming (NLP) for constructing robust open-loop impulsive trajectories. The algorithm utilizes a Gaussian Mixture Model (GMM) to approximate the state distribution, offering a high-fidelity UQ capability not seen in other methods. At the same time, the method relies on finite differencing inside the NLP solver; the complexity of the problem makes it seemingly difficult to obtain analytical gradients to further improve the method. 
Methods that use linear covariance analysis (LinCov) as a lower-level covariance propagation technique, in combination with various optimization techniques at the top level, have also been proposed in the literature. 
\highlight{For example, a genetic algorithm \cite{jinRobustTrajectoryDesign2020}, particle swarm optimization \cite{cavesmithAnglesOnlyRobustTrajectory2024}, Sweeping Gradient Method \cite{margolisRobustTrajectoryOptimization2024,margolisSweepingGradientMethod2023a}, and forward-backward shooting \cite{vargheseNonlinearProgrammingApproach2025} are used for the top-level optimization.} 
Several methods based on primer vector theory also exist; for example, \cite{oguri_stochastic-primer_2022} considers the covariance as an additional state variable, and \cite{sidhoumRobustLowThrustTrajectory2024} proposes a control law that explicitly modifies the costate. 

Methods based on sequential convex programming (SCP), dynamics linearization, and linear covariance propagation are proposed in \cite{ridderhof_chance-constrained_2020,benedikter_convex_2022,oguri_stochastic_2022,oguri_chance-constrained_2024}. Ridderhof et al. 
\cite{ridderhof_chance-constrained_2020} was one of the first works of the SCP approaches; 
Oguri and Lantoine \cite{oguri_stochastic_2022} additionally consider orbit determination and mass uncertainty, as well as multi-phase mission design. 
Benedikter et al.\cite{benedikter_convex_2022} proposes a different convexification scheme. 
Although its proof of lossless convexification was shown to be incorrect \cite{rapakoulias_comment_2023}, similar methods have been recently proposed for linear systems \cite{liu_optimal_2023,pilipovsky_computationally_2023}.
These linearization-based approaches can be broadly classified into two categories: 
methods that explicitly define the state covariance as decision variables \cite{benedikter_convex_2022}, 
which we term the \textit{full covariance} formulations, and methods that implicitly optimize for the covariance by constructing a large block matrix that is analogous to the Cholesky decomposition of covariance matrices \cite{ridderhof_chance-constrained_2020,oguri_stochastic_2022}
, which we term the \textit{block Cholesky} formulations. 
An inherent drawback in the \textit{block Cholesky} formulation is the required computational effort; although each subproblem is a convex program, it involves solving a large Linear Matrix Inequality,
the size of which scales quadratically with the number of discretization nodes. 
On the other hand, the \textit{full covariance} formulation subproblem is computationally more efficient due to the empirical linear growth of the computational complexity with respect to the number of nodes, which has been demonstrated for linear systems \cite{pilipovsky_computationally_2023}. 
Recent works on rendezvous and proximity operations consider rotational motion \cite{zhangStochasticTrajectoryOptimization2023,takuboMultiplicativeApproachConstrained2024} and navigation-enhancing maneuvers \cite{raChanceConstrainedSensingOptimalPath2024}.
\highlight{Although many works model unknown disturbances via additive Brownian motion,  
Gaussian processes \cite{garzelliGaussianProcessBasedCovariance2025} and multiplicative noise \cite{liuOptimalCovarianceSteering2024,benedikterCovarianceControlEarthtoMars2024} are also proposed.} 
\highlight{In addition to position and velocity uncertainties, mass uncertainty is considered in \cite{benedikterConvexApproachCovariance2023,benedikterConvexApproachStochastic2022,ridderhofMinimumFuelClosedLoopPowered2021,oguri_stochastic_2022}.}
Works that consider non-Gaussian uncertainties via Gaussian mixture models for linear systems are also emerging \cite{balciDensitySteeringGaussian2024,kumagaiChanceConstrainedGaussianMixture2024,knaupComputationallyEfficientCovariance2023}.
To summarize the rapidly growing number of SOC approaches, we provide a table for comparison in the Appendix.

SOC has been demonstrated in cislunar space, mainly for stationkeeping on periodic orbits in the Circular Restricted Three Body Problem (CR3BP)\highlight{\cite{oguri_chance-constrained_2024,gettatelliStationkeepingEarthMoonL22023}}. However, there is a gap in the literature on applications to general transfers. This paper adds to the literature of SCP-based SOC methods, with applications to designing a nominal robust trajectory and TCM policy for transfers in cislunar space. 
The main contributions are threefold. First, we show one of the first demonstrations of the SCP-based SOC method in cislunar space which is applicable for both low-thrust and impulsive propulsion systems. The method is capable of modeling unmodeled dynamics, maneuver execution errors, orbit determination errors, and initial insertion errors, and computing \highlight{a} state estimate-linear TCM policy that can explicitly account for the state covariance. The SOC approach is shown to effectively control the state uncertainty in an environment featuring high nonlinearity. 
Numerical examples of \highlight{transfer between Distant Retrograde Orbits (DROs) and another from 
a Near Rectilinear Halo Orbit (NRHO) to a halo orbit} demonstrate the effectiveness of the method. 
Second, we formulate a variant of the \textit{full covariance} formulation that formulates each subproblem as a convex problem which can be solved more efficiently than \textit{block Cholesky} formulations and is exact in the covariance propagation (cf. ~\cite{benedikter_convex_2022}). We show a method of imposing 2-norm upper bound constraints on the control input, which shows smaller conservativeness than recently proposed methods. 
Third, analyzing the nominal trajectory obtained from the optimization under uncertainties via the local Lyapunov exponent
shows that the resulting trajectory has better local stability compared to trajectory optimization that does not consider uncertainties. 

The paper is outlined as follows: 
\cref{sec:problem-formulation} outlines the force model (CR3BP), uncertainty models, and the $\Delta V_{99}$ minimization problem formulation.
\cref{sec:belief-space-dynamics} derives the filter model and the belief-space dynamics.
\cref{sec:solution_via_SCP} shows the solution to the $\Delta V_{99}$ minimization problem under the assumed models via convexification and SCP. 
\cref{sec:numerical_simulations} shows the results of applying the algorithms to several transfers in the cislunar region and an analysis via the local Lyapunov exponent.

\textit{Notation}:
$\R, \R^n, \R^{n\times m}, \Symmetric^n$ represent the space of real numbers, $n$-dimensional real vectors, $n$-by-$m$ real matrices, and symmetric $n$-by-$n$ matrices, respectively. 
$\P{\cdot}, \E{\cdot}, \cov{\cdot}, \norm{\cdot}, \norm{\cdot}_\infty, \lambda_{\max}(\cdot), \tr(\cdot)$ calculate the probability, expectation, covariance, 2-norm, infinity-norm, largest eigenvalue, and trace, respectively.\
 $A^{1/2}$ refers to any matrix such that $A^{1/2} (A^{1/2}) ^\top = A$.
 A random vector $\bm{\xi}$ with normal distribution of mean $\bm{\mu}$ and covariance matrix $\Sigma$ is denoted $\bm{\xi}\sim\normal(\bm{\mu}, \Sigma)$.
 The notation $\square_{1:k}$ denotes all elements from the first to the $k$-th element of the sequence.
 $\blkdiag{\cdot, \cdot, \cdots}$ denotes a block diagonal matrix with the arguments as diagonal blocks.
 $\mathcal{Q} (\xi; p)$ denotes the $p$-quantile of the random variable $\xi$.
\section{Problem Formulation} \label{sec:problem-formulation}
\begin{figure}[t]
    \centering
    \includegraphics[width=\textwidth]{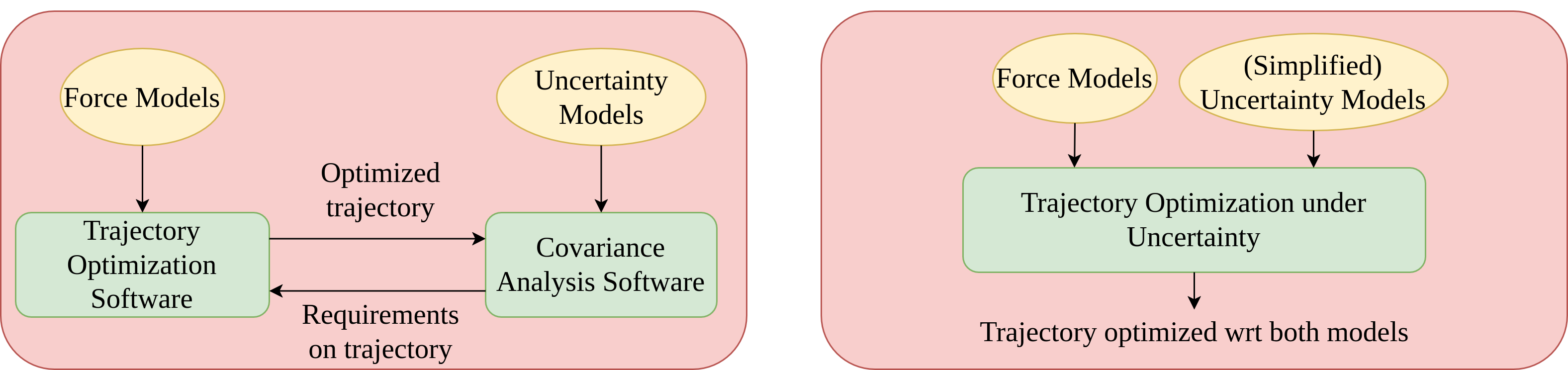}
    \caption{Comparison of a typical mission design process (left) vs trajectory optimization process that also takes into account the uncertainty models (right).}
    \label{fig:comparison}
\end{figure}

\margincomment{This paragraph has been moved from the introduction.}
\highlight{
The nature of our approach requires both the force model (e.g. gravity and thrust)
and the uncertainty model (e.g. navigation and control execution) to be considered in the optimization process.
}%
Figure \ref{fig:comparison} shows a comparison of a typical mission design process and trajectory optimization that takes into account the uncertainty models.
The traditional mission design process involves an iteration process between the trajectory design 
and uncertainty analysis, where uncertainty analysis may inform the trajectory design on 
where to add safety margins. The trajectory is then designed with the safety margins, and the
uncertainty analysis is performed again to check if the safety margins are sufficient;
generally, this process is repeated several times, which can be time-consuming.
On the other hand, the proposed method involves a trajectory optimization process that also takes into account the uncertainty models. 
While the uncertainty models that are used in the optimization process may not be of the same fidelity
as dedicated uncertainty analysis tools, the optimization process can provide a more direct way 
to design a trajectory that is robust to uncertainties without the need for an iterative process.
This section describes both the force and uncertainty models for our specific application.

\subsection{Force Model: Circular Restricted Three Body Model} \label{sec:dynamics}
For preliminary trajectory design in cislunar space, the Circular Restricted Three Body Problem (CR3BP) is a commonly-used model that can incorporate the gravitational effects of both the Earth and Moon on the spacecraft. We use this model for analysis in this work. 

The CR3BP describes the motion of a particle under the gravitational influence of two bodies $P_1$ and $P_2$, termed the primaries, with mass $m_1$ and $m_2$, respectively. The primaries are assumed to be moving in a circular orbit around their barycenter. The particle models the spacecraft, which has an approximately infinitesimal mass compared to the primaries. We take a right-handed Cartesian coordinate frame with the origin at the barycenter of the two primaries. The $X$-axis points towards $P_2$, the $Z$-axis is aligned with the orbital angular momentum vector of the primaries, and the $Y$-axis completes the right-hand system. As is done traditionally, we study the CR3BP dynamics in a length and time-nondimensionalized form. The second-order differential equations of motion (EOM) for the particle are \cite{pavlak_trajectory_2013}:
\begin{subequations}
\begin{align}
\ddot{x}-2 \dot{y}-x & =-\frac{(1-\mu)(x+\mu)}{d_1^3}-\frac{\mu(x-1+\mu)}{d_2^3} \\
\ddot{y}+2 \dot{x}-y & =-\frac{(1-\mu) y}{d_1^3}-\frac{\mu y}{d_2^3} \\
\ddot{z} & =-\frac{(1-\mu) z}{d_1^3}-\frac{\mu z}{d_2^3}
\end{align}
\end{subequations}
where
\begin{equation}
    d_1 = \sqrt{(x+\mu)^2 + y^2 + z^2}, \quad d_2 = \sqrt{(x-1+\mu)^2 + y^2 + z^2}, \quad \mu = \frac{m_2}{m_1+m_2}
\end{equation}
By defining the state $\bm{x}$ of the spacecraft as its position and velocity, i.e. $\bm{x} \in \mathbb{R}^{n_x}$, $\bm{x} = [x,y,z,\dot{x}, \dot{y}, \dot{z}]^{\top}$ with $n_x=6$, we can express the EOM in first-order differential form. When low-thrust propulsion applies  acceleration, the EOM can be written as the sum of the natural dynamics $\bm{f}_0(\bm{x})$ and the acceleration provided by the thrusters:
\begin{equation}
    \dot{\bm{x}} = \bm{f}(\bm{x}, \bm{u}) = \bm{f}_0(\bm{x}) + B \bm{u}
\end{equation}
where $B \in \R^{n_x \times n_u}$, and $\bm{u} \in \mathbb{R}^{n_u}$, with $n_u=3$ and 
\begin{align}
 \bm{f}_0(\bm{x}) = 
\begin{bmatrix}
    \dot{x}\\ \dot{y}\\ \dot{z} \\
    2 \dot{y}+ x - \frac{(1-\mu)(x+\mu)}{d_1^3}-\frac{\mu(x-1+\mu)}{d_2^3} \\
    - 2 \dot{x} + y - \frac{(1-\mu) y}{d_1^3}-\frac{\mu y}{d_2^3}\\
 -\frac{(1-\mu) z}{d_1^3}-\frac{\mu z}{d_2^3}
\end{bmatrix}, \quad 
 B = \begin{bmatrix}
     \bm{0}_{3\times3} \\ I_3
 \end{bmatrix}
\end{align}

\subsection{Uncertainty Models}

\subsubsection{Initial State Uncertainty} 
We assume knowledge about the distribution of the initial state, in terms of an initial \highlight{prior} state estimate $\hat{\bm{x}}_{0^-}$ \highlight{with}
estimation error $\tilde{\bm{x}}_{0^-}$, both Gaussian-distributed vectors:
\begin{equation} \label{eq:initial-distributional}
    \hat{\bm{x}}_{0^-} \sim \normal({\bm{\mu}}_0, \hat{P}_{0^-}), \quad \tilde{\bm{x}}_{0^-} \sim \normal(\bm{0}_{n_x}, \tilde{P}_{0^-})
\end{equation}
where ${\bm{\mu}}_0$ and $\hat{P}_{0^-}, \tilde{P}_{0^-} \succeq 0$ are known. 

\subsubsection{Unmodeled Accelerations}
We model unknown/unmodeled accelerations as a Brownian motion, which perturbs the spacecraft as an additional term in the dynamics.
The spacecraft dynamics are now expressed as a stochastic differential equation (SDE):
\begin{equation} \label{eq:SDE}
    \dd{\bm{x}} = \bm{f}(\bm{x},\bm{u}) \dd{t} + \bm{g}(\bm{x},\bm{u}) \dd{\bm{w}}
\end{equation}
where  $\dd{\bm{w}}(t) \in \mathbb{R}^{n_w}$ is standard Brownian motion, i.e., $\mathbb{E}[\diff \bm{w}]=\bm{0}_{n_w}$ and $\mathbb{E}\left[\diff \bm{w}(t) \diff \bm{w}^{\top}(t)\right]=I_{n_w} \diff t$.

\subsubsection{Maneuver Execution Error}
Regardless of the propulsion system, planned control maneuvers are subject to erroneous execution. In this work, we use the Gates model \cite{gates_simplified_1963} which involves four parameters $\{ \sigma_1, \sigma_2, \sigma_3, \sigma_4\}$. These respectively represent the fixed magnitude [\unit{km/s^2} (\unit{km/s})], proportional magnitude [\unit{\%}], fixed pointing [\unit{km/s^2}  (\unit{km/s})], and proportional pointing [\unit{deg}] errors for a spacecraft equipped with a low-thrust (impulsive) propulsion system. For a planned control input $\bm{u}_k$, it models the additional control error $\tilde{\bm{u}}_k$ as a Gaussian noise \cite{oguri_chance-constrained_2024}:
\begin{subequations}
\begin{align}
    &\tilde{\bm{u}}(t) = \tilde{\bm{u}}_k, \quad t \in [t_k, t_{k+1}) \\
    &\tilde{\bm{u}}_k = G_{\mathrm{exe}}(\bm{u}_k) \bm{w}_{\mathrm{exe},k}, \quad \bm{w}_{\mathrm{exe},k} \sim \mathcal{N}(0, I_3) \\
    &G_{\mathrm{exe}}(\bm{u}_k) = T(\bm{u}_k) P_{\mathrm{gates}}^{1/2} (\bm{u}_k) 
\end{align}
\end{subequations}
where the rotation matrix $T(\bm{u}_k)$ and the covariance matrix $P_{\mathrm{gates}}(\bm{u}_k)$ are defined as 
\begin{align}
    T(\bm{u}) = \begin{bmatrix}
        \hat{\bm{s}} & \hat{\bm{e}} & \hat{\bm{z}}
    \end{bmatrix}, \quad 
    \hat{\bm{z}}=\frac{\bm{u}}{\norm{\bm{u}}}, \quad \hat{\bm{e}}=\frac{[0,0,1]^{\top} \times \hat{\bm{z}}}{\norm{[0,0,1]^{\top} \times \hat{\bm{z}}}}, \quad \hat{\bm{s}}=\hat{\bm{e}} \times \hat{\bm{z}}
\end{align}
\begin{align}
    P_{\mathrm{gates}}(\bm{u}) = \begin{bmatrix}
        \sigma_p^2 & 0 & 0 \\ 
        0 & \sigma_p^2 & 0 \\
        0 & 0 & \sigma_m^2
    \end{bmatrix}, \quad
    \sigma_m^2 = \sigma_1^2 + \sigma_2^2 \norm{\bm{u}}^2, \quad 
    \sigma_p^2 = \sigma_3^2 + \sigma_4^2 \norm{\bm{u}}^2
\end{align}
\highlight{The orthonormal frame for the Gates model is defined with the normalized planned thrust direction $\hat{\bm{z}}$ as the third axis,
and the $\hat{\bm{s}}$ and $\hat{\bm{e}}$ axes lie in the plane perpendicular to the thrust direction. The rotation matrix $T(\bm{u})$ transforms 
the error vector $\bm{w}_{\mathrm{exe},k}$ from this frame to the frame in which the planned thrust direction
is represented, in this case, the CR3BP rotating frame.
}

\highlight{
    The Gates model is a specific case of control-dependent noise \cite{liuOptimalCovarianceSteering2024,benedikterCovarianceControlEarthtoMars2024};
    in this work, we model this control-dependence by leveraging the iterative solutions computed within the SCP framework.
    At each iteration, we compute the control error covariance with $G_{\mathrm{exe}}(\overline{\bm{u}}_k^*)$, where 
    $\overline{\bm{u}}_k^*$ is the \textit{reference} nominal control input at the $k$-th time step, obtained from the previous iteration.
    Strictly speaking, we are neglecting the effect of control error based on the feedback term in the control policy described in 
    \cref{sec:belief-space-dynamics}; however, this is numerically shown to be reasonable since the feedback term is small compared to the feedforward term.
}

\subsubsection{Orbit Determination Process}
We assume that the state information is partially available through a nonlinear discrete-time observation process at the 
observation times $t_0, t_1, \cdots, t_N$. The observation model is expressed as:
\begin{align} \label{eq:OD}
    \bm{y}_k= \bm{f}_{\mathrm{obs}} (\bm{x}_k)+ \bm{g}_{\mathrm{obs}} (\bm{x}_k) \bm{v}_k
\end{align}
where $\bm{y}_k \in \R^{n_y}$ is the observation, and $\bm{v}_k \in \R^{n_y}$ is an i.i.d. standard Gaussian vector.  
Based on the observation process, we can estimate the state and its distribution using a navigation filter.
The filtering process $\mathcal{F}_k$ is assumed to use all previous measurements $\bm{y}_{0:k}$, the initial state estimate $\hat{\bm{x}}_{0^-}$,
\highlight{and the initial prior estimation error covariance $\tilde{P}_{0^-}$}.
Then, the distribution of the state estimate $\hat{\bm{x}}_k$ is written as 
\begin{equation}
    \hat{\bm{x}}_k \sim \mathcal{F}(\hat{\bm{x}}_{0^-}, \mathhl{\tilde{P}_{0^-}}, \bm{y}_{0:k})
\end{equation}

\subsection{Control Policy}
The control policy is a function of the state estimate, 
since we can only observe the state through the navigation filter:
\begin{equation}
    \bm{u}(t) = \pi(\hat{\bm{x}}, t), \quad t \in [t_0, t_f]
\end{equation}

\subsection{Chance/Distributional Constraints for Safe Trajectory Design}

\subsubsection{Chance Constraints}
When we assume a stochastic system with unbounded uncertainties such as Gaussian distributions, it is difficult to guarantee constraints on the state or control in a deterministic manner. Thus, a chance constraint generalizes the constraint $f_{\mathrm{safe}}(\bm{x}, \bm{u}) \leq 0$ to be probabilistic in the form
\begin{equation}
    \P{f_{\mathrm{safe}}(\bm{x}, \bm{u}) \leq 0} \geq 1 - \epsilon_{cc}
\end{equation}
where $\epsilon_{cc}$ represents the probability of chance-constraint violation, chosen as $0 <\epsilon_{cc} \ll 1$. 

This work deals with chance constraints in the form of 2-norm control constraints:
\begin{subequations} \label{eq:chance-constraints}
\begin{align} 
    &\P{\norm{\bm{u}_k} \leq u_{\max}} \geq 1 - \epsilon_{u}, && k = 0, \cdots,  N-1 \label{eq:chance-constraints-control}
\end{align}
\end{subequations}
where $u_{\max}$ is a scalar upper bound on the control input 2-norm. 
We can also formulate more general forms of chance constraints; readers are referred to \cite{oguri_chance-constrained_2024}. For example, keep-out zones can be imposed on the state to ensure a low probability of collision with any object.

\subsubsection{Maximum Covariance Constraint}
To limit the growth of the state uncertainty, we can impose a maximum covariance constraint on the state covariance
at any time step $k$, in the form of a linear matrix inequality (LMI) constraint:
\begin{equation}
    \cov{\bm{x}_k} \preceq P_{\max}, \quad k = 0, \cdots, N
\end{equation}\
where $P_{\max} \succ 0$ is a user-specified positive-definite matrix.
This has another benefit of maintaining the validity of the linearization of the covariance propagation model, 
which is discussed in the next section.

\subsubsection{Final Distributional Constraint}
Similar to the chance constraints, the final state must be constrained 
in a probabilistic manner. Although we can impose chance constraints on the final state,
we choose to impose a distributional constraint for the final state.
We constrain the mean to match a desired state $\bm{\mu}_f$
and covariance of the final state $\bm{x}(t_f)$ 
to be within a user-specified covariance matrix $P_f \succ 0$, i.e.
\begin{equation} \label{eq:final-distributional}
    \E{\bm{x}(t_f)} = \bm{\mu}_f, \quad \cov{\bm{x}(t_f)} \preceq P_f
\end{equation}

\subsection{Objective Function}
Since our control policy is a function of the state estimate, which is a random variable, 
the cost function $J$ is also a random variable. We consider the $p$-quantile of the cost function $J$ as the objective function:
\begin{equation}
    \minimize_\pi \quad \mathcal{Q}(J; p)
\end{equation}
In this work, we consider the $99\%$ quantile of the total fuel consumption as the objective function.
This value is commonly termed the $\Delta V_{99}$ in spacecraft mission design:
\begin{equation}
    J = \int_{t_0}^{t_f} \norm{\bm{u}(t)} \dd{t}, \quad p = 0.99
\end{equation}

\subsection{Original Problem Formulation}
We can now formulate a chance-constrained stochastic optimal control problem as follows:

\begin{problem}{Chance-Constrained $\Delta V_{99}$ Optimization (CCDVO) under Uncertainties} \label{prob:main}
\begin{subequations}
    \begin{align}
    \minimize_\pi \quad & \mathcal{Q} \left(\int_{t_0}^{t_f} \norm{\bm{u}(t)} \dd{t}; 0.99 \right)\\
    \suchthat \quad 
    & \dd{\bm{x}} = \bm{f}(\bm{x},\bm{u}) \dd{t} + \bm{g}(\bm{x},\bm{u}) \dd{\bm{w}} && \triangleright \text{Nonlinear stochastic dynamics} \\
    & \hat{\bm{x}}_{0^-} \sim \normal({\bm{\mu}}_0, \hat{P}_{0^-}), \quad \tilde{\bm{x}}_{0^-} \sim \normal(\bm{0}_{n_x}, \tilde{P}_{0^-})
    && \triangleright \text{Initial state distribution}\\
    & \P{f_{\mathrm{safe}}(\bm{x}, \bm{u}) \leq 0} \geq 1 - \epsilon_{cc} && \triangleright \text{Chance constraint}\\
    & \cov{\bm{x}(t)} \preceq P_{\max} \ \forall t \in [t_0, t_f]&& \triangleright \text{Maximum covariance} \\
    & \E{\bm{x}(t_f)} = \bm{\mu}_f, \quad \cov{\bm{x}(t_f)} \preceq P_f && \triangleright \text{Terminal distribution}\\
    & \bm{y}_k= \bm{f}_{\mathrm{obs}} (\bm{x}_k)+ \bm{g}_{\mathrm{obs}} (\bm{x}_k) \bm{v}_k && \triangleright \text{Observation model} \\
    & \hat{\bm{x}}_k = \mathcal{F}(\hat{\bm{x}}_{0^-}, \mathhl{\tilde{P}_{0^-}}, \bm{y}_{0:k}) && \triangleright \text{Navigation filter} \\
    & \bm{u}(t) = \pi(\hat{\bm{x}}, t) && \triangleright \text{Control policy}
\end{align}
\end{subequations}
\end{problem}

Exactly solving this nonlinear stochastic optimal control problem can be computationally intractable;
for example, the propagation of the state distribution through the nonlinear dynamics
requires solving the Fokker-Planck equation, which can be computationally expensive.
Hence, we make assumptions on the uncertainty propagation model, controller design, and navigation filter, 
which are valid when the uncertainty is maintained to be small throughout the time horizon.

\section{Belief-Space Dynamics and Filtered Density Propagation} \label{sec:belief-space-dynamics}

In this section, we outline the assumptions on the control policy, dynamics linearization, and filtered covariance propagation model.
By assuming these models, we can tractably design a closed-loop controller that explicitly takes into account the state distribution evolution.

\subsection{Control Policy}
First, we restrict the class of control policies to piecewise constant affine controllers, i.e.,
the control policy $\pi$ is constituted by a piecewise constant control input, $\bm{u}_{0:N-1}$, 
where $\bm{u}(t) \equiv \bm{u}_k$ for $t \in [t_k, t_{k+1})$, and the control input 
is a linear function of the state estimate $\hat{\bm{x}}_k$: 
\begin{equation} \label{eq:feedback}
    \bm{u}_k = \overline{\bm{u}}_k + K_k (\hat{\bm{x}}_k - \overline{\bm{x}}_k)
\end{equation}
where $\overline{\bm{u}}_k$ is the nominal (feedforward) control, and $K_k \in \R^{n_u \times n_x}$ is the feedback gain matrix, with
the term $K_k (\hat{\bm{x}}_k - \overline{\bm{x}}_k)$ representing the correction to the nominal control based on the state estimate, 
more commonly known as TCMs in spacecraft mission design.
\footnote{For simplicity, 
in this work, we assume that the control discretization intervals are the same 
as the filter update intervals. It is straightforward to extend the formulation to cases
where the interval lengths are integer multiples of each other; for example, 
requiring that the control between two nodes are constant can be imposed as a 
simple linear equality constraint.}
\begin{remark}
    \normalfont The linear feedback policy \cref{eq:feedback} can be regarded as a sub-optimal
    policy in the class of nonlinear trajectory correction policies. For example, state-of-the-art
    ground operations for deep-space missions typically recompute the control by 
    solving a nonlinear trajectory optimization problem at each time step. 
    However, since the linear feedback policy is sub-optimal, it can provide a 
    computationally tractable upper bound on the quantile cost, and can be used
    to estimate statistical costs for more complicated policies.
\end{remark}

\subsection{Dynamics Linearization and Discretization}
The uncertainty model we consider in this work takes an Extended Kalman Filter (EKF)-like model, where the 
first two statistical moments are considered to represent the state distribution. 
The state estimate mean is propagated nonlinearly through the nonlinear dynamics, 
while its covariance is propagated according to a linearized model. 
This assumption allows us to tractably design closed-loop controllers by solving an 
optimization that explicitly takes into account the state distribution evolution.
First, we review
the process of dynamics linearization; \cref{eq:SDE} can be linearized around a reference trajectory and control $\{\bm{x}^*(t), \bm{u}^*(t) \}$ 
based on a first-order truncation of the Taylor expansion:
\begin{equation}
    \dd{\bm{x}} \approx
    \left[A(t) \bm{x} + B(t) \bm{u} + \bm{c}(t) \right]\dd{t} + G(t) \dd{\bm{w}}
\end{equation}
where
\begin{equation}
    A(t) = \eval{\pdv{\bm{f}}{\bm{x}}}_{\bm{x}^*, \bm{u}^*}, \quad 
    B(t) = \eval{\pdv{\bm{f}}{\bm{u}}}_{\bm{x}^*, \bm{u}^*}, \quad
    \bm{c}(t) = \bm{f}(\bm{x}^*, \bm{u}^*) - A \bm{x}^* - B \bm{u}^*, \quad
    G(t) = \bm{g}(\bm{x}^*, \bm{u}^*)
\end{equation}
With the zero-order hold control assumption, we can exactly discretize the linearized dynamics as 
\begin{equation} \label{eq:SDE-linear-discrete}
    \bm{x}_{k+1} = A_k \bm{x}_k + B_k \bm{u}_k + \bm{c}_k + G_k \bm{w}_k, \quad k = 0, \cdots, N-1
\end{equation}
where the discrete-time matrices can be evaluated by simultaneous numerical integration of the following ordinary differential equations from time 
$t_k$ to $t_{k+1}$ \cite{kayama_2022_low-thrust},\highlight{\cite[Ch. 4.9]{tapleyStatisticalOrbitDetermination2004}}:
\begin{subequations} \label{eq:discrete-matrices}
\begin{align}
    \dv{t} \bm{x}(t)&=\bm{f}\left(\bm{x}(t), \bm{u}_{k}^{*}\right), \quad &&\bm{x}\left(t_{k}\right)=\bm{x}_k^{*} \label{eq:state_EOM_LTV}\\
    \dv{t} \Phi_{A}\left(t, t_{k}\right)&=A(t) \Phi_{A}\left(t, t_{k}\right), \quad &&\Phi_{A}\left(t_{k}, t_{k}\right)=I_{n_{x}} \\
    \dv{t} \Phi_{B}\left(t, t_{k}\right)&=A(t) \Phi_{B}\left(t, t_{k}\right)+B(t), \quad &&\Phi_{B}\left(t_{k}, t_{k}\right)=0_{n_x\times n_u} \\
    \dv{t} \Phi_P(t, t_k) &= A(t)\Phi_P(t) + \Phi_P(t)A(t)^\top + G(t)G(t)^\top, \quad && \Phi_P(t_k, t_k) = 0_{n_x \times n_x} 
\end{align}
\end{subequations}
\highlight{where $\Phi_P$ represents the covariance matrix of the Brownian motion, with its differential equation 
derived similar to Eq. (4.9.35) of \cite{tapleyStatisticalOrbitDetermination2004}.}
Letting $\overline{\bm{\phi}}(t_{k+1})$ be the result of integrating \cref{eq:state_EOM_LTV} from $t_k$ to $t_{k+1}$, 
\begin{align}
    A_k &= \Phi_A(t_{k+1}, t_k), \quad
    B_k = \Phi_B(t_{k+1}, t_k), \quad
    \bm{c}_k = \overline{\bm{\phi}}\left(t_{k+1}\right)-A_k \bm{x}_{k}^{*} - B_k \bm{u}_{k}^{*} 
\end{align}
and $G_k$ is any matrix that satisfies $G_k G_k^\top = \Phi_P(t_{k+1}, t_k)$
\highlight{, obtained, for example, by Cholesky decomposition}.
Note that when the reference state and control are known at each time step,
the computation of the discrete matrices for each segment is independent and can be parallelized.

\subsection{Filtered Discrete Covariance Propagation}
Having derived the linearized covariance propagation model, we assume that the 
filtering process is a Kalman filter, since the Kalman filter is a 
minimum-variance estimator for linear systems with Gaussian noise \cite{tapleyStatisticalOrbitDetermination2004}.
The nonlinear observation process \cref{eq:OD} can also be linearized with respect to a reference state, as
\begin{equation}
    \bm{y}_k = C_k \bm{x}_k + D_k \bm{v}_k + \bm{c}_{\mathrm{obs},k}
\end{equation}
where
\begin{equation}
    C_k = \eval{\pdv{\bm{f}_{\mathrm{obs}}}{\bm{x}}}_{\bm{x}_k^*}, \quad 
    D_k = \bm{g}_{\mathrm{obs}}(\bm{x}_k^*), \quad 
    \bm{c}_{\mathrm{obs}} = \bm{f}_{\mathrm{obs}}(\bm{x}_k^*) - C_k \bm{x}_k^*
\end{equation}

\begin{figure}[h] 
    \centering
    \includegraphics[width=0.6\textwidth]{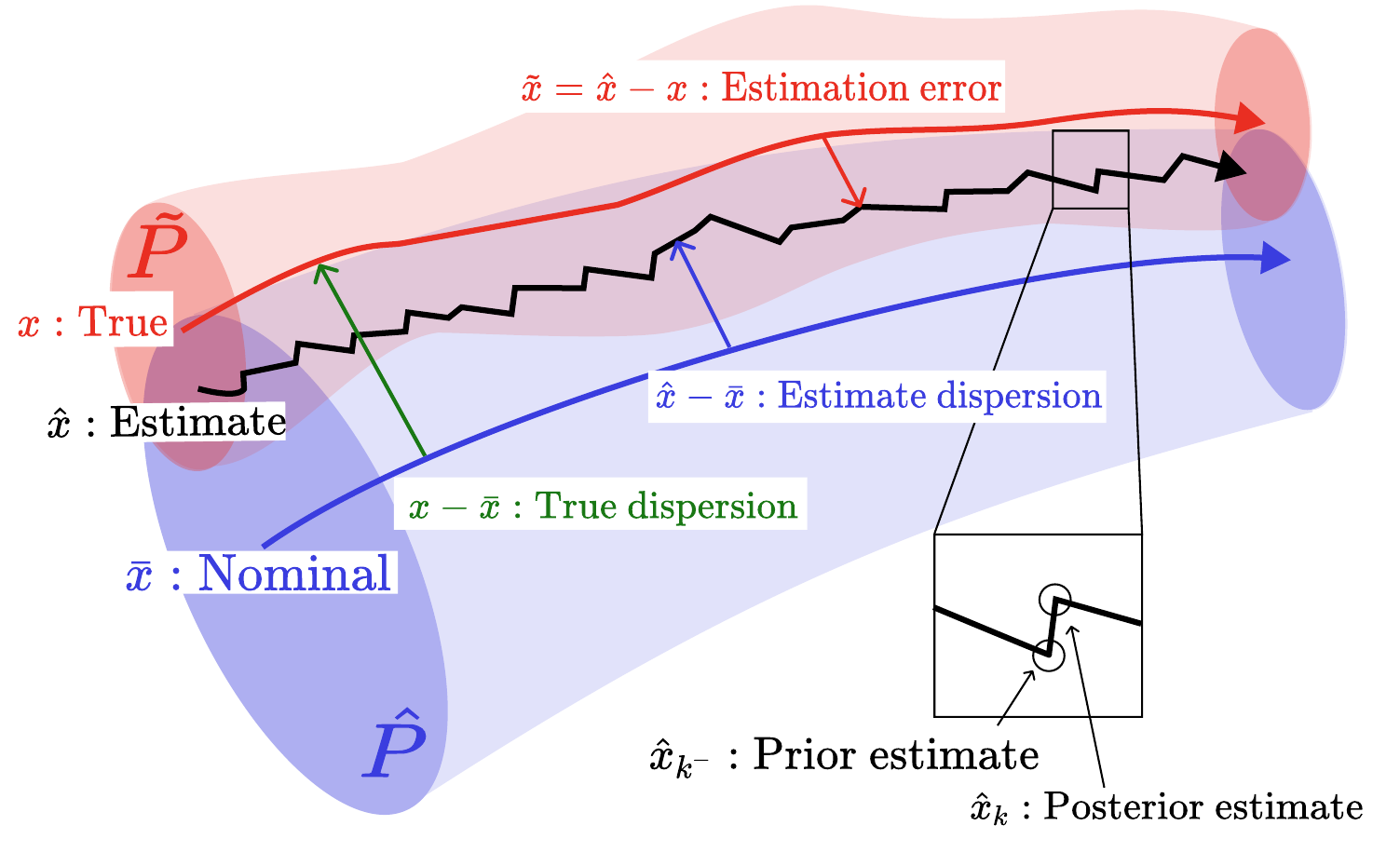}
    \caption{\highlight{Schematic of uncertainty-related terms}}
    \label{fig:schematic}
\end{figure}

\highlight{
We summarize necessary definitions for the derivations in this subsection, along with a visualization in \cref{fig:schematic}: 
$\bm{x}_k$ is the true state with mean $\overline{\bm{x}}_k := \E{\bm{x}_k}$. 
$\overline{\bm{x}}_k$ will also be referred to as the nominal state. 
$\hat{\bm{x}}_k$ is the state estimate; 
$\hat{\bm{x}}_{k^-}$ is the prior state estimate;
$\tilde{\bm{x}}_k := \bm{x}_k - \hat{\bm{x}}_k$ is the estimation error; 
$\tilde{\bm{x}}_{k^-} := \bm{x}_k - \hat{\bm{x}}_{k^-}$ is the prior estimation error.
Their covariance matrices are defined as
\begin{equation}
\begin{aligned}
    P_k &:= \E{(\bm{x}_k - \overline{\bm{x}}_k)(\bm{x}_k - \overline{\bm{x}}_k)^\top}, &&\\
    \hat{P}_k & := \E{ (\hat{\bm{x}}_k - \overline{\bm{x}}_k)(\hat{\bm{x}}_k - \overline{\bm{x}}_k)^\top},   
     &&\hat{P}_{k^-} := \E{ (\hat{\bm{x}}_{k^-} - \overline{\bm{x}}_k)(\hat{\bm{x}}_{k^-} - \overline{\bm{x}}_k)^\top} \\
    \tilde{P}_k &:= \E{ \tilde{\bm{x}}_k \tilde{\bm{x}}_k^\top}, 
    &&\tilde{P}_{k^-} := \E{ \tilde{\bm{x}}_{k^-} (\tilde{\bm{x}}_{k^-})^\top} 
\end{aligned}
\end{equation}
}
Between each time step, the Kalman filter updates the state estimate as
\begin{subequations} \label{eq:Kalman-updates}
\begin{align}
    &\hat{\bm{x}}_{k^-} = A_{k-1}\hat{\bm{x}}_{k-1} + B_{k-1} \bm{u}_{k-1} + \bm{c}_{k-1} && \triangleright \text{Time update}\\
    &\hat{\bm{x}}_k = \hat{\bm{x}}_k^- + L_k \tilde{\bm{y}}_k^- && \triangleright \text{Measurement update}
\end{align}
\end{subequations}
where the innovation process $\tilde{\bm{y}}_k^-$ and the Kalman gain $L_k$ are each defined as
\begin{align}
    \tilde{\bm{y}}_k^- &= \bm{y}_k - (\bm{c}_{\mathrm{obs},k} + C_k \hat{\bm{x}}_{k^-}) \\
    L_k &= \tilde{P}_{k^-} C_k^\top (C_k \tilde{P}_{k^-} C_k ^\top + D_k D_k^\top)^{-1}
\end{align}
The Kalman gain can be analytically calculated \textit{pre-optimization} \highlight{by recursively applying the following equations with 
$\tilde{P}_{0^-}$ given} \cite{tapleyStatisticalOrbitDetermination2004}:
\begin{subequations}\label{eq:update-covariances}
\begin{align}
    &\tilde{P}_{k^-} = A_{k-1} \tilde{P}_{k-1}A_{k-1}^\top + G_{\mathrm{exe}, k-1} G_{\mathrm{exe}, k-1}^\top + G_{k-1} G_{k-1}^\top, \quad 
    \mathhl{G_{\mathrm{exe}, k} = B_k G_{\mathrm{exe}}(\overline{\bm{u}}_k)} \\ 
    &\tilde{P}_k = (I-L_k C_k) \tilde{P}_{k^-}  (I-L_k C_k)^\top + L_k D_k D_k^\top L_k^\top
\end{align}
\end{subequations}
 The innovation process is Gaussian-distributed with zero mean and covariance:
\begin{equation}
    P_{\tilde{\bm{y}}_{k^-}} := \mathhl{ \E{ \tilde{\bm{y}}_{k^-} (\tilde{\bm{y}}_{k^-})^\top} }
    = C_k \tilde{P}_{k^-} C_k^\top + D_k D_k^\top
\end{equation}
Combining \cref{eq:Kalman-updates}, we get the discrete dynamics of the state estimate:
\begin{equation}
    \hat{\bm{x}}_{k+1} = A_k \hat{\bm{x}}_k + B_k \bm{u}_k + \bm{c}_k + L_{k+1}\tilde{\bm{y}}_{(k+1)^-}, \quad \hat{\bm{x}}_0 = L_0 \tilde{\bm{y}}_{0^-}
\end{equation}
Substituting the affine TCM policy \cref{eq:feedback} into the belief space dynamics, we get
\begin{equation}
    \hat{\bm{x}}_{k+1} = (A_k + B_k K_k) \hat{\bm{x}}_k \mathhl{- K_k \overline{\bm{x}}_k} + B_k \overline{\bm{u}}_k  + \bm{c}_k + L_{k+1}\tilde{\bm{y}}_{(k+1)^-}
\end{equation}
and its mean and covariance dynamics can be derived as
\begin{subequations} \label{eq:filtered-dynamics}
\begin{align}
&\overline{\bm{x}}_{k+1} = A_k \overline{\bm{x}}_k + B_k \overline{\bm{u}}_k + \bm{c}_k , &&\overline{\bm{x}}_0 = \bm{\mu}_0 \label{eq:filtered-dynamics-mean}\\ 
& \hat{P}_{k+1} = (A_k + B_k K_k) \hat{P}_k (A_k + B_k K_k)^\top + \hat{Q}_{k+1}, 
&& \hat{P}_0 = \hat{P}_{0^-} + \hat{Q}_0 \label{eq:filtered-dynamics-covariance}
\end{align}
\end{subequations}
where $\hat{Q}_k := L_{k} P_{\tilde{\bm{y}}_{k^-}} L_{k}^\top$ is defined for brevity.
\cite{ridderhof_chance-constrained_2020} shows that the state covariance can be 
obtained by simply adding the covariance of the state estimate and the estimation error:
\begin{equation}
    P_k = \tilde{P}_k + \hat{P}_k
\end{equation}
The developments in this section allow us to write \cref{prob:main} as follows. Note that 
while the linearized dynamics and the Kalman filter are used to approximate the covariance 
propagation, the mean of the state estimate is propagated through the nonlinear dynamics, 
as is done for the EKF. For the problem to be well-defined as an optimization problem,
and adhering to the spirit of the EKF, the linearized system matrices $A_k, B_k, \bm{c}_k, 
\hat{Q}_k$ should be written as functions of the \highlight{mean state} $\overline{\bm{x}}_k$.
For now, we assume that the linearization is performed around a given reference trajectory.

\begin{problem}{Discretized Chance-Constrained $\Delta V_{99}$ Minimization with Belief-Space Mean/Covariance Propagation} 
    \label{prob:main-discrete}
\begin{subequations}
    \begin{align}
    \minimize_{\mathhl{\overline{\bm{u}}_k, \overline{\bm{x}}_k, K_k, \hat{P}_k}} \quad & \mathcal{Q}\left(\sum\nolimits_{k=0}^{N-1} \norm{\bm{u}_k}; 0.99\right) \\
    \suchthat \quad
    & \overline{\bm{x}}_{k+1} = \overline{\bm{x}}_{k} + \int_{t_k}^{t_{k+1}}\bm{f}(\overline{\bm{x}}, \overline{\bm{u}}_k, t)\ \diff t  \label{eq:nonlinear_mean_eom} && \highlight{\triangleright \text{Discrete nonlinear mean dynamics}  }\\
    & \hat{P}_{k+1} = (A_k + B_k K_k) \hat{P}_k (A_k + B_k K_k)^\top + \hat{Q}_{k+1} \label{eq:linear_estimate_covariance_eom} && \highlight{\triangleright \text{Linear covariance dynamics}}\\
    & \P{\norm{\bm{u}_k} \leq u_{\max}} \geq 1 - \epsilon_{u} && \highlight{\triangleright \text{Control chance constraint}}\\
    & \hat{P}_k + \tilde{P}_k \preceq P_{\max} && \highlight{\triangleright \text{Covariance bound}}\\
    & \overline{\bm{x}}_N = \bm{\mu}_f, \quad \hat{P}_N + \tilde{P}_N \preceq P_f && \highlight{\triangleright \text{Final mean \& covariance constraints}}
    \end{align}
\end{subequations}
\end{problem}
\section{Solution via Convexification and Sequential Convex Programming} \label{sec:solution_via_SCP}

In this section, we present a numerical solution method to \cref{prob:main-discrete}.
The cost function and chance constraints are reformulated as deterministic second-order cone expressions, 
and the covariance dynamics\highlight{, originally nonlinear due to the coupling of the unknown feedback gains and covariance,} 
are reformulated as affine equality constraints \highlight{and linear matrix inequalities} via
lossless convexification. This allows us to solve the problem via sequential convex programming.

\subsection{Deterministic Convex Formulation}

\subsubsection{Objective Function Upper Bound}
To reformulate the optimization problem as a tractable deterministic
optimization, we first provide a deterministic upper bound on the objective function.
Under the linear covariance propagation and Gaussian initial state estimate distribution, each
control input is Gaussian-distributed and its mean and covariance 
are given by 
\begin{equation}
    \E{\bm{u}_k} = \overline{\bm{u}}_k, \quad \cov{\bm{u}_k} = K_k \hat{P}_k K_k^\top
\end{equation}
Hence, our objective function is a sum of the Euclidean norms 
of Gaussian-distributed random variables. To the authors' knowledge,
there is no closed-form expression (in terms of the variable's mean and covariance)
for the quantile function of the Euclidean norm of a nonzero-mean Gaussian random variable.
However, we can provide a deterministic upper bound; using the 
triangle inequality, \cite{oguri_chance-constrained_2024} provides the following bound:
\begin{equation} \
    \mathcal{Q}(\norm{\bm{u}_k}; p) \leq 
    \norm{\E{\bm{u}_k}} + \sqrt{\mathcal{Q}_{\chi_{n_u}^2}(p)} \sqrt{\lambda_{\max} (\cov{\bm{u}_k})} := \tilde{\mathcal{Q}}(\norm{\bm{u}_k}; p)
\end{equation}
where $\mathcal{Q}_{\chi_{n_u}^2}(\cdot)$ is the inverse cumulative distribution function for a chi-squared distribution with $n_u$ degrees of freedom. 
In later simulations, we see that this bound is quite tight for the problem settings we consider. 
Notice that the objective function is now separated into a sum of 
a term involving the mean dynamics and a term involving the covariance dynamics, which we can
write as $J_\mu$ and $J_{\Sigma}$.

\subsubsection{Case with No Chance Constraints}
In the absence of chance constraints, the mean steering problem and covariance steering problem are decoupled \cite{liu_optimal_2023}. Hence, we can solve two separate optimization problems:
\begin{problem}[Mean steering problem] \label{pr:mean-steering}
Minimize $J_\mu$ subject to the filtered mean dynamics \cref{eq:filtered-dynamics-mean} with mean terminal constraint of \cref{eq:final-distributional}.
\end{problem}
\begin{problem}[Covariance steering problem]\label{pr:covariance}
Minimize $J_\Sigma$ subject to the filtered covariance dynamics \cref{eq:filtered-dynamics-covariance} with covariance terminal constraint of \cref{eq:final-distributional}.
\end{problem}
While the mean steering problem is a traditional nonlinear trajectory optimization problem for which various solution methods exist, 
the covariance steering problem involves an equality \highlight{with bilinear ($K_k \hat{P}_k$) and cubic ($K_k \hat{P}_k K_k^\top$) expressions of the variables.}
However, \cite{liu_optimal_2023,pilipovsky_computationally_2023} show that we can equivalently formulate the problem as a convex optimization
for a certain class of objective functions, by utilizing \textit{lossless convexification}, i.e.
the solution obtained from solving a convex optimization problem gives the solution for the 
original nonconvex problem. The procedure for lossless convexfication is as follows:
1. replace $K_k$ by a set of new decision variables $U_k \in \R^{n_u\times n_x}$ such that
$ U_k := K_k \hat{P}_k $, 
2. introduce additional decision variables $Y_k \in \R^{n_u \times n_u}$, 
3. rewrite the covariance equations with $U_k$, replacing the expression $U_k {\hat{P}_k}^{-1} U_k^\top$ that appears, with $Y_k$,
4. additionally impose the linear matrix inequality (LMI) constraint
\begin{equation}
    \begin{bmatrix}
        \hat{P}_k & U_k^\top \\ U_k & Y_k
    \end{bmatrix} \succeq 0
\end{equation}
Then, the covariance steering problem is equivalently formulated as
\begin{subequations}
    \begin{align}
        \minimize_{Y_k, U_k, \hat{P}_k} \quad &J_{\Sigma} \\
        \suchthat \ 
        & \hat{P}_{k+1}  = A_k \hat{P}_k A_k^\top + B_k U_k A_k^\top + A_k U_k^\top B_k^\top + B_k Y_k B_k^\top + \hat{Q}_{k+1} \\ 
        &\begin{bmatrix}
            \hat{P}_k & U_k^\top \\ U_k & Y_k
        \end{bmatrix} \succeq 0 \\
        & \hat{P}_0 = \hat{P}_{0^-} + \hat{Q}_0 , \quad P_f - \hat{P}_N - \tilde{P}_N \succeq 0
    \end{align}
\end{subequations}
The feedback gains can be retrieved post-optimization by $K_k = U_k \hat{P}_k^{-1}$, and 
$Y_k$ represents the control covariance matrix at optimality. \cite{liu_optimal_2023,pilipovsky_computationally_2023} prove the validity of the lossless convexification 
for objectives with expectation-of-quadratic cost; essentially, the proof for lossless convexification relies on the objective function satisfying
\begin{equation} \label{eq:condition_for_lossless}
    \pdv{J_{\Sigma}}{Y_k} \succ 0
\end{equation}
which is true for the expectation-of-quadratic case. Unfortunately, for quantile 2-norm minimization, since \cite{petersenMatrixCookbook2012}
\begin{equation}
    \pdv{\sqrt{\lambda_{\max}(Y_k)}} {Y_k} = \frac{1}{2 \sqrt{\lambda_{\max}(Y_k)}} \pdv{\lambda_{\max}(Y_k)}{Y_k} = \frac{1}{2 \sqrt{\lambda_{\max}(Y_k)}} \bm{v}_{\max} \bm{v}_{\max}^\top
\end{equation}
where $\bm{v}_{\max}$ is the eigenvector corresponding to the maximum eigenvalue $\lambda_{\max}(Y_k)$,
this cost function only satisfies $\pdv{\sqrt{\lambda_{\max}(Y_k))}} {Y_k} \succeq 0$ (not strict linear matrix inequality), 
hence we cannot guarantee lossless convexification.  
Thus, we add an additional term in the objective:
\begin{equation}
    J_{\Sigma} := \sum_k \sqrt{\mathcal{Q}_{\chi_{n_u}^2}(p)} \sqrt{\lambda_{\max}(Y_k)} + \epsilon_Y \tr(Y_k)
\end{equation}
where $\epsilon_Y > 0$ is a small scalar to be chosen by the user. Since $\partial \tr(Y_k) / \partial Y_k = I \succ 0$, 
adding this term guarantees that \cref{eq:condition_for_lossless} is satisfied, hence the convexification is lossless.
\highlight{Theoretically, $\epsilon_Y$ can be an arbitrarily small positive number, but in practice, 
we should not use a value that is too small (e.g. below $10^{-4}$) since the solver operates with finite precision and termination tolerance.
Conceptually speaking, minimizing the trace of $Y_k$ is equivalent to minimizing the volume of 
the input covariance ellipsoid, which is different from $\Delta V_{99}$ but similarly reduces the 
statistical input.}

\subsubsection{Chance-Constrained Problem}
When chance constraints are present, the mean and covariance are coupled in the constraints. Here, we focus on the 2-norm upper bound on the control input.
With the same derivation as for the objective function, we can formulate a conservative deterministic form as 
\begin{equation} \label{eq:2-norm-deterministic}
    \norm{\overline{\bm{u}}_k} + \sqrt{\mathcal{Q}_{\chi_{n_u}^2}(1 - \epsilon_{u})} \sqrt{\lambda_{\max} (Y_k)} \leq u_{\max}
\end{equation}
since we already introduced the decision variable $Y_k$ which will be equal to the covariance of the control at optimality, i.e.
$ Y_k = K_k \hat{P}_k K_k^\top $.
$\lambda_{\max}(\cdot)$ is a convex operator for symmetric matrices \cite{boydConvexOptimization2004}. 
However, the square root over the maximum eigenvalue operator makes this deterministic formulation nonconvex. 
In order to handle this nonconvexity in the framework of SCP, we decompose this constraint into two inequalities, 
one of which remains nonconvex. The nonconvexity is treated via slack variables and Augmented Lagrangian penalization
in the objective function, which is a method with proven guarantees on optimality.
We introduce a set of variables $\tau_{0:N-1} \geq 0$, such that
\begin{equation} \label{eq:Y_tau}
    \lambda_{\max}(Y_k) - \tau_k^2 \leq 0
\end{equation}
Then, \cref{eq:2-norm-deterministic} is implied by the second-order cone constraint:
\begin{equation}\label{eq:2-norm-convex}
     \norm{\overline{\bm{u}}_k} + \sqrt{\mathcal{Q}_{\chi_{n_u}^2}(p)} \cdot \tau_k \leq u_{\max} 
\end{equation}
Although \cref{eq:Y_tau} is a difference-of-convex form and thus nonconvex, we simply linearize the concave term around some reference for $\tau$, which we denote as ${\tau}^{\mathrm{ref}}$:
\begin{equation}
    \lambda_{\max}(Y_k) - (\tau_k^{\mathrm{ref}})^2 - 2 \tau_k^{\mathrm{ref}} (\tau_k - \tau_k^{\mathrm{ref}}) \leq \zeta_k
\end{equation}
where $\zeta_k \geq 0 $ is a slack variable, which is penalized in the objective function via the \texttt{SCvx*} framework
if it takes positive values.
\begin{remark}[Conservativeness of the control upper bound relaxation]
    \normalfont \cref{eq:Y_tau} is not guaranteed to be satisfied with equality, hence this is a potentially conservative approximation to \cref{eq:2-norm-deterministic} (which is already a conservative approximation of \cref{eq:chance-constraints-control}). 
    However, this method has performed most reliably and predictably compared to
    globally overapproximating the square root via its linearization \cite{rapakouliasDiscreteTimeOptimalCovariance2023a}
    or exactly imposing difference-of-convex forms without slack variables\cite{pilipovsky_computationally_2023}.
    The numerical simulations also show little conservativeness in the control profile.
\end{remark}

We finally show the subproblem that is solved at each iteration of the SCP algorithm.
The objective function is a second-order-cone expression; the constraints are 
either linear matrix inequalities, second-order cone constraints, or affine equalities, since
$\infty$-norm constraints can be written as affine inequalities via auxiliary variables, and 
$\lambda_{\max}(\cdot)$ inequalities can be written as LMIs. Hence, the subproblem is a 
semidefinite programming (SDP), which can be solved via off-the-shelf solvers such as MOSEK \cite{mosek_aps_mosek_2023}.
For simplicity, we do not show the reformulation of these constraints into standard SDP form, since
parsers such as YALMIP \cite{lofberg_yalmip_2004} automate this process; see also \cite{boydConvexOptimization2004}
for details.
The slack variables are penalized in the objective function via a function $J_{\mathrm{pen}}(\cdot)$\highlight{, defined 
later in \cref{eq:J_pen}}. 
For notational simplicity, define $\bm{z} := \{ \overline{\bm{x}}_{0:N}$, $\hat{P}_{0:N}$, $\overline{\bm{u}}_{0:N-1}$, $U_{0:N-1}$, $Y_{0:N-1}$, $\tau_{0:N-1} \}, \bm{\xi} := [\bm{\xi}_0^\top, \cdots, \bm{\xi}_{N-1}^\top]^\top, \bm{\zeta} := [\zeta_0, \cdots, \zeta_{N-1}]^\top$, and 
\begin{equation}
J_0 (\bm{z}) := \sum_{k=0}^{N-1} \norm{\overline{\bm{u}}_k} + \sqrt{\mathcal{Q}_{\chi_{n_u}^2}(0.99)} \cdot \tau_k  + \epsilon_Y \tr(Y_k)     
\end{equation}
\begin{problem}[SDP subproblem solved at each SCP Iteration]
\label{prob:subproblem}
\begin{subequations} \label{eq:subproblem}
\begin{align}
    &\minimize_{\bm{z}, \bm{\xi}, \bm{\zeta} \geq 0}
    \quad J_{\mathrm{aug}}(\bm{z}, \bm{\xi}, \bm{\zeta}) := J_0 (\bm{z}) + J_{\mathrm{pen}}(\bm{\xi}, \bm{\zeta}) \label{eq:objective_augmented}\\
    &\quad \suchthat \nonumber\\
    &\text{for } k = 0,1, \cdots, N-1: \nonumber \\
    &\quad \overline{\bm{x}}_{k+1} - (A_k \overline{\bm{x}}_k + B_k \overline{\bm{u}}_k + \bm{c}_k) = \mathhl{\bm{\xi}_{k}} \label{eq:mean-LTV} && \highlight{\triangleright \text{Linearized mean dynamics}}\\ 
    &\quad \hat{P}_{k+1} = A_k \hat{P}_k A_k^\top + B_k U_k A_k^\top + A_k U_k^\top B_k^\top+ B_k Y_k B_k^\top + \hat{Q}_{k+1} \label{eq:covariance_dynamics} && \highlight{\triangleright \text{Linear covariance dynamics}}\\ 
    &\quad \begin{bmatrix}
        \hat{P}_k & U_k^\top \\ U_k & Y_k
    \end{bmatrix} \succeq 0 \label{eq:LMI-lossless} && \highlight{\triangleright \text{LMI for lossless convexification}}\\
    &\quad \lambda_{\max}(Y_k) - (\tau_k^{\mathrm{ref}})^2 - 2 \tau_k^{\mathrm{ref}} (\tau_k - \tau_k^{\mathrm{ref}}) \leq \zeta_{k} && \highlight{\triangleright \text{Control chance constraint 1}}\\
    &\quad \norm{\overline{\bm{u}}_k} +\sqrt{\mathcal{Q}_{\chi_{n_u}^2}(1 - \epsilon_u)} \cdot \tau_k - u_{\max} \leq 0  \label{eq:2-norm-convex-again} && \highlight{\triangleright \text{Control chance constraint 2}}\\
    &\quad Y_k \in \Symmetric^{n_u},\ Y_k\succeq 0, \ \tau_k \geq 0 && \highlight{\triangleright \text{Control covariance-related definitions}}\\
    & \text{for } k = 0, 1, \cdots, N: \nonumber \\
    &\quad \hat{P}_k + \tilde{P}_k \preceq P_{\max} && \highlight{\triangleright \text{Covariance bounds}}\\
    &\quad \hat{P}_k \in \Symmetric^{n_x} && \highlight{\triangleright \text{State covariance definition}}\\
    & \overline{\bm{x}}_0 = \bm{\mu}_0, \quad \overline{\bm{x}}_N = \bm{\mu}_N , \quad  
    \hat{P}_0 = \hat{P}_{0^-} + \hat{Q}_0, \quad 
    P_f - \hat{P}_N - \tilde{P}_N \succeq 0 && \highlight{\triangleright \text{Boundary mean \& covariance}}\\
    &s_{\bm{\overline{x}}}\| \bm{\overline{X}} - \bm{\overline{X}}^{\mathrm{ref}}\|_{\infty} \leq \Delta_{\mathrm{TR}}, \quad
    s_{\bm{\overline{u}}}\| \bm{\overline{U}} - \bm{\overline{U}}^{\mathrm{ref}}\|_{\infty} \leq \Delta_{\mathrm{TR}}, \quad
    s_{\bm{\tau}}\| \bm{\tau} - \bm{\tau}^{\mathrm{ref}}\|_{\infty} \leq \Delta_{\mathrm{TR}} \label{eq:trust-region} && \highlight{\triangleright \text{Trust region constraints}}
\end{align}
\end{subequations}
\end{problem}
While the initial reference trajectory for obtaining the linearized dynamics can be obtained by solving the deterministic version of the problem through 
some trajectory optimization algorithm, an initial reference $\tau_{0:N-1}^{\mathrm{ref}}$ is also required to initialize the SCP algorithm. 
This can be obtained by solving the covariance steering problem with $\sum \tr(Y_k)$ cost, 
where the linear system matrices are obtained by linearization around the initial reference trajectory.
Then, we can set ${\tau}_k^{\mathrm{ref}} = \sqrt{\lambda_{\max} (Y_k)}$ for each $k$.

\cref{eq:trust-region} imposes trust region constraints to maintain the validity
of the linearization by restricting the variation from the linearization reference to be sufficiently small.
$\Delta_{\mathrm{TR}} $ is the trust region, which is updated according to the algorithm's framework. Since the scaling of variables are different, we use a different scaling for the trust region, in the form of
and $s_{(\cdot)}$ is the scaling value for each variable. 
\highlight{
$\bm{\overline{X}} := [\bm{\overline{x}}_0^\top, \bm{\overline{x}}_1^\top, \cdots, \bm{\overline{x}}_N^\top]^\top$, 
$\bm{\overline{U}} := [\bm{\overline{u}}_0^\top, \bm{\overline{u}}_1^\top, \cdots, \bm{\overline{u}}_{N-1}^\top]^\top$, and 
$\bm{\tau} := [\tau_0, \cdots, \tau_{N-1}]^\top$ are defined for notational convenience. }
The $\mathrm{ref}$ superscript denotes the reference value, which is the value of the variable at a previous iteration.
We choose $s_{\bm{x}} = 1$ and $s_{\bm{u}} = u_{\max}$. $s_{\tau}$ is chosen based on the initial solve, set as $s_\tau = \max_k \tau_k$.

\subsection{Sequential Convex Programming Scheme}
We are ready to present the algorithm for the SCP-based CCDVO framework. The pseudoalgorithm is shown in \cref{alg:SCP}, 
based on the \texttt{SCvx*} algorithm \cite{oguri_successive_2023}. 
$\bm{z}^{\mathrm{ref}}$ is the reference solution from a previous iteration.
\begin{algorithm}[tb]
\caption{Robust Trajectory Optimization under Uncertainty via \texttt{SCvx*}}
\label{alg:SCP}
\begin{algorithmic}[1]
    \Require initial reference trajectory and control $\{\overline{\bm{X}}^{\mathrm{ref}}, \overline{\bm{U}}^{\mathrm{ref}}\}$
    \State Solve covariance steering problem with trace objective to get initial reference values ${\tau}_{0:N-1}^{\mathrm{ref}}$
    \While {iterations don't exceed the maximum}
        \If{first iteration or the reference was updated in previous iteration}
            \State Compute \{$A_k, B_k, \bm{c}_k, \hat{Q}_k, \tilde{P}_k\}$ from $\{\overline{\bm{X}}^{\mathrm{ref}}, \overline{\bm{U}}^{\mathrm{ref}}\}$
        \EndIf
        \State $\bm{z} \gets $ solve \cref{prob:subproblem}
        \If{convergence criteria met} \Comment{Line 2 of \cite{oguri_successive_2023}}
            \State \Return $\bm{z}$
        \EndIf
        \If{acceptance conditions met}
            \Comment{\cref{eq:acceptance}}
            \label{line:acceptance}
            \State 
            $\bm{z}^{\mathrm{ref}} \gets \bm{z}$
            \Comment{Solution update}
            \State Update penalty weight $w$, Lagrange multipliers $\bm{\lambda}, \bm{\mu}$
            \Comment{Line 13-15, Algorithm 1 of \cite{oguri_successive_2023}}
        \EndIf
        \State Update trust region $\Delta_{\mathrm{TR}}$%
        \Comment{\cref{eq:TRupdate}}
    \EndWhile
\end{algorithmic}
\end{algorithm}
Although we use the algorithm almost unchanged from the reference, 
there are a few aspects that warrant some discussion, including small algorithmic modifications to extend the algorithm to be applicable to trajectory optimization under uncertainty.
These aspects are highlighted in the following. Readers are referred to the reference for further details on the algorithm. 

\subsubsection{Augmented Lagrangian Penalty Function}

One notable difference is the augmented penalty in the objective function; we use the form
\begin{equation} \label{eq:J_pen}
    J_{\mathrm{pen}}(\bm{\xi}, \bm{\zeta} ; w, \bm{\lambda}, \bm{\mu}) = \bm{\mu} \cdot \bm{\xi} + \frac{w}{2} \bm{\xi} \cdot \bm{\xi} + \sqrt{w} \norm{\bm{\xi}}_1
     + \bm{\lambda} \cdot \bm{\zeta} + \frac{w}{2} \max(0, \bm{\zeta}) \cdot \max(0, \bm{\zeta}) + \sqrt{w} \norm{\max(0, \bm{\zeta})}_1
\end{equation}
where $\bm{\mu}, \bm{\lambda}$ are the Augmented Lagrangian multipliers for the equality and inequality constraints, respectively, and $w$ is the penalty weight. 
Note that $\bm{\xi}$ and $\bm{\zeta}$ are not limited to our earlier definitions, but can be generalized to any equality/inequality slack variables.
In contrast to the original \texttt{SCvx*} penalty, the penalty function additionally penalizes the 1-norm of the constraint violation; as also noted in \cite{oguriLosslessControlConvexFormulation}, adding this additional term empirically results in convergence with fewer iterations.

\subsubsection{Step Acceptance and Trust Region Update}

A key step in \cref{alg:SCP} is determining \textit{solution acceptance} in line \ref{line:acceptance}.
Like in \cite{oguri_successive_2023,maoSuccessiveConvexificationNonconvex2016}, the condition for solution acceptance is based on the ratio of the nonlinear cost reduction $\Delta J$ to the approximated cost reduction $\Delta L$, as follows:
\begin{align}
\rho = \frac{\Delta J}{\Delta L} = 
\frac{J_{\mathrm{aug}}^{\mathrm{NL}}(\bm{z}^{\mathrm{ref}} ) - J_{\mathrm{aug}}^{\mathrm{NL}}(\bm{z}^*)}
{J_{\mathrm{aug}}^{\mathrm{NL}}(\bm{z}^{\mathrm{ref}}) - J_{\mathrm{aug}}(\bm{z}^*, \bm{\xi}^*, \bm{\zeta}^*)}
\label{eq:rho}
\end{align}
where $J_{\mathrm{aug}}$ is given in \cref{eq:objective_augmented} while $J_{\mathrm{aug}}^{\mathrm{NL}}$ represents the cost of the nonlinear ``original'' problem with augmented Lagrangian,
\begin{equation}
    J_{\mathrm{aug}}^{\mathrm{NL}}(\bm{z}) := J_0(\bm{z}) + J_{\mathrm{pen}}(\bm{g}(\bm{z}), \bm{h}(\bm{z}); w, \bm{\lambda}, \bm{\mu})
\end{equation}
with $\bm{g}, \bm{h}$ being the original nonconvex expressions for the equality/inequality constraints which are linearized in the subproblem, in this case, 
\begin{equation}
    \bm{g}(\bm{z}) = \begin{bmatrix}
    \overline{\bm{x}}_{1} - \left(\overline{\bm{x}}_{0} + \int_{t_0}^{t_{1}}\bm{f}(\overline{\bm{x}}, \overline{\bm{u}}_0, t)\ \diff t\right) \\
    \vdots \\
    \overline{\bm{x}}_{N} - \left(\overline{\bm{x}}_{N-1} + \int_{t_{N-1}}^{t_{N}}\bm{f}(\overline{\bm{x}}, \overline{\bm{u}}_{N-1}, t)\ \diff t\right) 
\end{bmatrix}, \quad
    \bm{h}(\bm{z}) = \begin{bmatrix}
    \lambda_{\max}(Y_0) - \tau_0^2 \\
    \vdots \\
    \lambda_{\max}(Y_{N-1}) - \tau_{N-1}^2
\end{bmatrix}
\end{equation}
Note the same values of $\mathhl{w}, \bm{\lambda}, \bm{\mu}$ are used when calculating all the $J_{\mathrm{aug}}^{\mathrm{NL}}$ and $J_{\mathrm{aug}}$ in \cref{eq:rho}.

The value of $\rho$ is then used to determine whether to accept the current iteration and how to update the trust region radius $\Delta_{\mathrm{TR}} $.
In this work, we employ a different scheme than the one proposed in \cite{oguri_successive_2023}.
For $1 \geq \eta_0 > \eta_1 > \eta_2 > 0$, the step acceptance criterion is defined as:
\begin{align}
\begin{cases}
\text{accept}& \text{if}\ \rho \in [1 - \eta_0, 1 + \eta_0]
\\
\text{reject}& \text{else}
\end{cases}
\label{eq:acceptance}
\end{align}
while the trust region is updated as follows:
\begin{align}
\Delta_{\mathrm{TR}} \gets
\begin{cases}
\mathhl{\min}{\{\alpha_2 \Delta_{\mathrm{TR}}, \Delta_{\mathrm{TR}, \mathhl{\max}}\}}	& \text{if}\ \rho \in [1 - \eta_2, 1 + \eta_2]
\\
\Delta_{\mathrm{TR}}	& \text{elseif}\ \rho \in [1 - \eta_1, 1 + \eta_1]
\\
\mathhl{\max}{\{\Delta_{\mathrm{TR}} / \alpha_1, \Delta_{\mathrm{TR}, \mathhl{\min}}\}} & \text{else}
\end{cases}
\label{eq:TRupdate}
\end{align}
where $\alpha_1>1$ and $\alpha_2>1$ determine the contracting and enlarging ratios of $\Delta_{\mathrm{TR}} $, respectively, and $0 < \Delta_{\mathrm{TR}, \min} < \Delta_{\mathrm{TR}, \max} $ are the lower and upper bounds of the trust region radius.

This reformulation of the acceptance criterion and trust region update scheme is motivated by the fact that the convex subproblem presented in \cref{prob:main-discrete}, and hence \cref{prob:subproblem}, is based on \textit{inexact linearization} of the state covariance propagation
(see \cref{remark:inexact-linearization}).
This inexactness affects the analysis in the proof of Lemma 3 of \cite{oguri_successive_2023}, allowing $\Delta L $ defined in \cref{eq:rho} to take negative values (Lemma 3 of \cite{oguri_successive_2023} proves $\Delta L \geq 0 $ assuming exact linearization).
This observation implies that, under inexact linearization, a solution to the convex subproblem may yield $\Delta L < 0$, and that a solution with $\Delta J < 0$ is certainly a valid step as long as $\Delta L \approx \Delta J $.
Thus, this work proposes the criterion given in \cref{eq:acceptance} to accept a solution whenever $\rho$ is close to unity, which happens when the nonlinear and linear cost reductions take similar values.
This algorithmic choice is in a similar spirit to a successful trajectory optimization algorithm, hybrid differential dynamic programming (HDDP) \cite{lantoineHybridDifferentialDynamic2012a}, which applies augmented Lagrangian framework to DDP.

Another implication of inexact linearization is that \cref{alg:SCP} does not fully inherit the convergence guarantee provided by the \texttt{SCvx*} algorithm.
The convergence guarantee of the original \texttt{SCvx*} algorithm relies on the non-negativity of $\Delta L$, which, again, is not guaranteed here.
Note that this inexact linearization issue is present in other formulations for SCP-based trajectory optimization under uncertainty 
(e.g., \cite{ridderhof_chance-constrained_2020,oguri_stochastic_2022,benedikter_convex_2022,kumagaiSequentialChanceConstrainedCovariance2024}), although none of these studies seem to have pointed out this fact.
Rigorously analyzing the convergence behavior of SCP algorithms under inexact linearization, especially those for trajectory optimization under uncertainty, is left for important future work.

\begin{remark}
\label{remark:inexact-linearization}
\normalfont The linearized filtered state covariance dynamics given by \cref{eq:linear_estimate_covariance_eom} do not provide exact linearization of the state covariance propagation of the original system given in \cref{eq:SDE}.
As apparent in \cref{eq:linear_estimate_covariance_eom}, propagation of $P_k (= \hat{P}_k + \tilde{P}_k ) $ depends on the evaluation of the matrices $A_k,B_k, \bm{c}_k, \hat{Q}_k $, etc, which are obtained by linearizing the system about the reference mean $\{\overline{\bm{X}}^{\mathrm{ref}}, \overline{\bm{U}}^{\mathrm{ref}}\} $.
On the other hand, the current formulation in \cref{eq:linear_estimate_covariance_eom} does not capture the effect of changes in the mean trajectory $\{\overline{\bm{X}}, \overline{\bm{U}}\} $ to those in $A_k,B_k, \bm{c}_k, \hat{Q}_k $, etc, resulting in inexact linearization.
To exactly linearize the problem with respect to variables, one needs to capture the effect of change in $\{\overline{\bm{X}}, \overline{\bm{U}}\} $, which however requires the second-order sensitivity, i.e., state transition tensor (STT).
Although using STTs would address the inexact linearization issue, it would require a significantly heavier computational effort, for likely marginal improvement of the overall SCP performance.
\end{remark}

\subsection{Scaling System Matrices and Variables to Resolve Numerical Issues}
An issue inherent in space trajectory optimization problems is that the dispersion of the variables is much smaller than the value of the mean. 
If we solve the above problem without scaling, off-the-shelf solvers struggle with numerical issues, 
resulting in either warnings by the solver and/or 
constraints being violated with unignorable numerical errors.
Thus, it is necessary to introduce scaling of the constraints relating to the covariance. 
Although it is difficult to entirely resolve bad scaling in all constraint coefficients and variables, and
most off-the-shelf solvers apply their own internal scaling scheme to mediate numerical issues, 
scaling strategies are still important in avoiding the aforementioned issues.

\subsubsection{Proposed Scaling Method}

The scaling strategy is quite simple: introduce a scalar $d > 1$ and 
substitute the variables that relate to the covariance $\hat{P}_k, U_k, Y_k$, with 
\begin{equation}
    \mathhl{\hat{P}_k \leftarrow d^2 \cdot \hat{P}_k, \quad U_k \leftarrow d^2 \cdot U_k, \quad Y_k \leftarrow d^2 \cdot Y_k}
\end{equation}
One can see that the covariance dynamics equation \cref{eq:covariance_dynamics} remains the same, but 
with the coefficients $\hat{Q}_{k}, P_{\max}, P_f, \tilde{P}_k$ replaced with $d^2 \cdot \hat{Q}_{k}, d^2 \cdot P_{\max}$, and so on.
The linear matrix inequality \cref{eq:LMI-lossless} remains the same since the LHS is multiplied by a positive scalar.
The result of this scaling results in a ill-scaled reformulation for the objective and chance constraints, the objective function without penalty function augmentation becomes
\begin{equation}
    \sum_{k=0}^{N-1} \norm{\overline{\bm{u}}_k} + \sqrt{\mathcal{Q}_{\chi_{n_u}^2}(0.99)} \cdot \tau_k / d + \epsilon_Y \tr(Y_k) / d^2
\end{equation}
Note that $\tau_k$ is replaced with $d \cdot \tau_k$, from \cref{eq:Y_tau}, so it is divided by $d$.
We can similarly rescale the chance constraint. Although the objective and chance constraints become
ill-scaled, this did not cause noticable numerical issues.
After solving the optimization, the variables in the original scaling can be recovered by reversing the scaling.
Another benefit of introducing the scaling is that the use of slack variables in the linearized constraints
will be similarly penalized, leading to faster reduction of the slack variable values.

\subsubsection{Comparison with a Different Scaling Method}
In our previous paper \cite{kumagaiSequentialChanceConstrainedCovariance2024}, we implemented another scaling scheme based on \cite{benedikter_convex_2022}. 
Choosing a scaling matrix $D \in \R^{n_x \times n_x}$, 
both the coefficients of the covariance dynamics and the variables were scaled as
\highlight{$A_k \leftarrow D A_k D^{-1}$, $B_k \leftarrow D B_k$, $\hat{Q}_{k} \leftarrow D \hat{Q}_k D^\top$,
$U_k \leftarrow U_k D$, $\hat{P}_k \leftarrow D \hat{P}_k D^\top$.}
However, we observed through further numerical experiments that this numerical scaling scheme
is more prone to numerical issues compared to the above method. 
Our intuition for this is that the previous scheme introduces a wide range of scaling in the covariance-related 
variables, as $U_k$ is scaled `up' by $D$ and $\hat{P}_k$ is scaled `up' by $D^2$, while $Y_k$ is not
scaled at all, due to the $D$ matrix canceling out in the losslessness-inducing LMI constraint.
In addition, the coefficients of the covariance dynamics also possess a wide range of scaling, due to 
the quadratic term $B_k Y_k B_k$ now being scaled `up' by a factor $D^2$ while the term $A_k \hat{P}_k A_k^\top$ retains the same 
factor of scaling since $A_k \leftarrow D A_k D^{-1}$ approximately keeps the coefficient scaling. 
This wide range of scaling in the constraints and variables is likely the cause of the numerical issues.
\highlight{
    On the other hand, we note that the scaling method proposed in this paper
	has the drawback of using a scalar scaling parameter. This can be an issue when
	different elements (e.g. position and velocity) of the state covariance matrix have different orders of magnitude. 
    Generalizing the scaling method to allow different scaling for each state is left for future work.
}

\begin{remark}[Comparisons between block Cholesky and full covariance formulations]
    \normalfont Compared to the \textit{block} \textit{Cholesky} type algorithms proposed in \cite{ridderhof_chance-constrained_2020,oguri_stochastic_2022}, \textit{full covariance}-type approaches, 
    including our method, require significantly less computational effort. The former approaches construct a large linear matrix inequality (LMI) with size on the order of $N n_x$-by-$N n_u$, 
    while the latter constructs $N$ LMI's on the order of $(n_x + n_u)$-by-$(n_x+n_u)$ (see \cref{eq:LMI-lossless}). Assuming that the main computational burden of solving the subproblem comes from the LMI's, we can say that the computational complexity is $\order{N^2 n_x n_u}$ vs $\order{N (n_x + n_u)}$. As the number of discretization nodes increases, \textit{full covariance} methods greatly outperform the \textit{block Cholesky} methods. See also \cite{pilipovsky_computationally_2023} for a comparison of solution times for the problem under linear dynamics (i.e. for the linearized subproblem). As noted later in the numerical results, our method for 200 discretization nodes requires ten minutes, while \textit{block Cholesky} approaches would require several hours. 
    
    Another advantage that \textit{full covariance} approaches allows is the capability to perform feedback directly on the newest state estimate. 
    While \cref{eq:feedback} seems to be a standard affine feedback policy, it has not been the norm for \textit{block Cholesky} methods, such as state history feedback \cite{ridderhof_chance-constrained_2020,oguri_stochastic_2022} which requires storing all state histories and larger feedback gain matrices, and feedback based on a non-physical stochastic process \cite{okamoto_optimal_2019,oguri_chance-constrained_2024}. 
\end{remark}

\begin{remark}[Difference between ours and a similar method]\label{remark:differences}
    \normalfont While \cite{benedikter_convex_2022} proposes a similar \textit{full covariance} approach, the lossless convexification proposed for the covariance propagation is lossy \cite{rapakoulias_comment_2023}. 
    Although the numerical examples in \cite{benedikter_convex_2022} seem lossless, when implementing the algorithm for this work, we often encountered very lossy covariance propagations.
\end{remark}
\section{Numerical Simulations} \label{sec:numerical_simulations}
\begin{table}[htbp]
    \fontsize{9}{9}\selectfont

    \caption{Uncertainty parameters for DRO--DRO (NRHO--Halo) transfer}
    \centering
    \begin{tabular}{llll} \toprule 
        Category & Error Type & Symbol & Value  \\\midrule
        \multirow{3}{*}{Initial state estimate} & Position std dev & $\tilde{\sigma}_{r_0^-}$ & 50 \unit{km} \\
                                                & Velocity std dev & $\tilde{\sigma}_{v_0^-}$ & 1 \unit{m/s} \\ 
                                                & Covariance matrix & $\tilde{P}_{0^-}$ &  $\blkdiag{\tilde{\sigma}_{r_0^-}^2 I_3, \tilde{\sigma}_{v_0^-}^2 I_3}$  \\ \midrule
        \multirow{3}{*}{Final state} & Position std dev & $\sigma_{r_f} $ & 20  \unit{km} \\
                                        & Velocity std dev & $\sigma_{v_f}$ & 0.1 \unit{m/s} \\ 
                                        & Covariance matrix & $P_f$ &  $\blkdiag{\sigma_{r_f}^2 I_3, \sigma_{v_f}^2 I_3}$  \\ \midrule
        \multirow{4}{*}{Gates model}  & Fixed magnitude err std dev             & $\sigma_1$ & \highlight{$10^{-3}$  \unit{mm/s^2} } \\
                                      & Proportional magnitude error std dev    & $\sigma_2$ & \highlight{1 \%}    \\
                                      & Fixed pointing error std dev            & $\sigma_3$ & \highlight{$10^{-3}$   \unit{mm/s^2} } \\
                                      & Proportional pointing error std dev     & $\sigma_4$ & \highlight{0.5 \unit{deg} }  \\ \midrule
         \multirow{2}{*}{Navigation} & Position std dev & $\sigma_r^{\mathrm{nav}}$ & 10 \unit{km} \\
                                    & Velocity std dev & $\sigma_v^{\mathrm{nav}}$ & 0.1 \unit{m/s} \\ \midrule
         Chance constraints & Allowed control violation probability & $\epsilon_u$ & 0.01  \\ \midrule
         \multirow{3}{*}{Maximum covariance} & Position std dev & $\sigma_{r,\max} $ & $\infty$ (500)  \unit{km} \\
         & Velocity std dev & $\sigma_{v,\max}$ & $\infty$ (3) \unit{m/s} \\ 
         & Covariance matrix & $P_{\max}$ &  $\blkdiag{\sigma_{r,\max}^2 I_3, \sigma_{v,\max}^2 I_3}$  \\ \midrule
         Unmodeled acceleration & Stochastic acceleration & $\sigma_{a}$ &   $ 1 \times 10^{-10} (10^{-11})$ \unit{km/s^{3/2}} \\
         \bottomrule 
    \end{tabular}
    \label{tab:parameters}    
\end{table}
This section demonstrates \cref{alg:SCP} on two transfers between periodic orbits in the Earth-Moon CR3BP. 
The values of uncertainties for each scenario are shown in \cref{tab:parameters}. 
The stochastic term in the SDE has a constant covariance matrix, with 
noise entering the state dynamics through change in velocity:
\begin{equation}\nonumber
    \bm{g}(\bm{x}, \bm{u}) = \left[0_{3 \times 3}, \ \sigma_a I_3 \right]^\top
\end{equation}
The orbit determination process is a full-state observation with Gaussian noise:
\begin{equation} \nonumber
    \bm{f}_{\mathrm{obs}}(\bm{x}) = \bm{x}, \quad \bm{g}_{\mathrm{obs}}(\bm{x}) = \blkdiag{ \sigma_r^{\mathrm{nav}} I_3, \sigma_v^{\mathrm{nav}} I_3}
\end{equation}
\begin{table}[ht]
    \centering
    \fontsize{9}{9}\selectfont
    \caption{\texttt{SCvx*} parameters}
    \begin{tabular}{ccccccccccc}\toprule
        Parameter & $\{\epsilon_{\mathrm{opt}}, \epsilon_{\mathrm{feas}}\}$ & $\{\eta_0, \eta_1, \eta_2\}$ & $\{\alpha_1, \alpha_2\}$ & $\beta$ & $\gamma$ & \{$\Delta_{\mathrm{TR},\min}, \Delta_{\mathrm{TR},\max}$\} & $\Delta_{\mathrm{TR}}^{(1)}$ & $w_\max$ & $w^{(1)}$ \\\midrule
        Value     & \{$10^{-4}, 10^{-6}$\} & $\{1, 0.85, 0.1\}$           & $\{2, 3\}$               & 2       & 0.9      & \{$10^{-6}, 1  $\}     & 0.3       &  $10^{8}$   & 1000 \\\bottomrule
    \end{tabular}
    \label{tab:scp_parameters}
\end{table}
The parameters for the \texttt{SCvx*} algorithm are shown in \cref{tab:scp_parameters}. 
The mass ratio is $\mu = 0.01215059$ for the CR3BP. 
Units for nondimensionalization are $\mathrm{LU} = 384748 $ \ \unit{km} and $\mathrm{TU} = 375700 \ \unit{s}$. 
$\epsilon_Y = 10^{-4}$ is used in \cref{eq:objective_augmented} and $d = 100$ for scaling.
We use the \texttt{ode78} propagator in MATLAB for all numerical integration, with absolute and relative tolerances set to $10^{-13}$. 
Convex optimization problems are solved with YALMIP \cite{lofberg_yalmip_2004} and MOSEK \cite{mosek_aps_mosek_2023}. All simulations are performed with MATLAB on a standard personal laptop.

\subsection{DRO--DRO Transfer}
\begin{table}[h]
    \fontsize{9}{9}\selectfont
    \caption{Initial conditions for DROs ($y = z = \dot{x} = \dot{z} = 0$) }
    \centering
    \begin{tabular}{cccccc}\toprule
        Orbit & $x$ (ND) & $\dot{y}$ (ND) & Period (ND) & Period (days) \\ \midrule
        DRO \#1 & 0.58041127991124 & 0.973651613293327 & 5.71743682447432 & $\approx$ 25.3\\
        DRO \#2  & 0.233114246213419 &  2.41810511614024 & 6.2574913469559279 & $\approx$ 27.7 \\ \bottomrule
    \end{tabular}
    \label{tab:DROs}
\end{table}
\begin{table}[h]
    \fontsize{9}{9}\selectfont
    \caption{Parameters for DRO-to-DRO transfer}
    \centering
    \begin{tabular}{ccc} \toprule
        Parameter & Value & Unit \\ \midrule
        Time-of-flight &  25 & days\\
        Maximum s/c acceleration & $0.5$ & \unit{mm/s^2} \\  
        Discretization nodes & 50 & NA \\
        \bottomrule
    \end{tabular}
    \label{tab:DRO}
\end{table}

\begin{figure}[h!]
    \centering
    \begin{subfigure}[c]{0.25\textwidth}
    \centering
    \includegraphics[width=\linewidth]{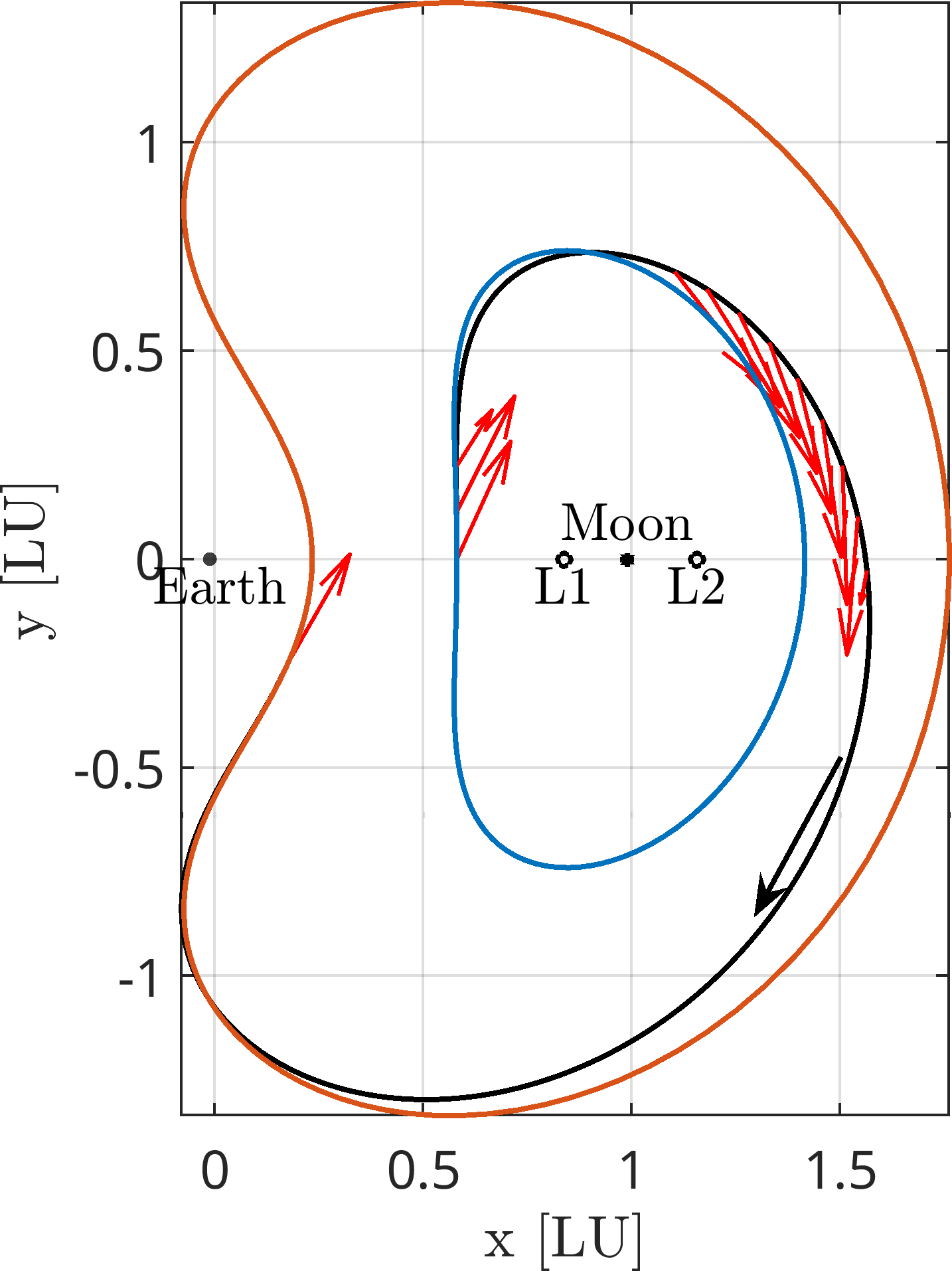}
    \caption{Trajectory}
    \label{fig:DROtoDRO-reference-traj}
    \end{subfigure}%
    \hspace{0.1\textwidth}
    \begin{subfigure}[c]{0.35\textwidth}
    \centering
        \includegraphics[width=\linewidth]{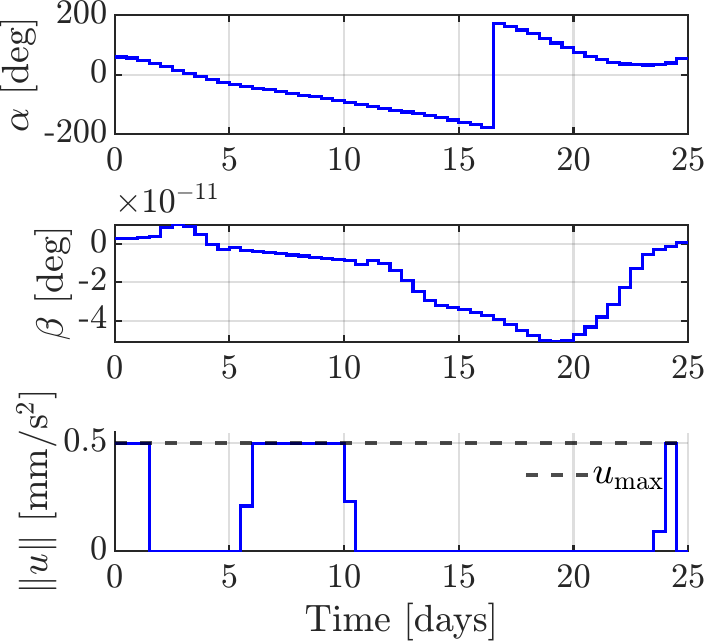}
    \caption{Control profile}
    \label{fig:DROtoDRO-reference-control}  
    \end{subfigure}%
    \caption{Reference fuel-optimal DRO-to-DRO transfer from DRO \#1 (inner) to DRO \#2 (outer)}
    \label{fig:DROtoDRO-reference}
\end{figure}

\begin{figure}[thbp]
    \centering
    \begin{subfigure}[c]{0.4\textwidth}
        \centering
        \includegraphics[width=\linewidth]{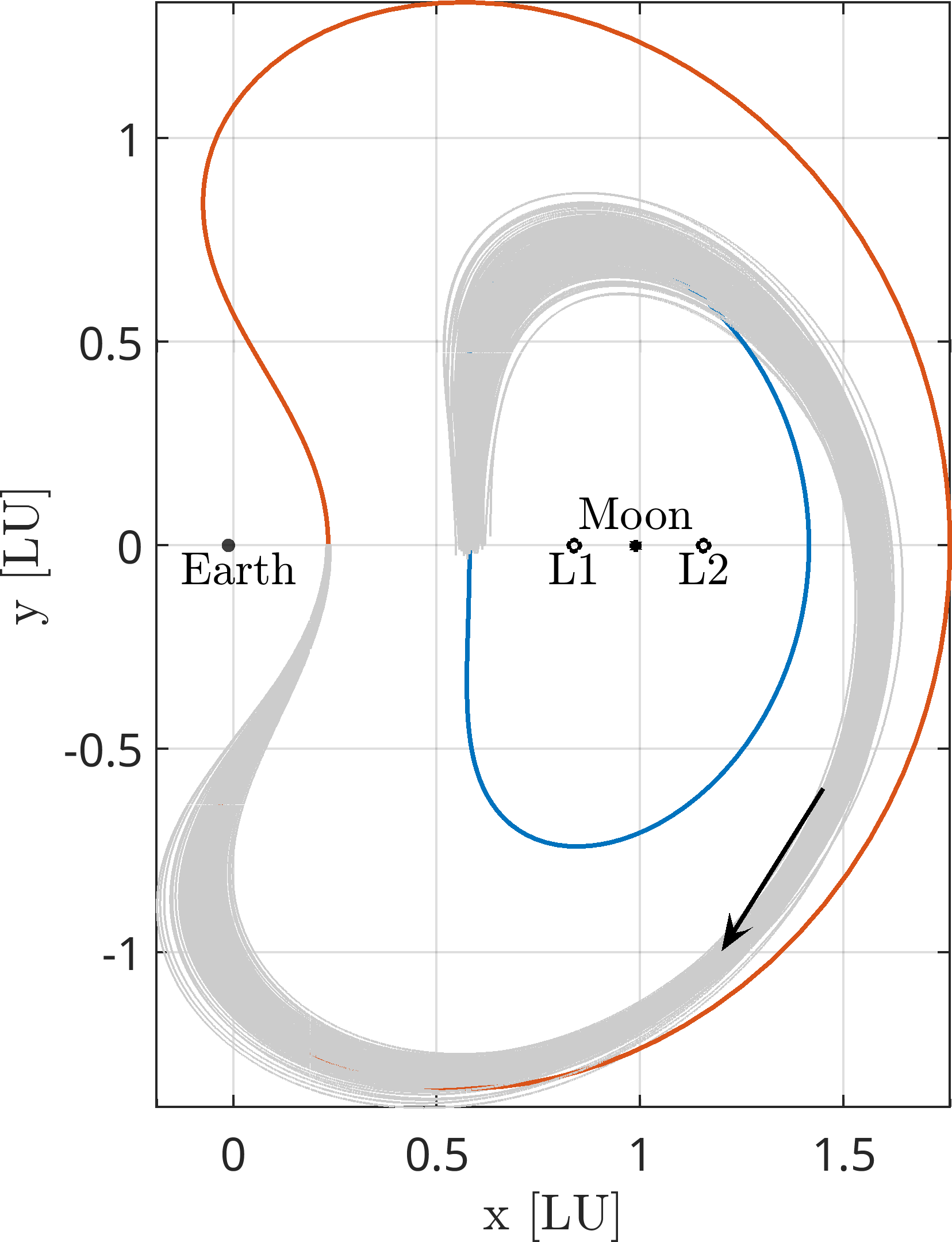}
        \caption{Trajectories (Deviations enlarged 50x)}
        \label{fig:DROtoDRO-robust-traj}
    \end{subfigure}\hfill%
    \begin{subfigure}[c]{0.5\textwidth}
        \centering
        \includegraphics[width=\linewidth]{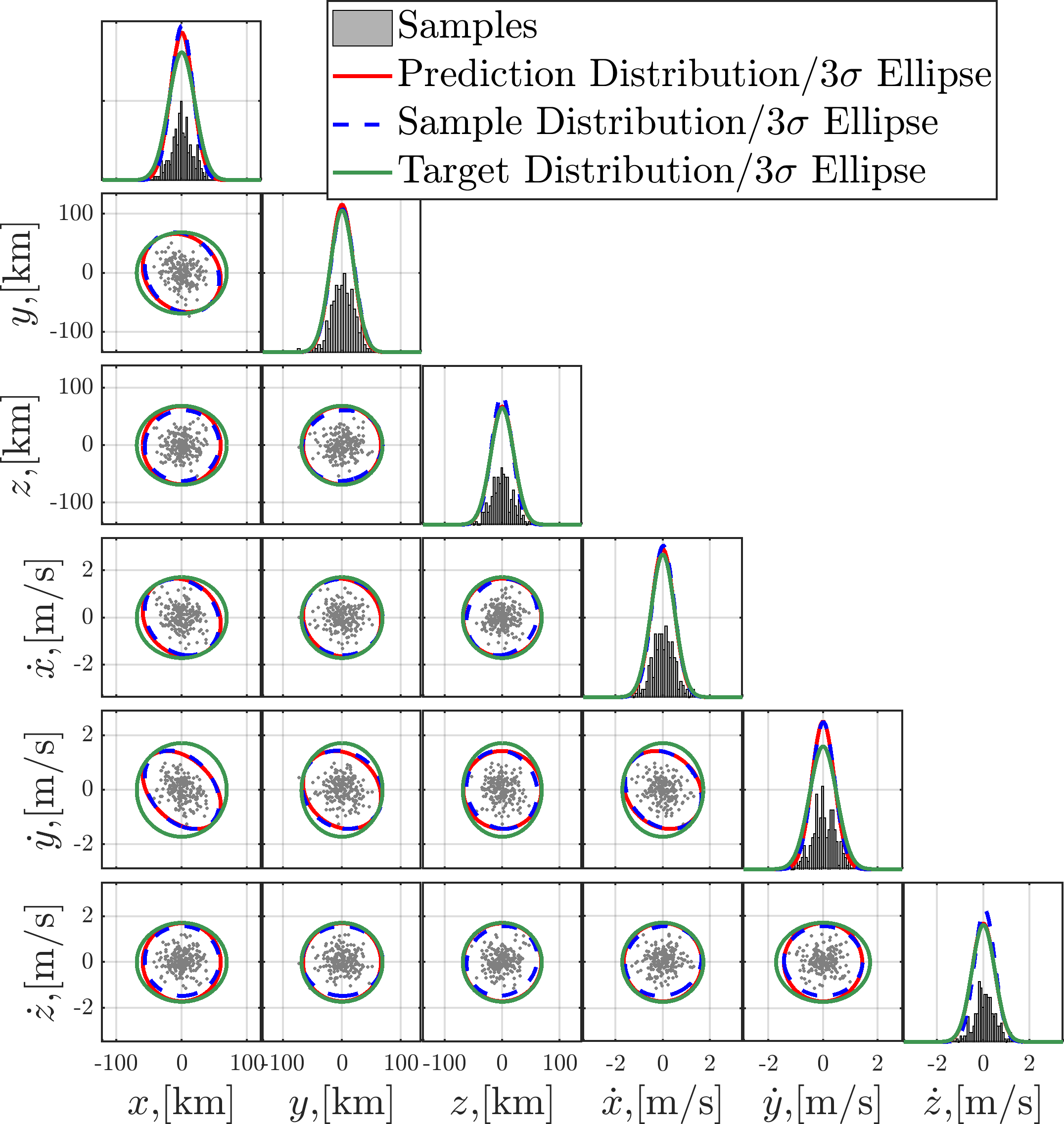}
        \caption{State dispersion at the final node}
        \label{fig:DRO-DRO-final-dispersion}
    \end{subfigure}
    \begin{subfigure}[c]{0.4\textwidth}
        \centering
        \includegraphics[width=\linewidth]{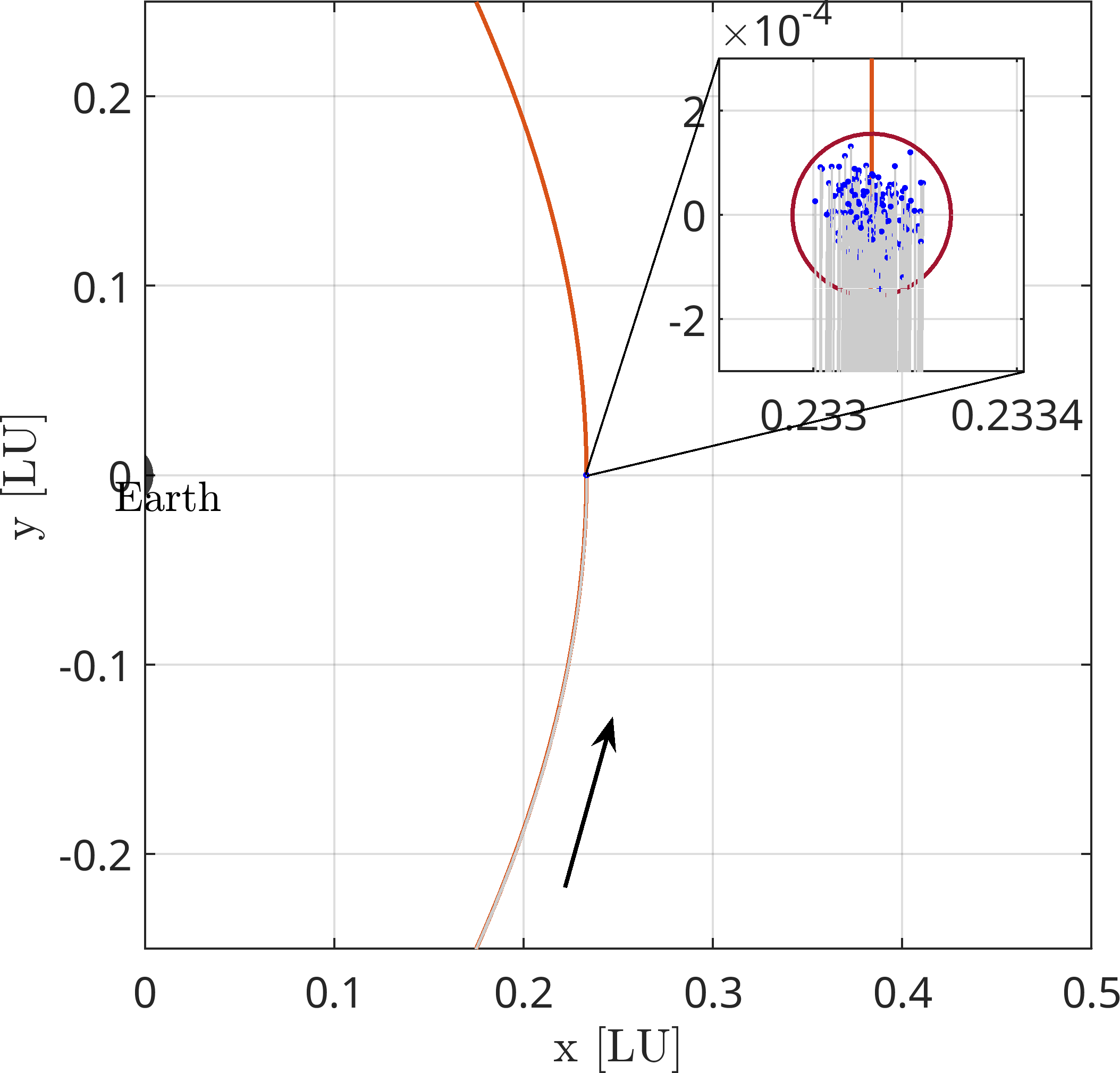}
        \caption{\highlight{In-plane dispersions at the final node (closed-loop)}}
        \label{fig:DROtoDRO-robust-traj-endpoint}
    \end{subfigure}\hfill%
    \begin{subfigure}[c]{0.4\textwidth}
        \centering
        \includegraphics[width=\linewidth]{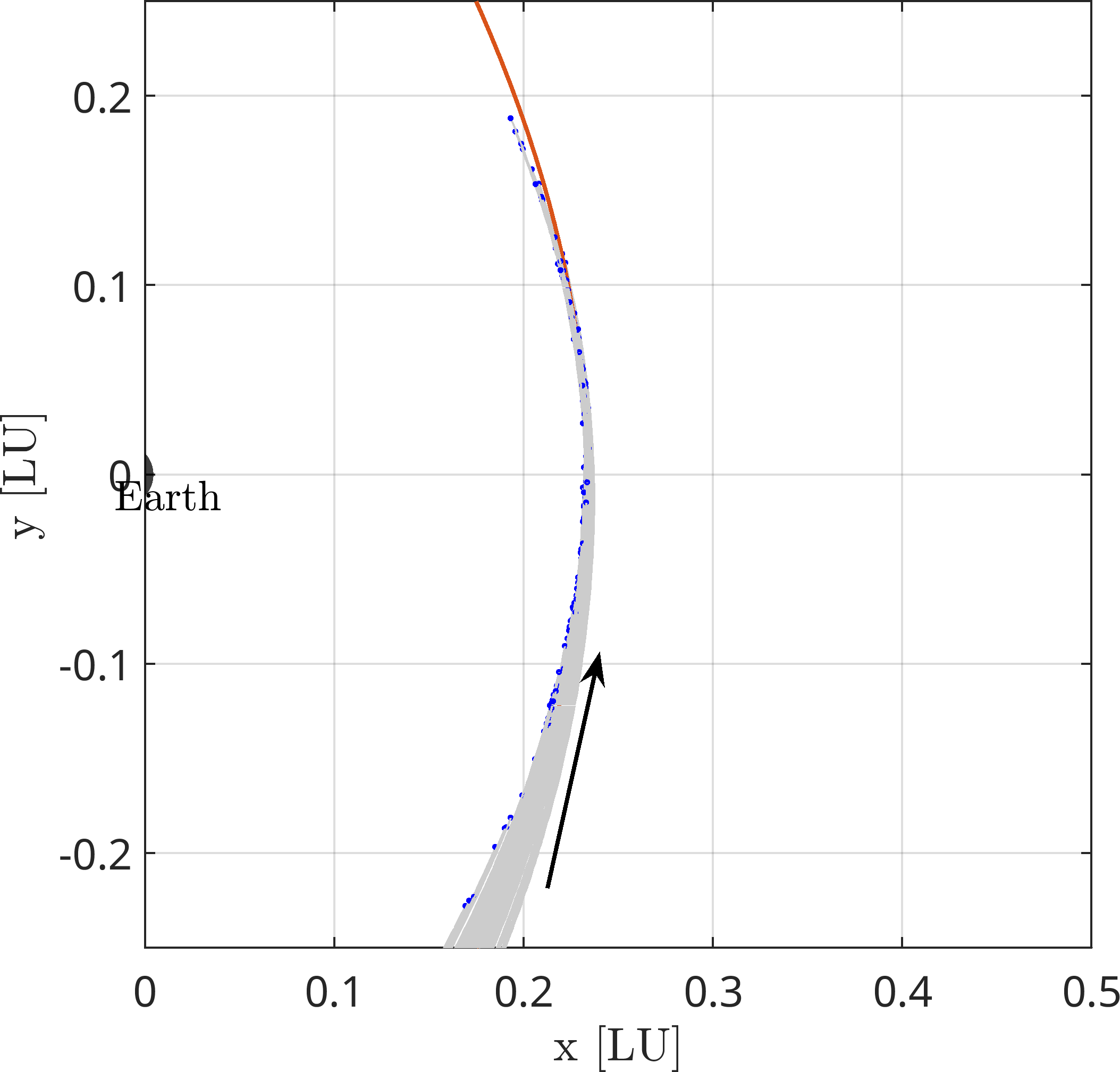}
        \caption{\highlight{In-plane dispersions at the final node (open-loop)}}
        \label{fig:DROtoDRO-robust-traj-openloop-endpoint}
    \end{subfigure}
    \caption{Monte Carlo simulation trajectories for DRO--DRO transfer.}
    \label{fig:DROtoDRO-trajectories}
\end{figure}

\begin{figure}[htbp]
    \begin{subfigure}[c]{0.5\textwidth}
        \includegraphics[width=\linewidth]{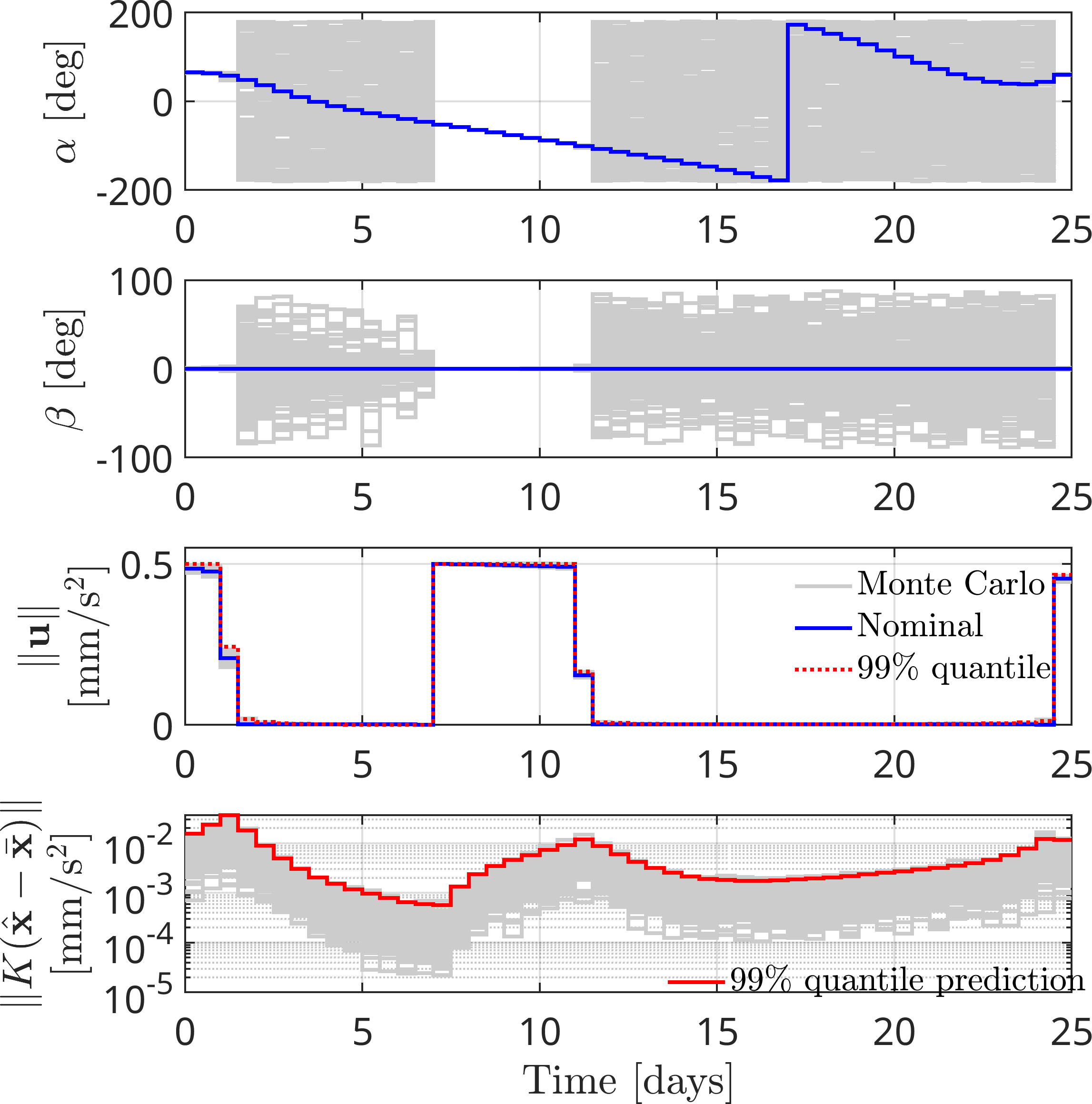}
        \caption{Control profile}
        \label{fig:DROtoDRO-control-profile}  
    \end{subfigure}\hfill%
    \begin{subfigure}[c]{0.45\textwidth}
        \includegraphics[width=\linewidth]{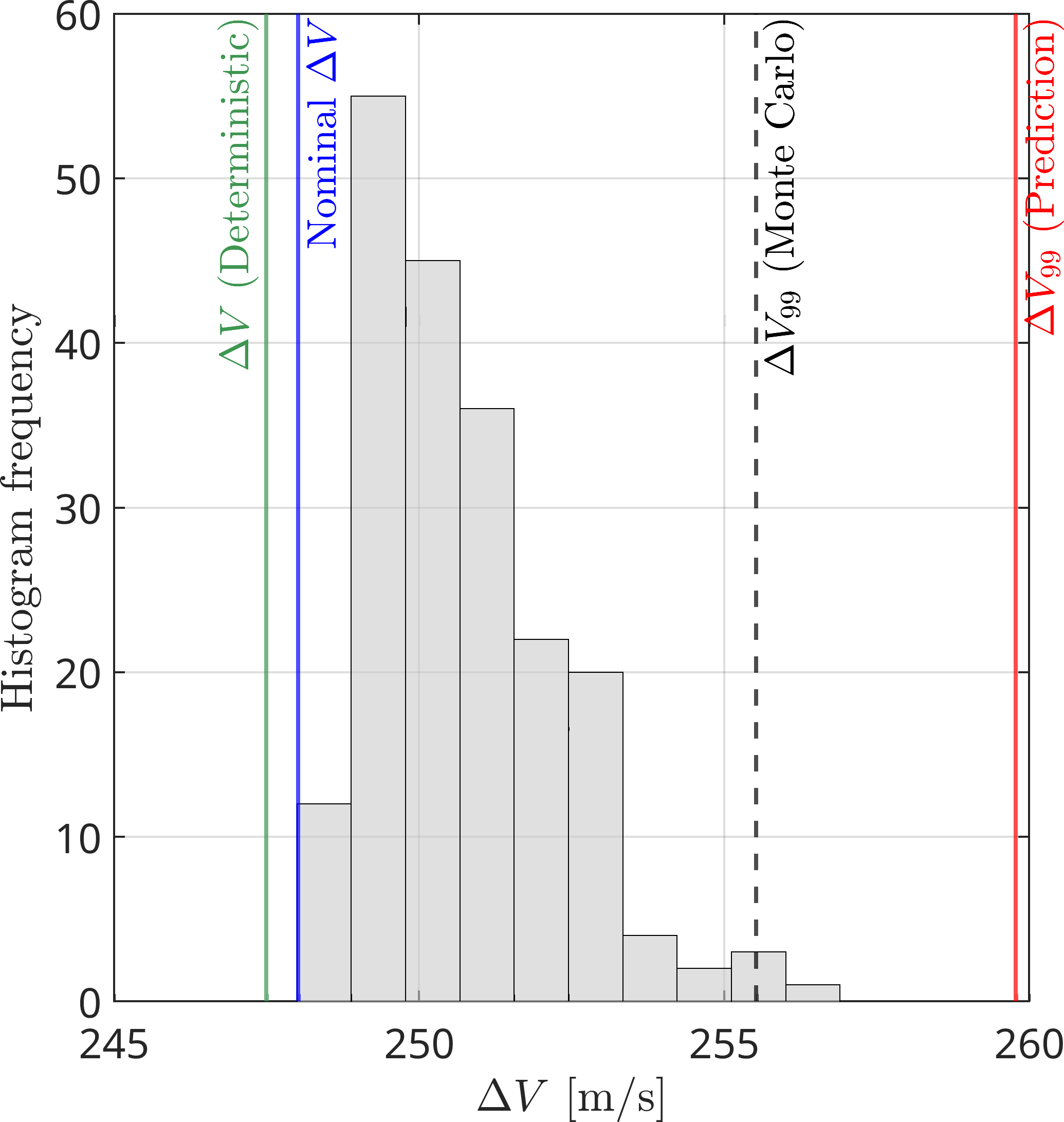}
        \caption{$\Delta V$ comparisons}
        \label{fig:DROtoDRO-DV}
    \end{subfigure}
    \caption{Monte Carlo simulation control statistics for DRO--DRO transfer.}
    \label{fig:DROtoDRO-control}
\end{figure}

As our first example, we consider a planar transfer between two Distant Retrograde Orbits (DROs), both around the L2 equilibrium. The initial conditions of the DROs are shown in 
\cref{tab:DROs}. Depending on one's working environment, differential corrections may be required to obtain a precisely closed orbit \cite{pavlak_trajectory_2013}. The problem parameters are summarized in \cref{tab:DRO}. The initial reference trajectory is obtained as follows. First, the minimum-energy trajectory with unbounded control input is obtained with the indirect method, using the MATLAB nonlinear equation solver \texttt{fsolve}. With the implementation of the analytical derivatives of the shooting function (such as outlined in \cite{arya_hyperbolic-tangent-based_2019}), the indirect method with minimum-energy cost can be found easily, with initial random costate guesses. Next, the trajectory obtained is used as an initial reference trajectory in the sequential convex programming-based nonconvex optimization algorithm, \texttt{SCvx*} \cite{oguri_successive_2023}. The generation of optimal reference trajectories is itself a vast topic, and we do not intend to claim that this is the best method; rather our focus is on `robustifying' the nominal trajectories. The reference trajectory and control are shown in \cref{fig:DROtoDRO-reference}.
The angles $\alpha, \beta$ are defined as $\alpha = \operatorname{arctan2} (u_y / u_x)$ and $\beta = \arcsin (u_z / \sqrt{u_x^2+u_y^2+u_z^2})$, with $\operatorname{arctan2}$ representing the four-quadrant inverse tangent function.
Next, we apply \cref{alg:SCP} to obtain optimal control policies to correct for operational and dynamics uncertainties. 
The algorithm converges in 6 iterations in 10 seconds.
Note that this is orders of magnitude faster than the \textit{block Cholesky} methods, which we have implemented in \cite{oguri_stochastic_2022} and can require tens of minutes for convergence; 
this is consistent with the observations made in the linear covariance steering literature \cite{pilipovsky_computationally_2023}. 
To verify our solution, we also perform a Monte Carlo simulation with 200 samples.
\cref{fig:DROtoDRO-trajectories,fig:DROtoDRO-control} show the statistics prediction based on the 
linear covariance model, as well as the Monte Carlo simulation results.
The trajectories are shown in \cref{fig:DROtoDRO-robust-traj}, 
and the control profiles in \cref{fig:DROtoDRO-control}.
In \cref{fig:DROtoDRO-robust-traj}, the sample trajectories stay bounded throughout the transfer
thanks to the trajectory corrections, and the dispersion shrinks rapidly towards the end of the transfer.
Figure \ref{fig:DRO-DRO-final-dispersion} shows the statistics at the final epoch for each combination of the state variable,
along with the distribution for the prediction ($P_N = \hat{P}_N + \tilde{P}_N$), samples, and target ($P_f$).
The plots on the diagonals show the dispersions of each state component in terms of a histogram.
We see that 1) the prediction distribution (red) accurately predicts the sample distribution (blue) and
2) the sample distribution (blue) 3-sigma ellipses are within the target ellipse (green) in all state combinations, 
indicating that the final covariance constraint is satisfied in the Monte Carlo samples.
\highlight{\cref{fig:DROtoDRO-robust-traj-endpoint,fig:DROtoDRO-robust-traj-openloop-endpoint} 
show the final state dispersion in blue, at the last node for the closed-loop and open-loop cases, respectively.
Open-loop means that the the feedback gain $K$ is set to zero. We see that, in contrast to the closed-loop case, 
the open-loop case shows large deviations along the final velocity direction.}

Figure \ref{fig:DROtoDRO-control-profile} shows the control profile for the robust trajectories.
The first two plots show the in-plane and out-of-plane control profiles for the 
nominal trajectory (blue) and the Monte carlo samples (grey). Comparing against 
the third plot, which shows the control acceleration magnitude, we see that during 
the thrusting arcs, the control is dominated by the feedforward term, resulting in 
small dispersion in the control angles. On the other hand, during the coasting arcs,
the TC term dominates, resulting in larger dispersion in the control angles.
The third plot also shows the 99\% bound, computed from the left-hand side of \cref{eq:2-norm-deterministic}.
The bound accurately bounds the Monte Carlo samples, suggesting that the 
linear covariance model is also successful in predicting the control distribution.
The fourth plot shows the control acceleration magnitude for the TC term only.
We see that the TC effort is largest \highlight{towards the end of the first and second thrusting arc,
which ``cleans up'' the effect of maneuver execution error that accumulates during this arc.
There is another peak in TC effort towards the end of the transfer which starts before the third short 
thrusting arc, most likely in order to satisfy the final covariance constraint.
} 
Figure \ref{fig:DROtoDRO-DV} shows the $\Delta V$ comparisons for the deterministic and robust trajectories.
The deterministic trajectory has a lower $\Delta V$ compared to that of the nominal robust trajectory.
We expect that this is due to the robust trajectory being optimized to also minimize the 
TC effort and the chance/distributional constraints. 
The $\Delta V_{99}$ prediction (red) successfully bounds that of the Monte Carlo samples (black broken line),
suggesting that the $\Delta V_{99}$ prediction with our model can provide upper bounds
for the Monte Carlo samples (without running the computationally expensive Monte Carlo simulation).
\highlight{
    \begin{remark}
        \normalfont
    We note that, as evident in the control profile, this work assumes a constant bound on the 
thrust \textit{acceleration}, implicitly neglecting the mass dynamics. However, when considering mass dynamics, 
monotonically decreasing mass implies a monotonically increasing available thrust acceleration.
This implies that the computed control policy will be \textit{feasible} when using the initial mass to 
compute the acceleration bound. Second, for a typical ion-engine spacecraft, the consumed mass is 
a few percent of the total mass, and thus the loss in \textit{optimality} is small.
    \end{remark}
}

\subsection{L2 NRHO to L1 Halo Transfer}
\begin{table}[ht]
    \fontsize{9}{9}\selectfont
    \caption{Initial conditions for L2 NRHO to L1 Halo transfer. ($y = \dot{x} = \dot{z} = 0$) }
    \centering
    \begin{tabular}{cccccc}\toprule
        Orbit & $x$ (ND) & $z$ (ND) & $\dot{y}$ (ND) & Period (days) \\ \midrule
        L2 NRHO & 1.018826173554963 & -0.179797844569828 & -0.096189089845127 &  $\approx$ 6.4 \\
        L1 Halo &0.823383959653906 & 0.010388134109586 &  0.128105259453086 &  $\approx$ 12.0 \\ \bottomrule
    \end{tabular}
    \label{tab:NRHO-Halo-ICs}
\end{table}%
\begin{table}[ht]
    \fontsize{9}{9}\selectfont
    \caption{Parameters for L2 NRHO to L1 Halo transfer}
    \centering
    \begin{tabular}{ccc} \toprule
        Parameter & Value & Unit \\ \midrule
        Time-of-flight &  57.4 & days\\
        Maximum spacecraft acceleration & $0.273$  & \unit{mm/s^2} \\  
        Number of discretization nodes & 200 & NA \\ 
        \bottomrule
    \end{tabular}
    \label{tab:NRHO}
\end{table}%
Next, we demonstrate a transfer from an L2 Near Rectilinear Halo Orbit (NRHO) to an L1 Halo orbit. 
The initial conditions for the periodic orbits are in \cref{tab:NRHO-Halo-ICs} and the problem parameters in \cref{tab:NRHO}. 
The NRHO is the 9:2 NRHO for the ongoing ARTEMIS program; 
the departure is from the apolune (southernmost point) of the NRHO, and the arrival is at one of the intersections between the L1 Halo and the $xz$-plane. 
The reference trajectory and control profile are shown in \cref{fig:NRHOtoHalo-reference}. 
The process of obtaining the reference trajectory is similar to the previous example; 
however, since this transfer features a close pass to the moon,
the indirect method \textit{regularizes} the equations of motion based on \cite{oguriRegularizationCircularRestricted2024}
to avoid the singularity in the inverse-square law. 
Furthermore, since the spacecraft's thrust capability is lower,
we convert the minimum-energy trajectory to a fuel-optimal trajectory
 by setting a high-thrust magnitude and gradually lowering the magnitude, 
for each step solving the deterministic trajectory optimization problem with \texttt{SCvx*}
using the previous solution as the initial guess.

\begin{figure}[htbp]
    \centering
    \begin{subfigure}[c]{\textwidth}
    \includegraphics[width=\linewidth]{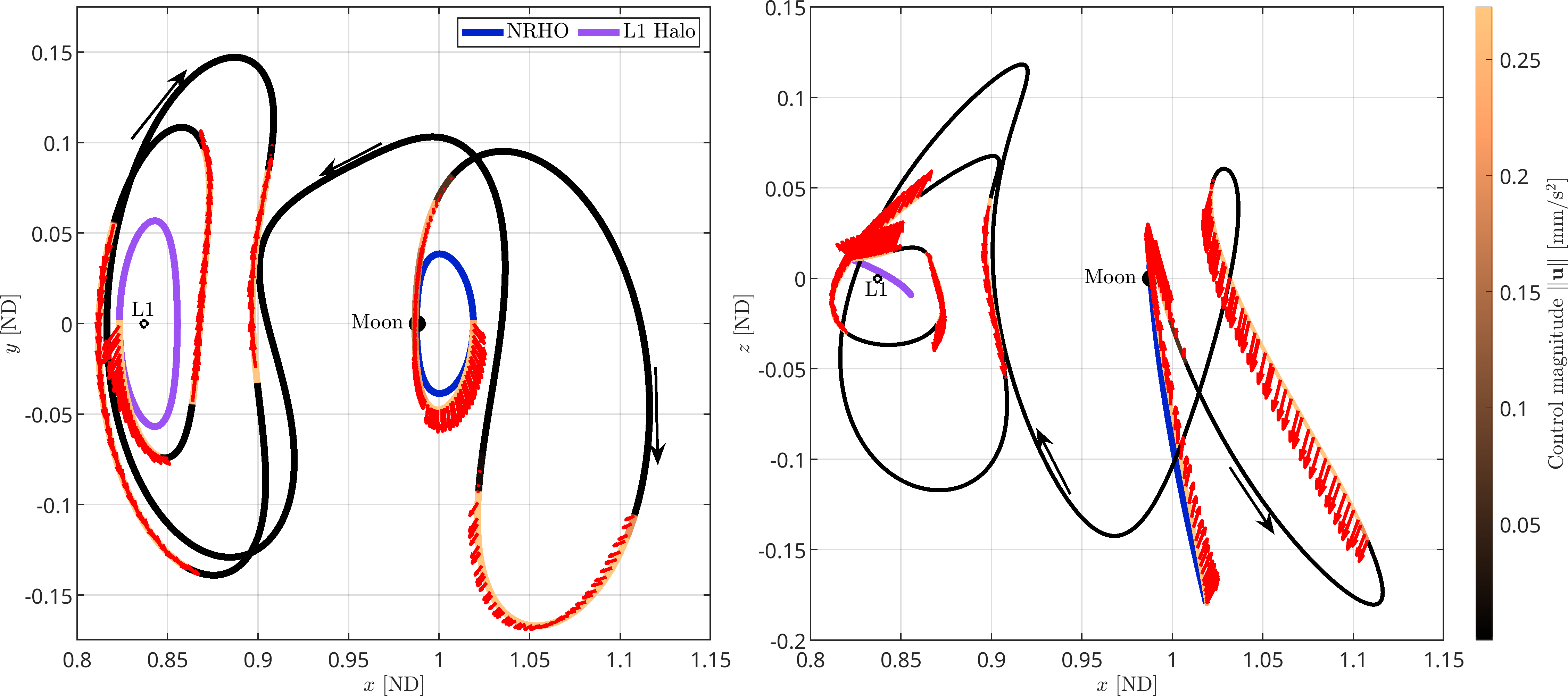}
    \subcaption{Trajectory $xy$ and $xz$ projection} 
    \end{subfigure}
\end{figure}%
\begin{figure}[htbp]\ContinuedFloat
    \begin{subfigure}[c]{0.53\textwidth}
        \includegraphics[width=\linewidth]{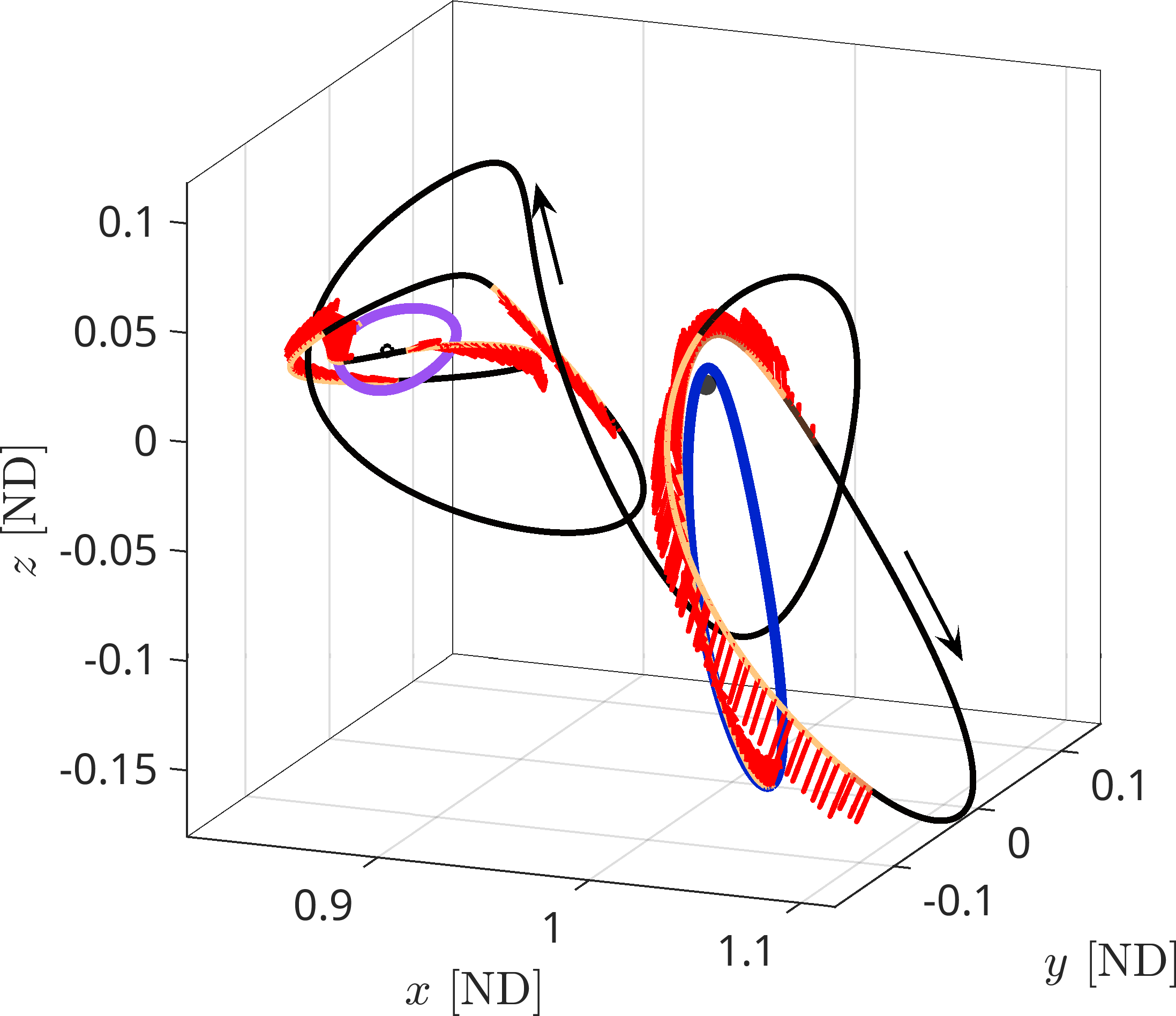}
        \subcaption{Trajectory 3D view}
    \end{subfigure}%
    \hfill%
    \begin{subfigure}[c]{0.45\textwidth}
    \includegraphics[width=\linewidth]{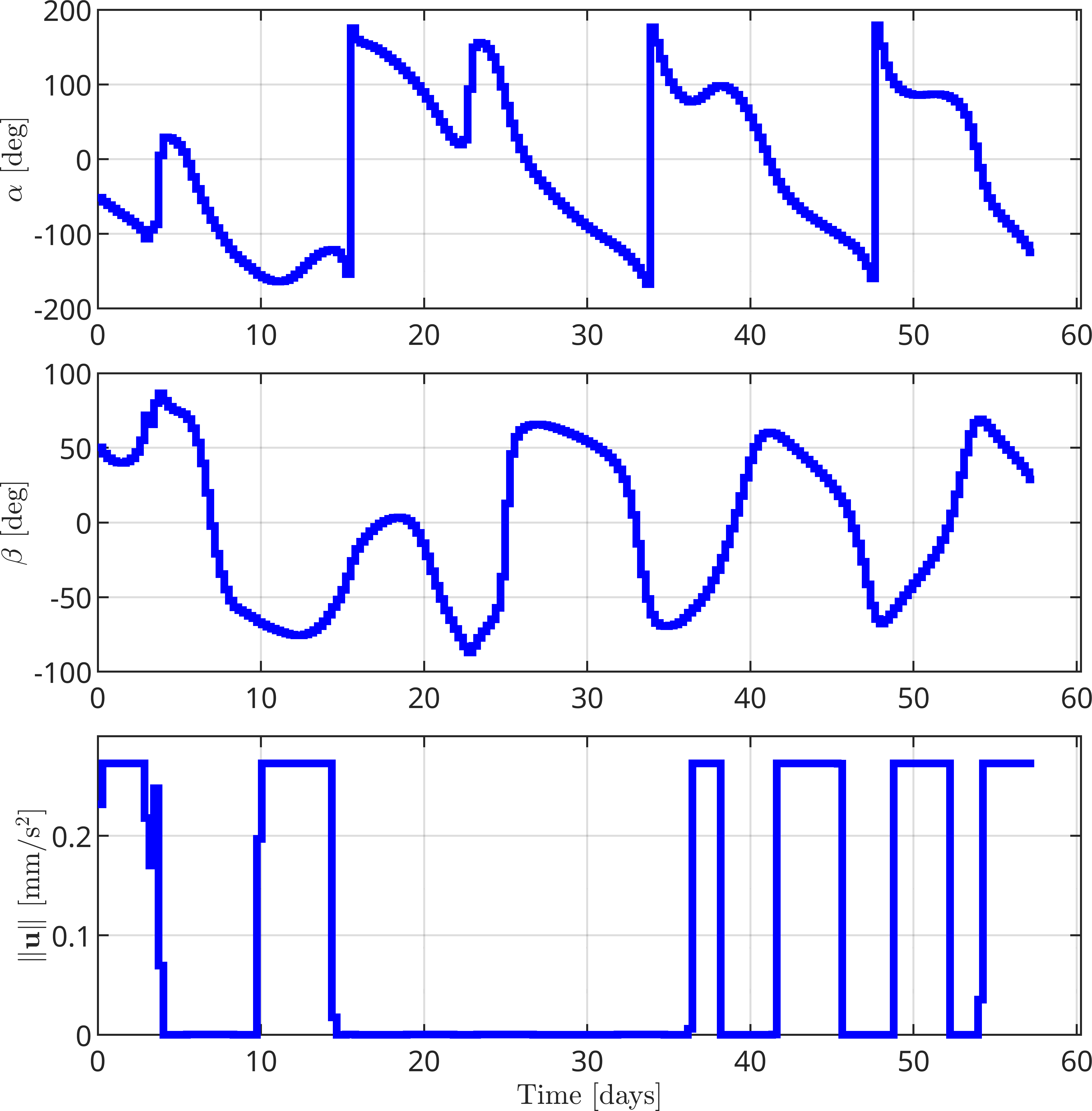}
    \subcaption[]{Control}
    \end{subfigure}
    \caption{Reference fuel-optimal L2 NRHO to L1 Halo transfer.}
    \label{fig:NRHOtoHalo-reference}
\end{figure}

First, we remark on some characteristics of the reference trajectory. 
The trajectory can be divided into three phases, 
1. NRHO departure, 2. heteroclinic connection, and 3. Halo arrival. 
Visual inspection tells us that the first phase utilizes dynamical structures near the NRHO 
and resembles the P4HO2 family of periodic orbits \cite[Fig. 9]{zimovan-spreenRectilinearHaloOrbits2020}
 that bifurcate from the NRHO. 
 The second phase is a coasting phase, mirroring a heteroclinic connection. 
 The final phase is the `spiral down' towards the target Halo orbit. 
 This phase seems to utilize the Lissajous orbit structures around the L1 point, 
 and is similar to the one used in the ARTEMIS mission \cite[Fig. 10]{sweetserARTEMISMissionDesign2014}, 
 although the transfer here has larger inclination changes and features less revolutions.
 We remark that the trajectory optimization with SCP finds a trajectory that leverages
 the dynamical structures of the CR3BP without specific user inputs or dividing 
 the transfer into multiple segments.

Next, we apply \cref{alg:SCP} to the reference trajectory. In this scenario, 
we impose a maximum covariance constraint in order to maintain the validity
of the linear covariance model.
The algorithm converges in 60 iterations in 10 minutes. 
The larger number of iterations required compared to the previous example can be attributed to the higher nonlinearity
in the lunar flyby phase and as a result, the noticeable change in the nominal solution from the initial reference; this transfer shows a larger discrepancy between the deterministic and robust trajectories compared to the DRO--DRO transfer.
Figure \ref{fig:NRHOtoHalo-trajcomparison} shows comparisons between the initial reference trajectory and the obtained nominal (mean) trajectory. 
We see notable differences between the two. 
As we see later, the robust trajectory devotes a significant amount of TC effort to 
shrink the state dispersion at the beginning of the transfer, especially since the 
lunar flyby phase is soon after the departure from the NRHO. Since the entire control 
term is bounded by the maximum control constraint, the nominal thrust magnitude is
smaller than that of the deterministic trajectory. This is also reflected in the control profile in \cref{fig:NRHOtoHalo-control-profile}.
As a result, the robust trajectory especially has a larger movement in the $y$ direction 
after the lunar flyby, and the heteroclinic connction uses a trajectory that is 
closer to the $xz$-plane. The difference in trajectories becomes smaller as 
the spacecraft approaches the target Halo orbit.
\begin{figure}[htbp]
    \centering
    \begin{subfigure}[c]{\textwidth}
    \includegraphics[width=\linewidth]{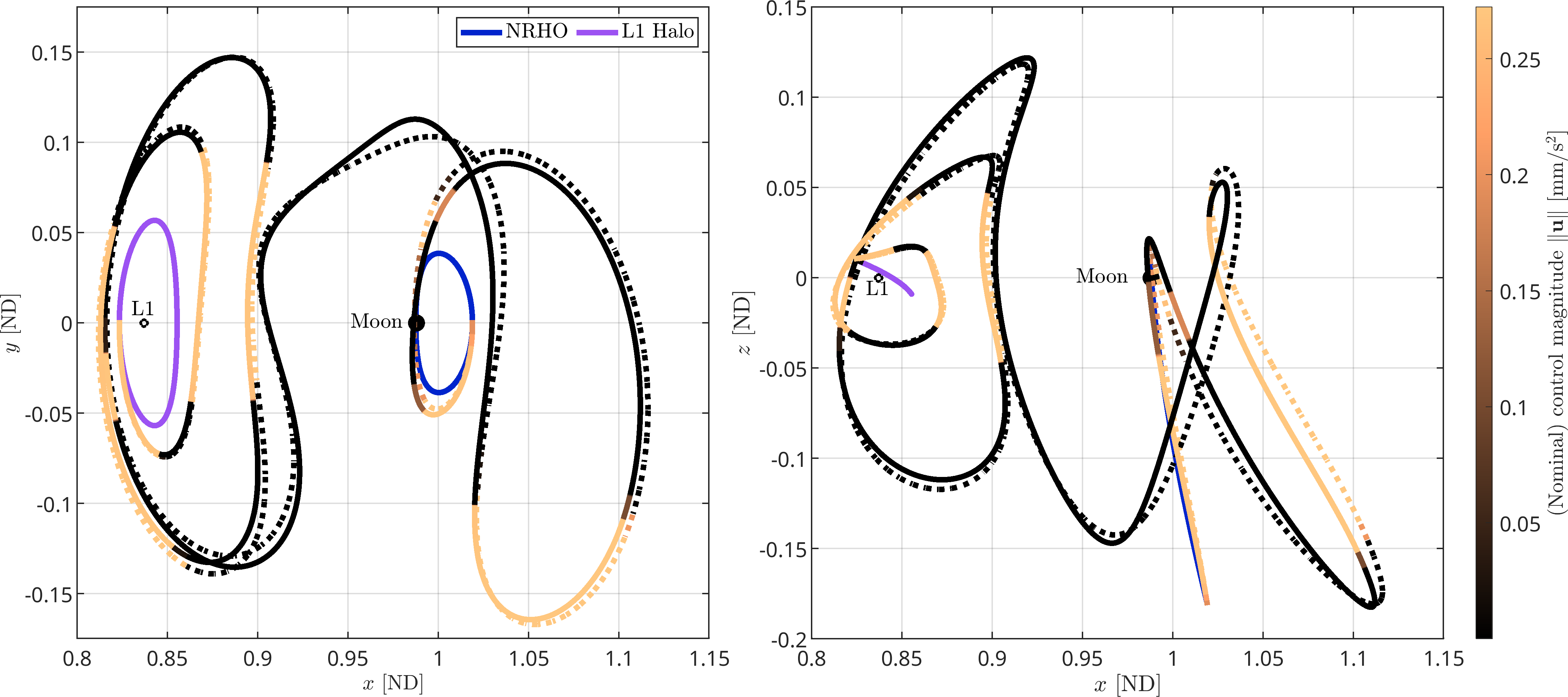}
    \caption{$xy$ and $xz$ projection. Broken lines show the trajectory from \cref{fig:NRHOtoHalo-reference}; 
    solid lines show the robust trajectory.}
    \end{subfigure}
\end{figure}%
\begin{figure}[htbp]\ContinuedFloat
    \begin{subfigure}[c]{0.25\textwidth}
    \includegraphics[width=\linewidth]{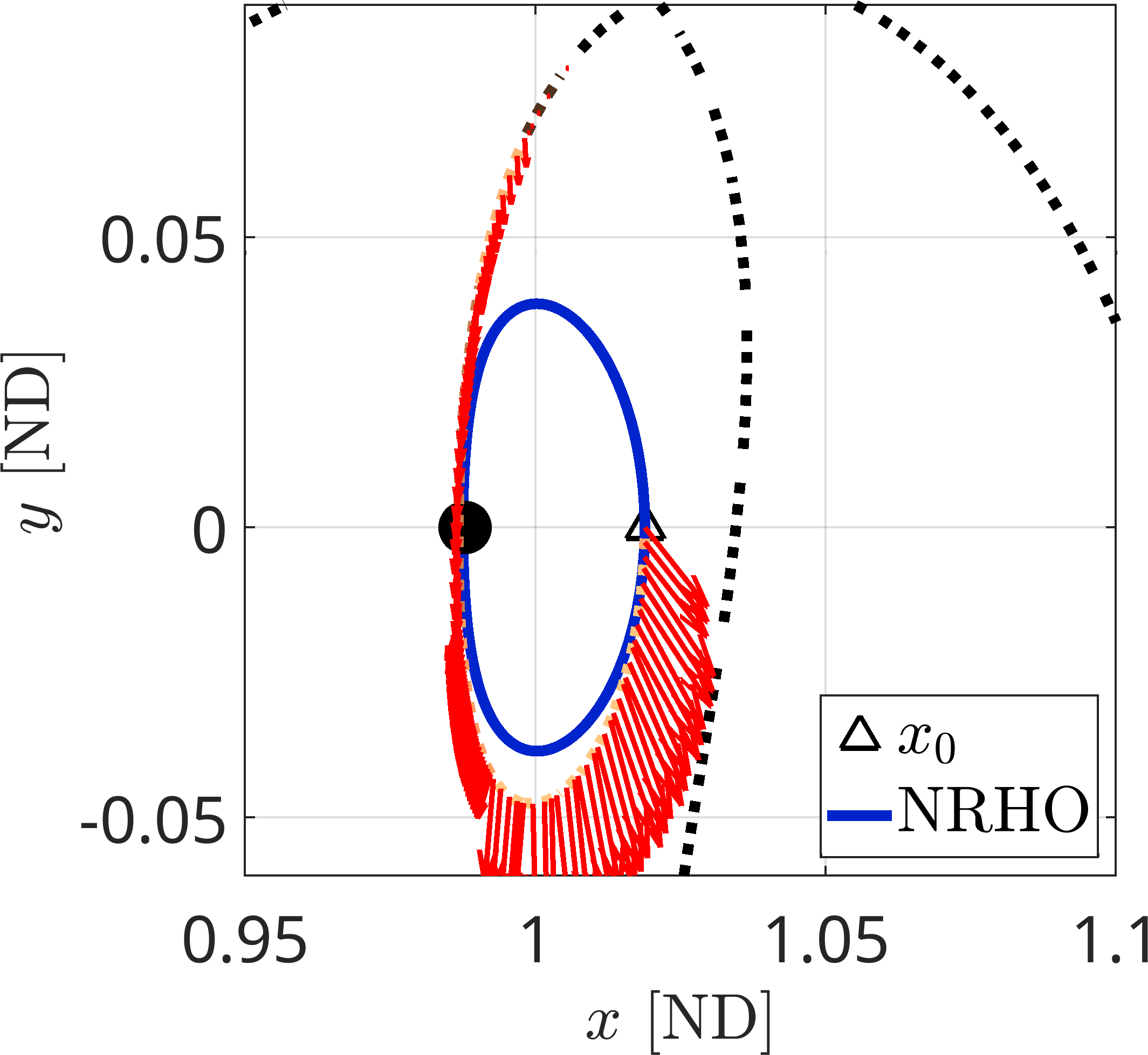}
    \caption{$xy$, deterministic}    
    \end{subfigure}%
    \begin{subfigure}[c]{0.25\textwidth}
    \includegraphics[width=\linewidth]{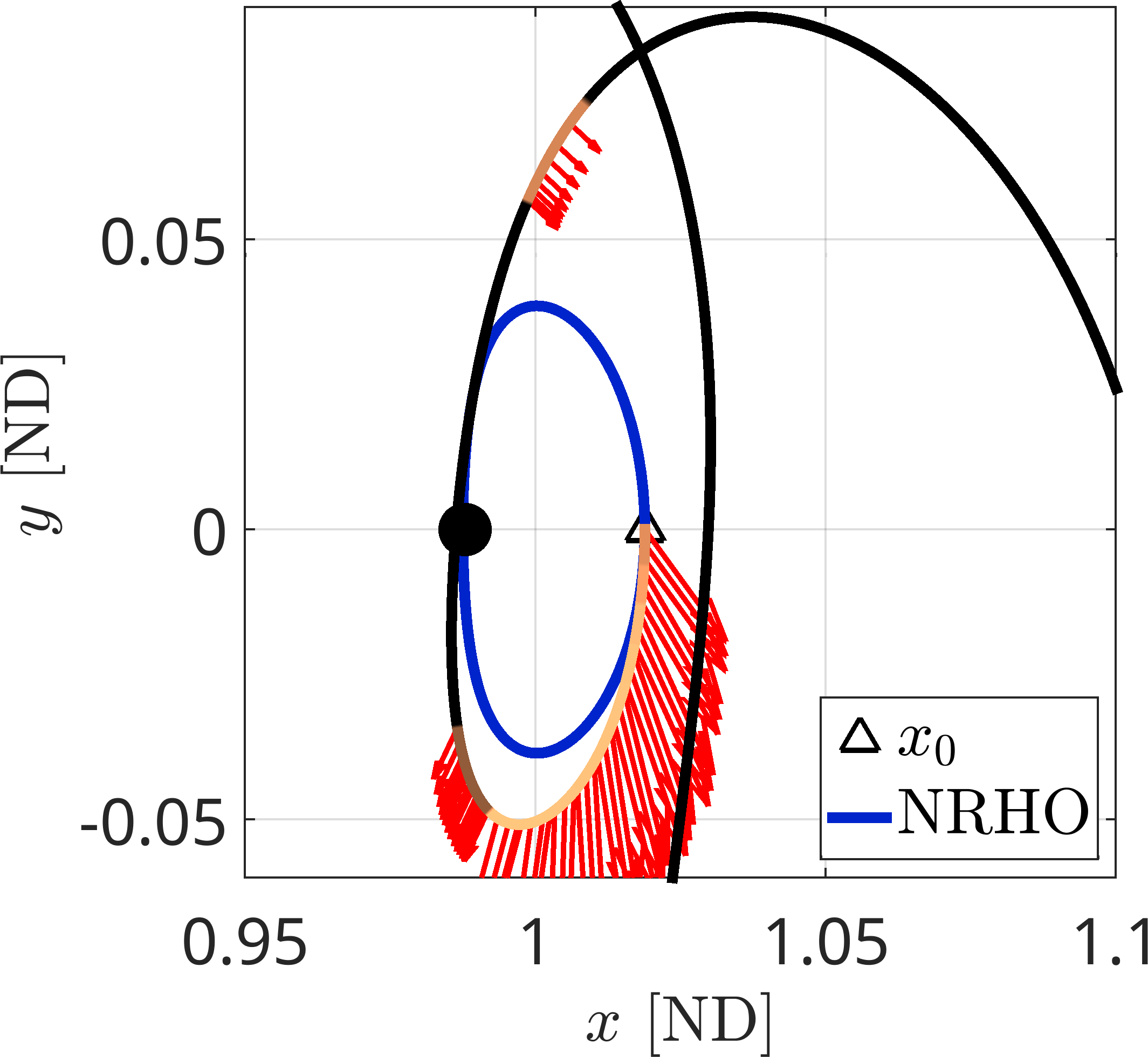}
    \caption{$xy$, robust}    
    \end{subfigure}%
    \begin{subfigure}[c]{0.25\textwidth}
    \includegraphics[width=\linewidth]{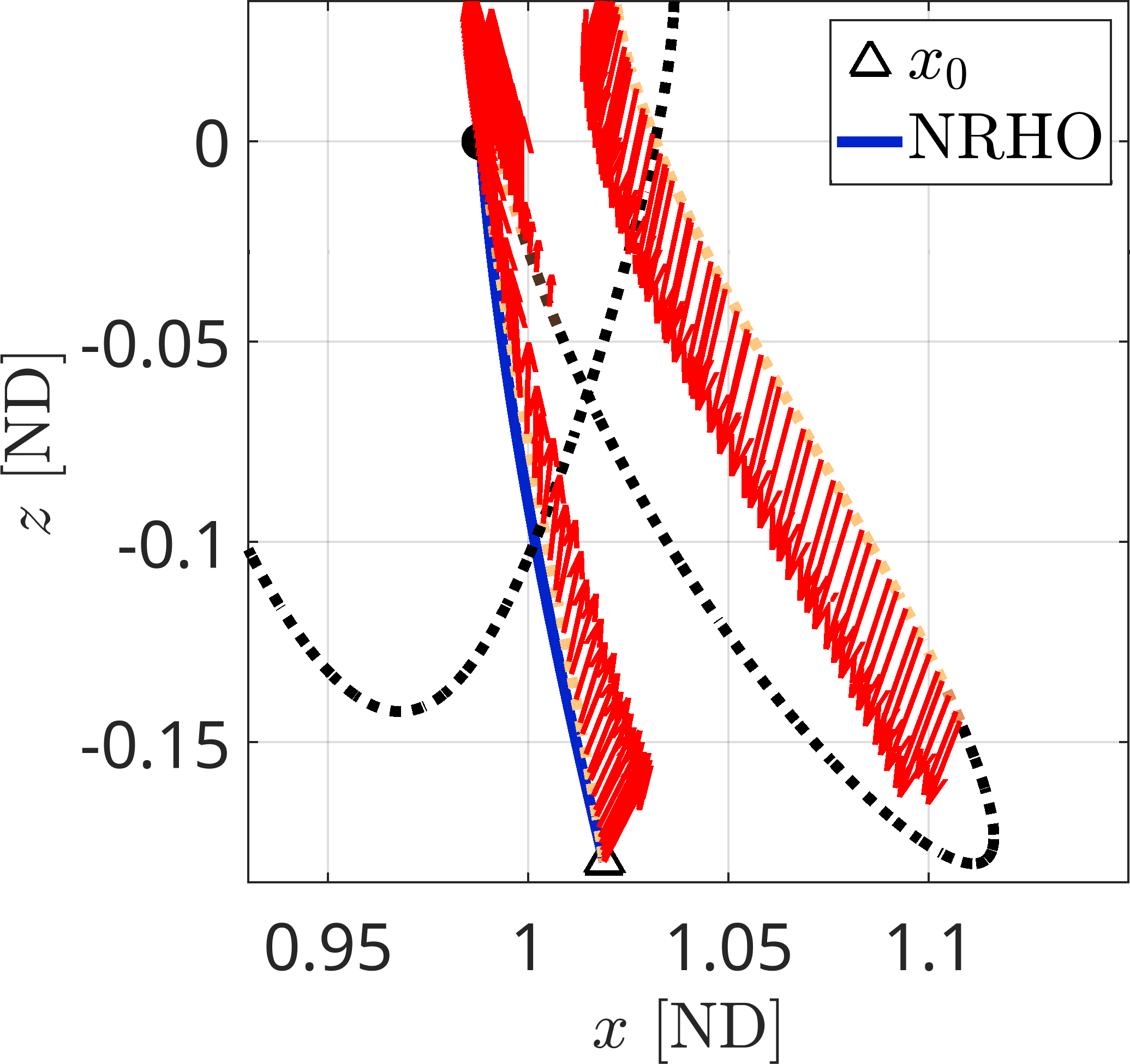}
    \caption{$xz$, deterministic}    
    \end{subfigure}%
    \begin{subfigure}[c]{0.25\textwidth}
    \includegraphics[width=\linewidth]{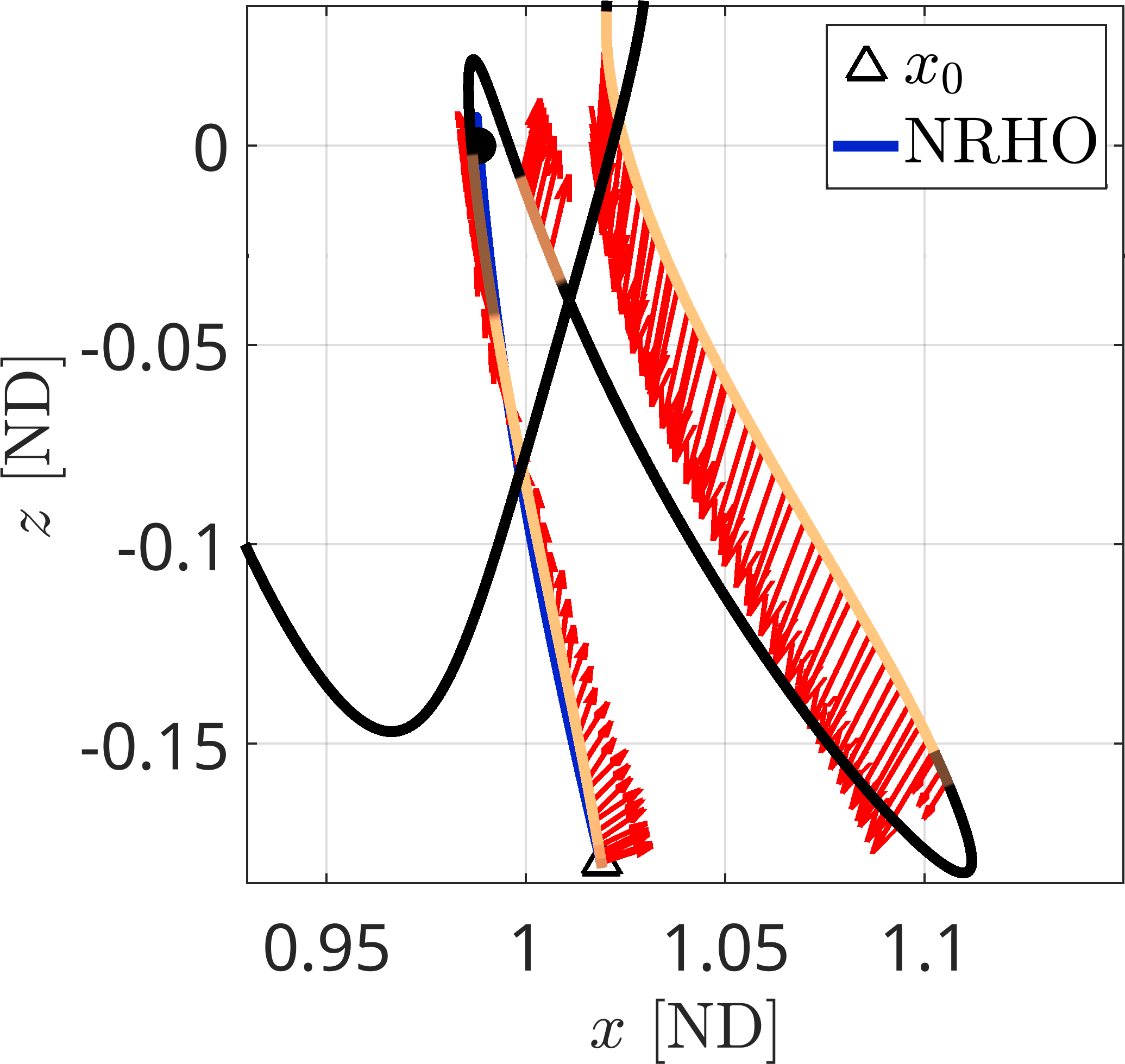}
    \caption{$xz$, robust}    
    \end{subfigure}
    \caption{Comparison of the deterministic vs robust trajectory. 
    (a): Projections of both trajectories. 
    (b)-(e): Zoomed in views of the trajectories and thrust vectors around the lunar flyby}
    \label{fig:NRHOtoHalo-trajcomparison}
\end{figure}

\begin{figure}[htbp]
    \centering
    \includegraphics[width=\linewidth]{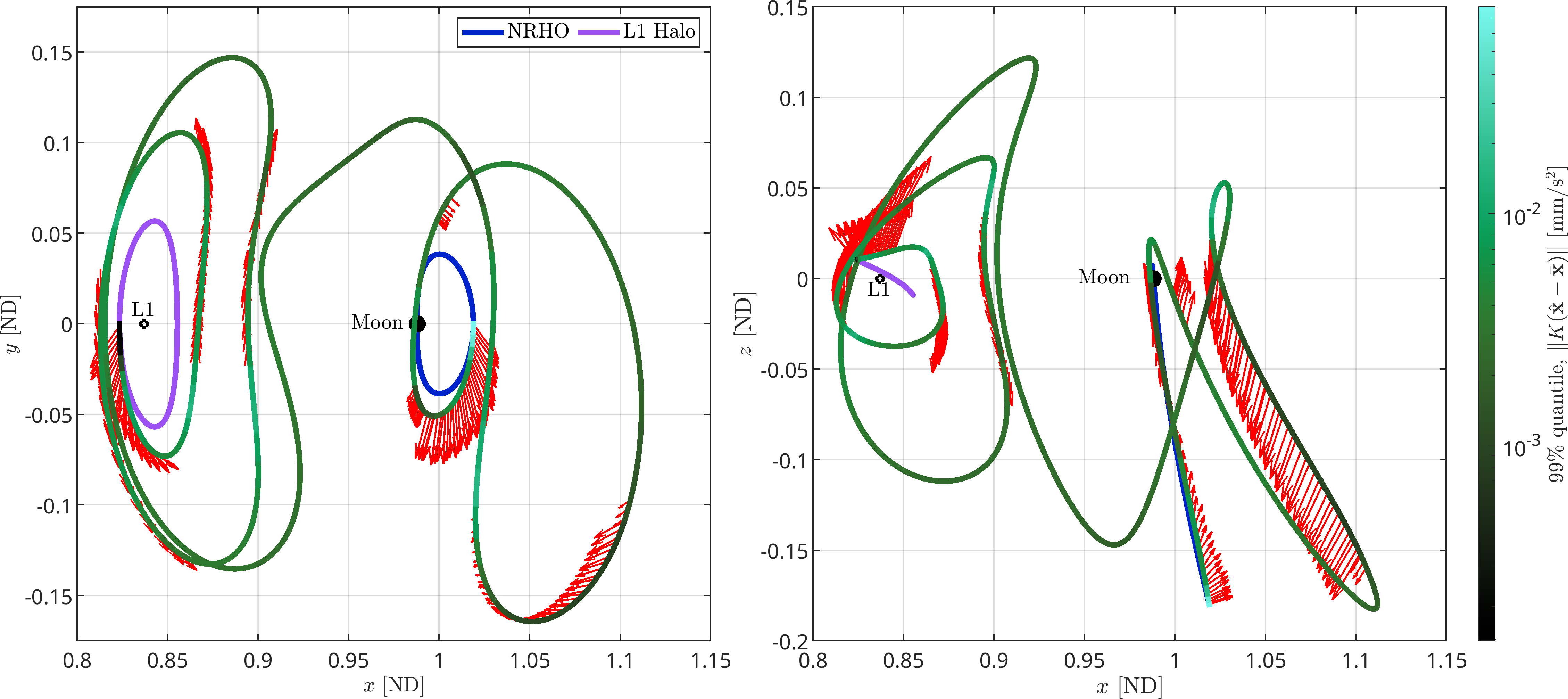}
    \caption{Quantile trajectory correction magnitudes along the trajectory; red arrows show the nominal thrust acceleration}
    \label{fig:NRHOtoHalo-feedback_along_trajectory}
\end{figure}

\begin{figure}[htbp]
    \centering
    \includegraphics[width=0.7\linewidth]{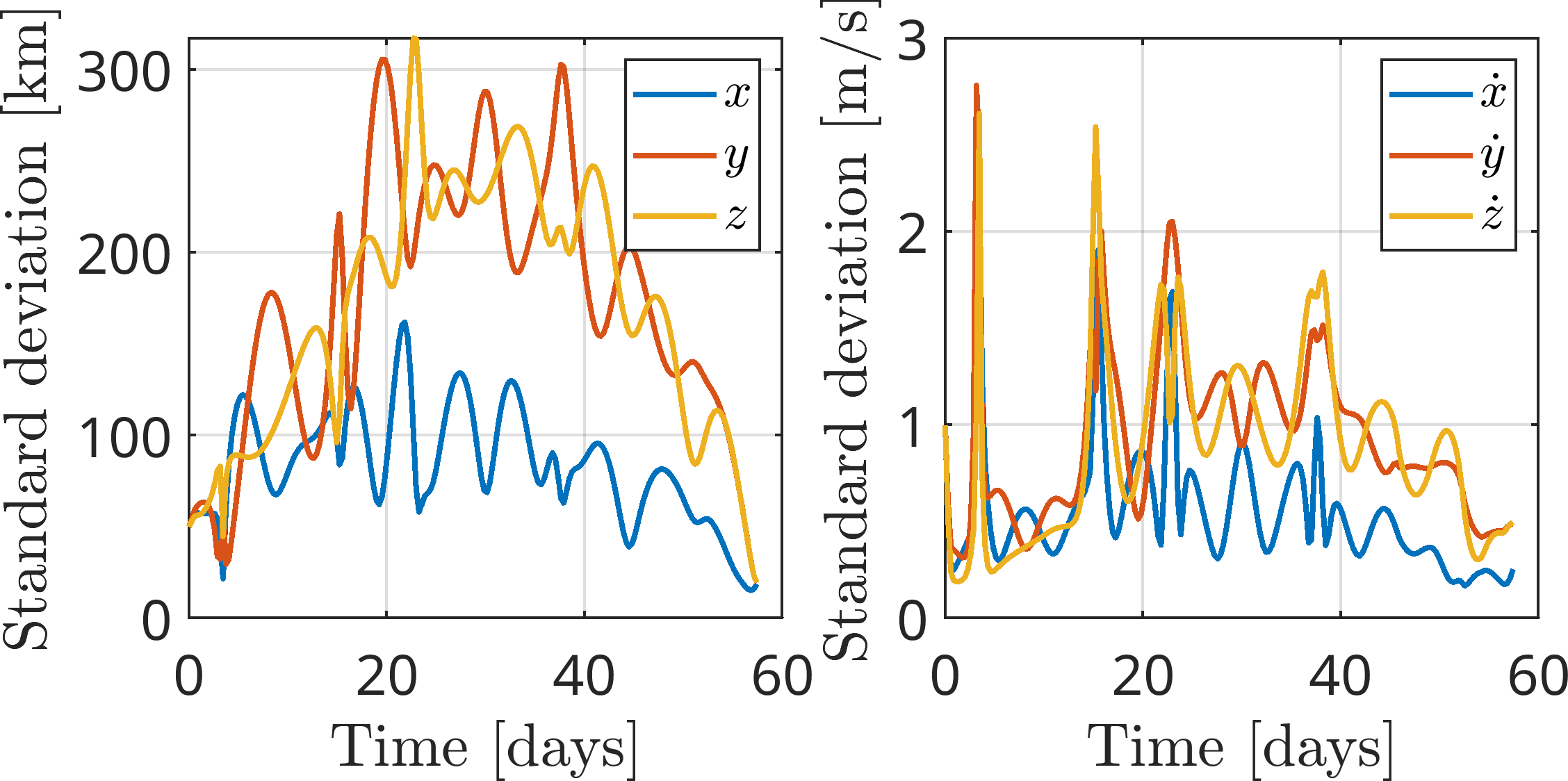}
    \caption{State standard deviations (Square root of diagonals of $P$) for NRHO--Halo transfer}
    \label{fig:NRHOtoHalo-standard-deviations}
\end{figure}

Next, we observe the TC magnitudes along the trajectory.
Figure \ref{fig:NRHOtoHalo-feedback_along_trajectory} shows the TCM magnitudes (in log scale) along the trajectory ($\norm{K (\hat{\bm{x}} - \overline{\bm{x}})}$).
The largest TC magnitudes are, in order, the initial departure from the NRHO, the lunar flyby, 
and the final spiral down to the Halo orbit. 
To understand this, we observe the state standard deviations in \cref{fig:NRHOtoHalo-standard-deviations}.
During the lunar flyby phase, the maximum covariance constraint is active in the velocity components, i.e.
the diagonals of the covariance matrix that correspond to the velocity closely approach the 
corresponding indices of the maximum covariance matrix.
As also verified in \cref{sec:LLE}, the trajectory passes through a highly error-sensitive region
quickly after the lunar flyby, which explains the large TC magnitude leading up to this phase.
The heteroclinic connection has small TC magnitudes as well as the nominal thrusting magnitude.
Towards the end of the transfer, the TC magnitudes increase to satisfy the final covariance constraint.

\begin{figure}[htbp]
    \centering
    \begin{subfigure}[c]{0.49\textwidth}
        \includegraphics[width=\textwidth]{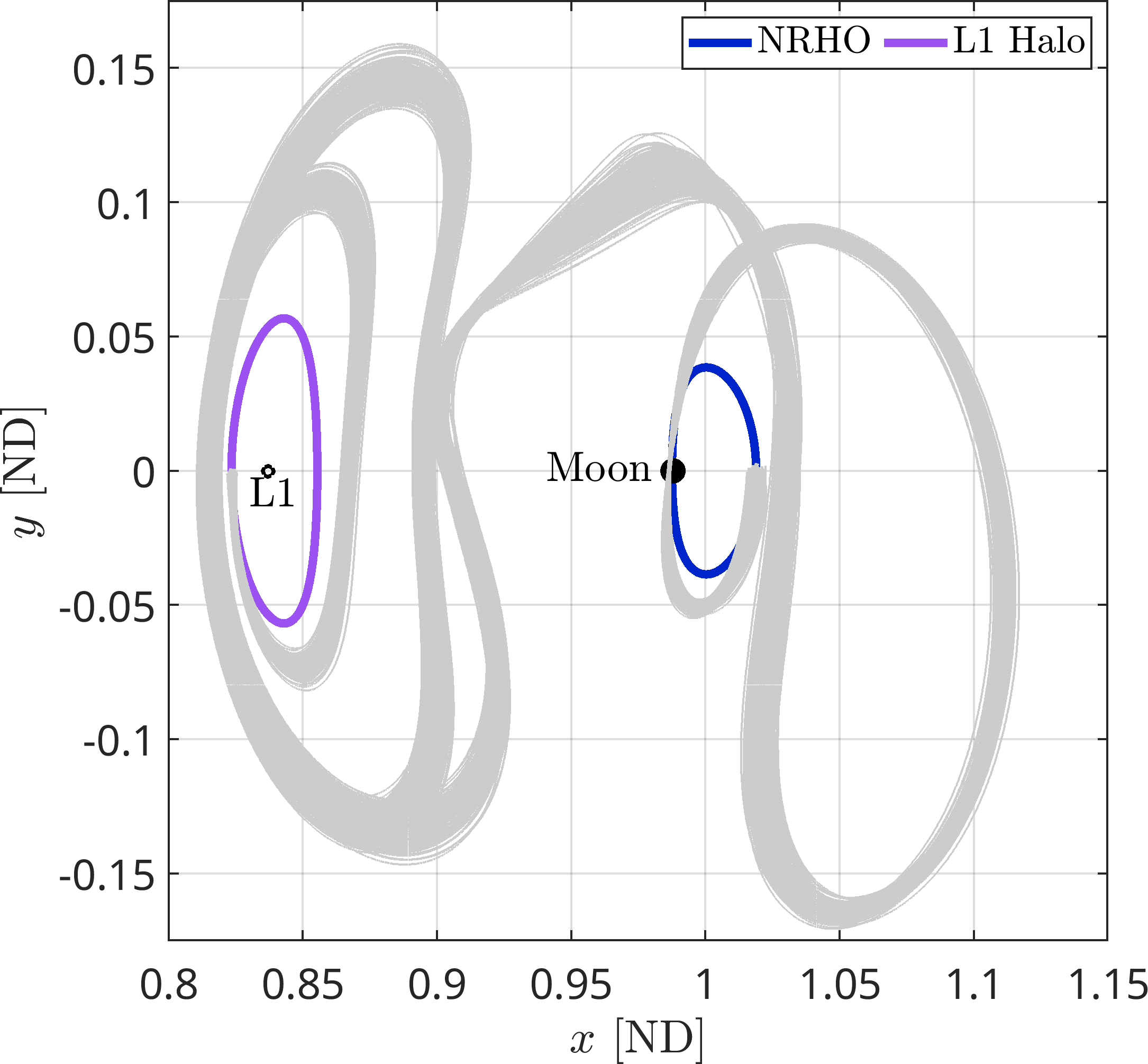}
        \caption{Trajectories ($xy$ projection, Deviations enlarged 10x)}
        \label{fig:NRHOtoHalo-robust-traj-xy}
    \end{subfigure}
    \hfill%
    \begin{subfigure}[c]{0.49\textwidth}
        \includegraphics[width=\textwidth]{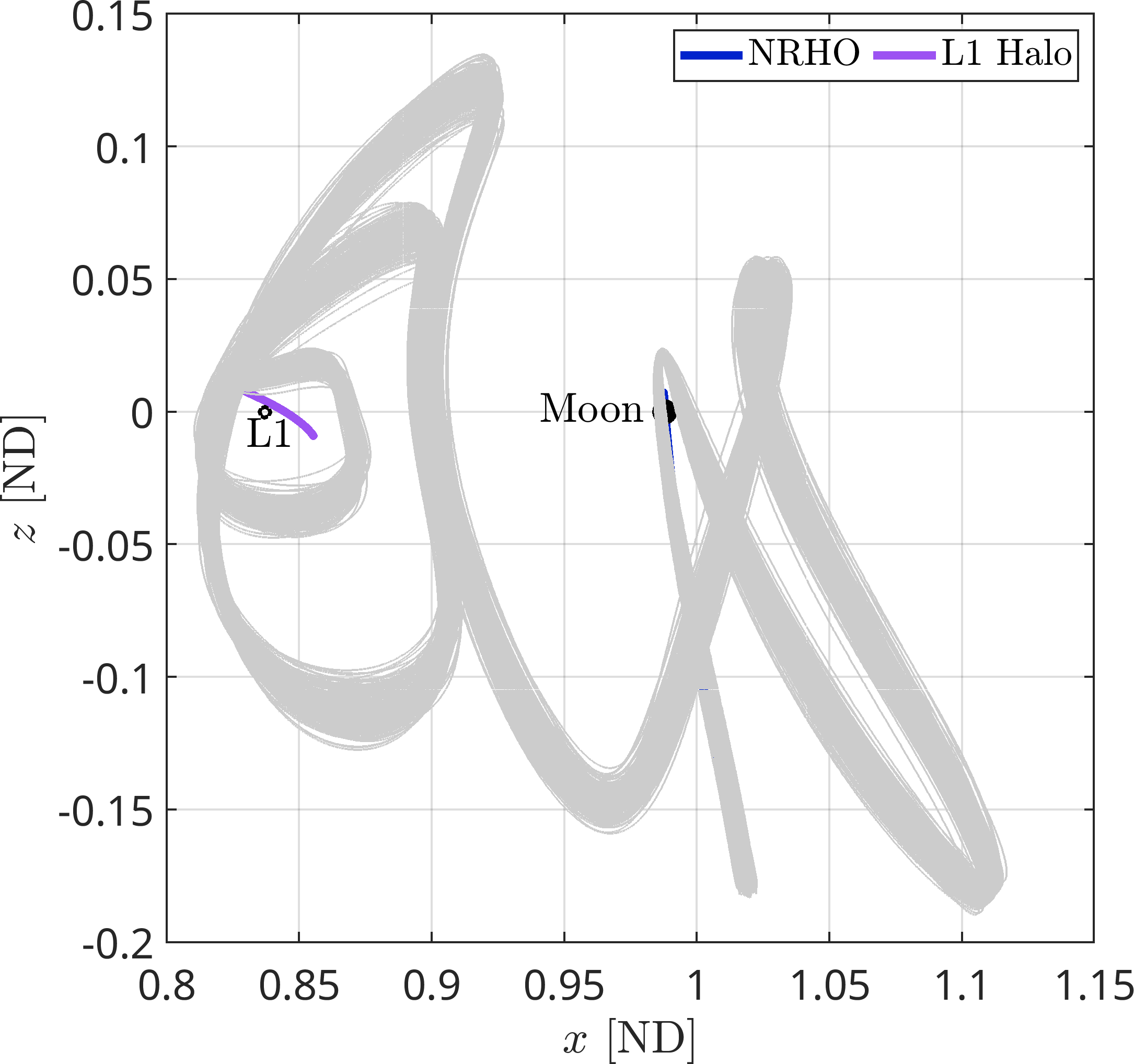}
        \caption{Trajectories ($xz$ projection, Deviations enlarged 10x)}
        \label{fig:NRHOtoHalo-robust-traj-xz}
    \end{subfigure}
    \begin{subfigure}[c]{0.49\textwidth}
        \includegraphics[width=\textwidth]{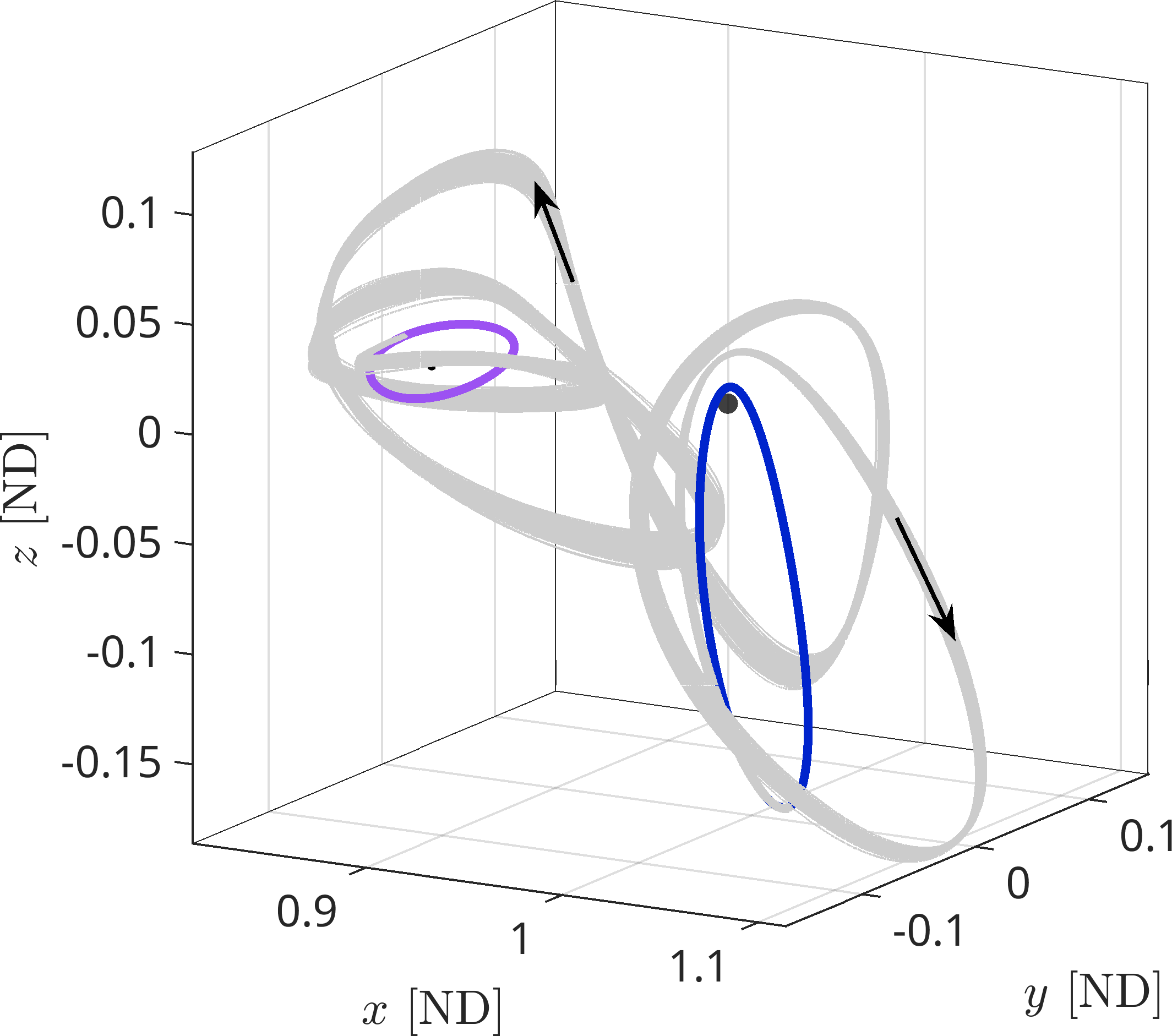}
        \caption{Trajectories ($xz$ projection, Deviations enlarged 5x)}
        \label{fig:NRHOtoHalo-robust-traj-3d}
    \end{subfigure}%
    \hfill%
    \begin{subfigure}[c]{0.49\textwidth}
        \includegraphics[width=\textwidth]{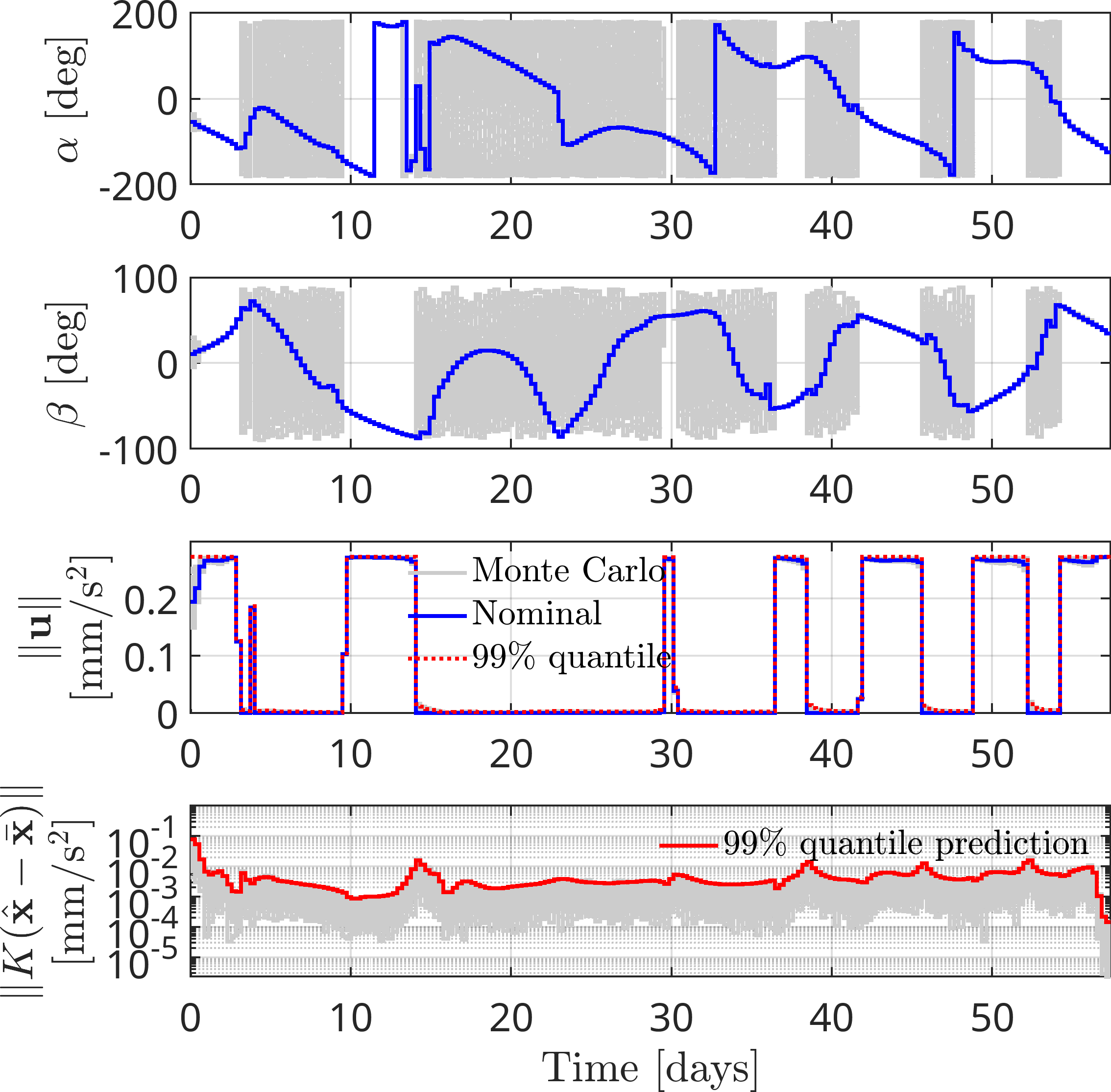}
        \caption{Control profile}
        \label{fig:NRHOtoHalo-control-profile}
    \end{subfigure}
\end{figure}

\begin{figure}[htbp] \ContinuedFloat
    \centering
    \begin{subfigure}[c]{0.49\textwidth}
        \includegraphics[width=\textwidth]{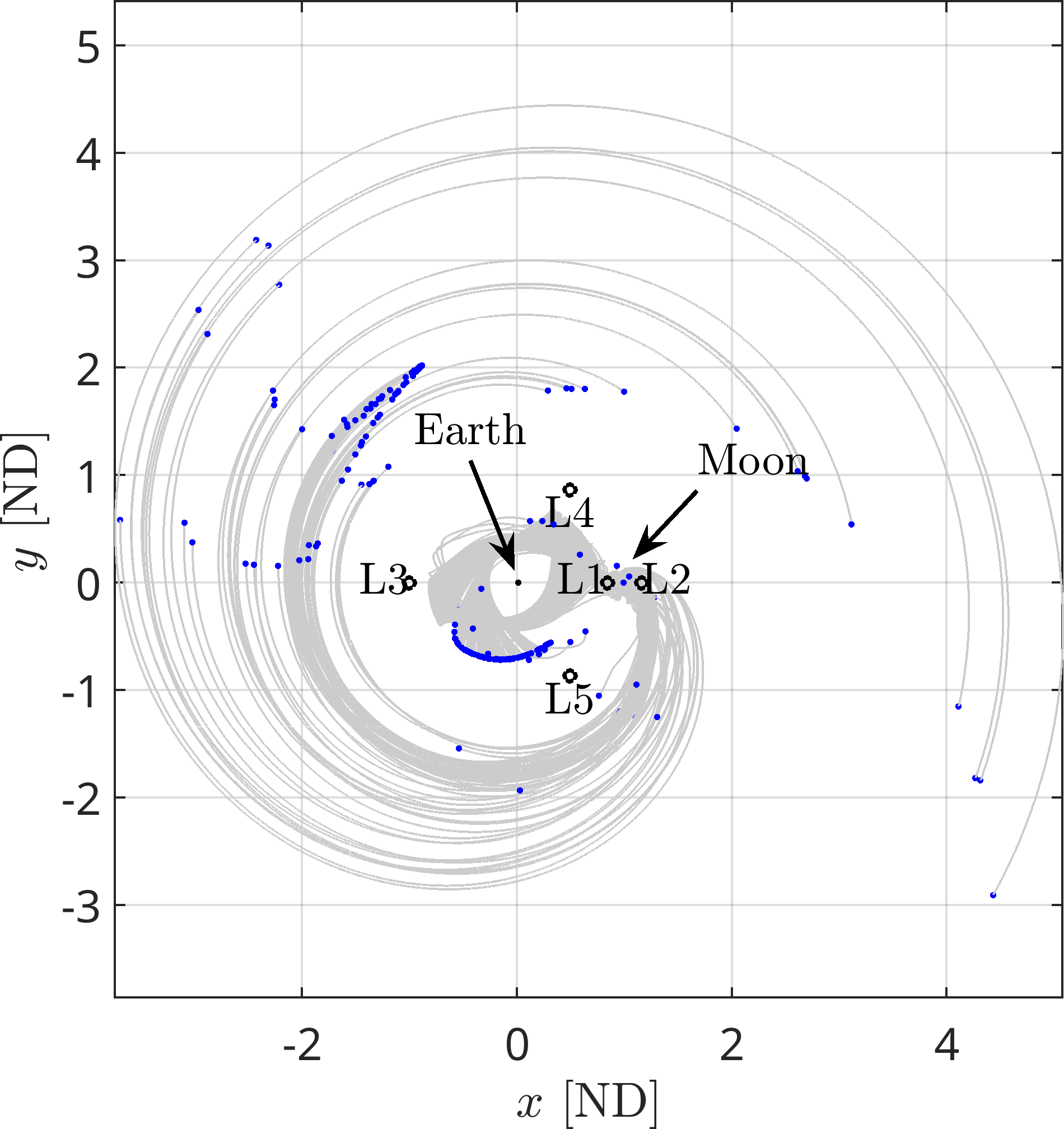}
        \caption{\highlight{Trajectories (open-loop, $xy$ projection); L1-L5 show the Lagrange points}}
        \label{fig:NRHOtoHalo-robust-traj-openloop-xy}
    \end{subfigure}%
    \caption{Monte Carlo simulation results for NRHO--Halo transfer.}
    \label{fig:NRHO-Halo-robust}
\end{figure}

Again, we perform Monte Carlo simulations with 200 samples.
Figure \ref{fig:NRHO-Halo-robust} shows the results along with the derived control policy. 
Figures \ref{fig:NRHOtoHalo-robust-traj-xy},\ref{fig:NRHOtoHalo-robust-traj-xz}, and \ref{fig:NRHOtoHalo-robust-traj-3d} show the trajectories of the Monte Carlo samples.
Around the lunar flyby, the samples are tightly bound thanks to the TCs.
Recall that the maxmimum covariance constraint is active during this period of the transfer.
Note that in in this problem setting, setting the largest covariance constraint is crucial for the success of the algorithm;
without this constraint, the state dispersion can grow outside the region where linear covariance propagation is 
a valid approximation of the Fokker-Planck equation, resulting in the Monte Carlo samples escaping the cislunar region 
even if the SCP process converges.
We also observe that the deviation expands near the heteroclinic connection, where the TC magnitude is \highlight{relatively small}.
This is an interesting phenomenon; it may be due to a dynamical structure, or it may simply be due to the
heteroclinic connection being in the middle of the transfer, where efforts are not focused on satisfying the final covariance constraint.
Observing the control profile in \cref{fig:NRHOtoHalo-control-profile}, we see that, as 
observed in \cref{fig:NRHOtoHalo-feedback_along_trajectory}, the TC magnitude is largest at the beginning of the transfer
\highlight{leading up to the lunar flyby at around day 4. We also see the same trend as observed in the DRO--DRO transfer, 
where the TC effort grows towards the end of each thrusting arc.}
\highlight{\cref{fig:NRHOtoHalo-robust-traj-openloop-xy} shows the open-loop case. Compared to the DRO--DRO transfer,
the effect of not applying trajectory corrections is severe, even though the additive noise is an order of magnitude smaller (see \cref{tab:parameters}). 
Many of the samples escape the cislunar region via the L2 gateway.
}
\subsection{Analysis via Local Lyapunov Exponents} \label{sec:LLE}
So far, we have observed that the CCDVO framework successfully provides optimized nominal (mean) trajectories as well as TCM policies. 
Our hypothesis is that the trajectories obtained with this framework are in general more robust regardless of the trajectory correction policy that is used. 
Here, we provide some analysis of the stablity along the nominal trajectory computed from \cref{alg:SCP} via the local Lyapunov exponent.
The finite-time local Lyapunov exponent (LLE) is defined as \cite{andersonApplicationLocalLyapunov2003}
\begin{equation}
    \Lambda = \frac{1}{\Delta t} \ln\|\Phi(t + \Delta t, \mathhl{t}) \|_2 = \frac{1}{\Delta t} \ln \sqrt{\lambda_{\max} \Phi(t + \Delta t, \mathhl{t}) \Phi^\top(t + \Delta t, t)}
\end{equation}
with $\Phi(t + \Delta t, t)$ representing the state transition matrix from time $t$ to $t + \Delta t$. 
While the horizon length $\Delta t$ is a parameter which will effect the value of LLE,
 we have seen that the general trends in the LLE along the trajectory are fairly robust to the choice of $\Delta t$. 
 For our analysis we choose $\Delta t$ to be equal to the time step of the discretization.

\begin{figure}[htbp]
    \centering
    \begin{subfigure}[c]{0.49\textwidth}
        \includegraphics[width=\textwidth]{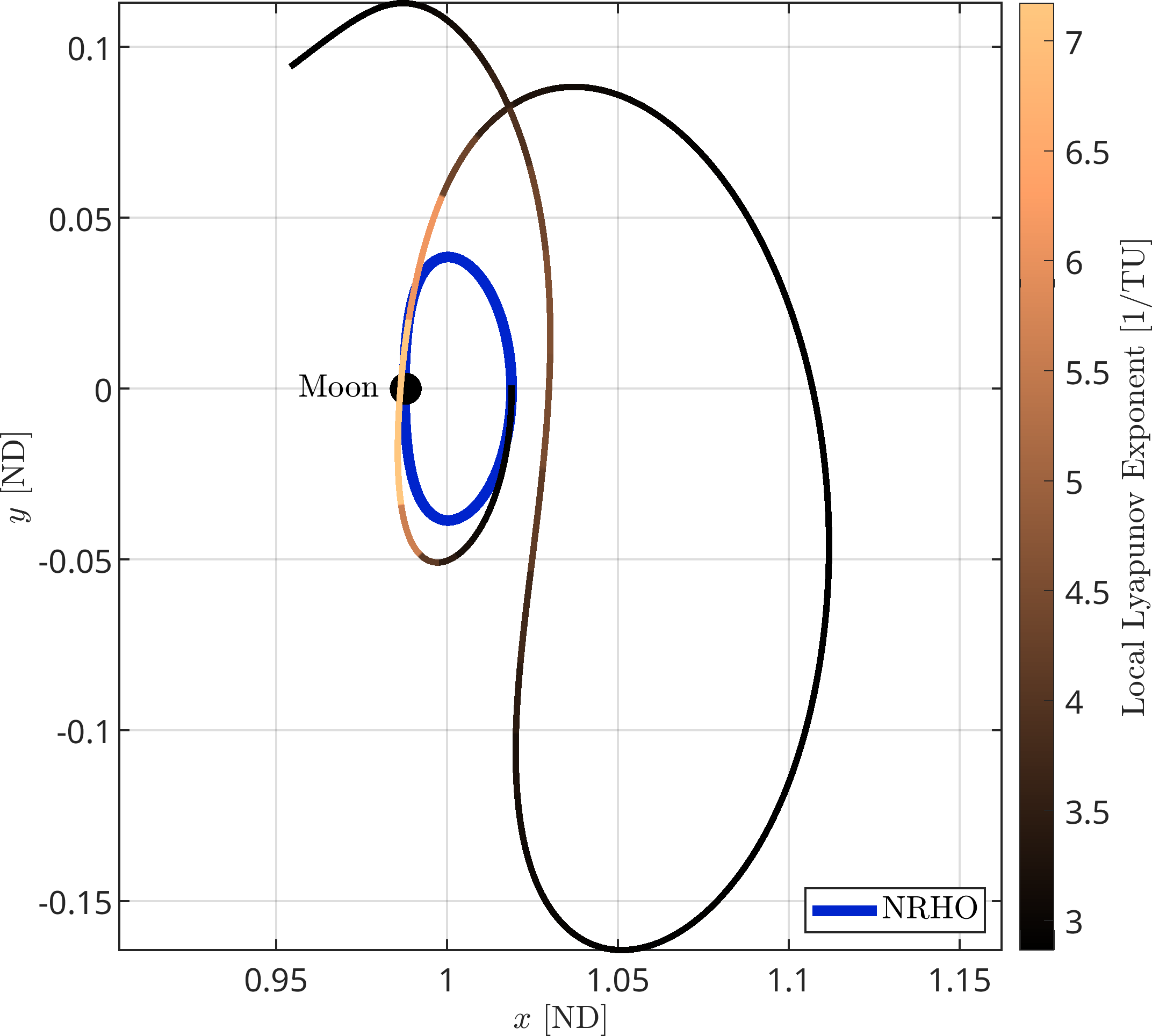}
        \caption{First half of the trajectory, $xy$ projection}
    \label{fig:LLE_firsthalf}
    \end{subfigure}%
    \hfill
    \begin{subfigure}[c]{0.49\textwidth}
        \includegraphics[width=\textwidth]{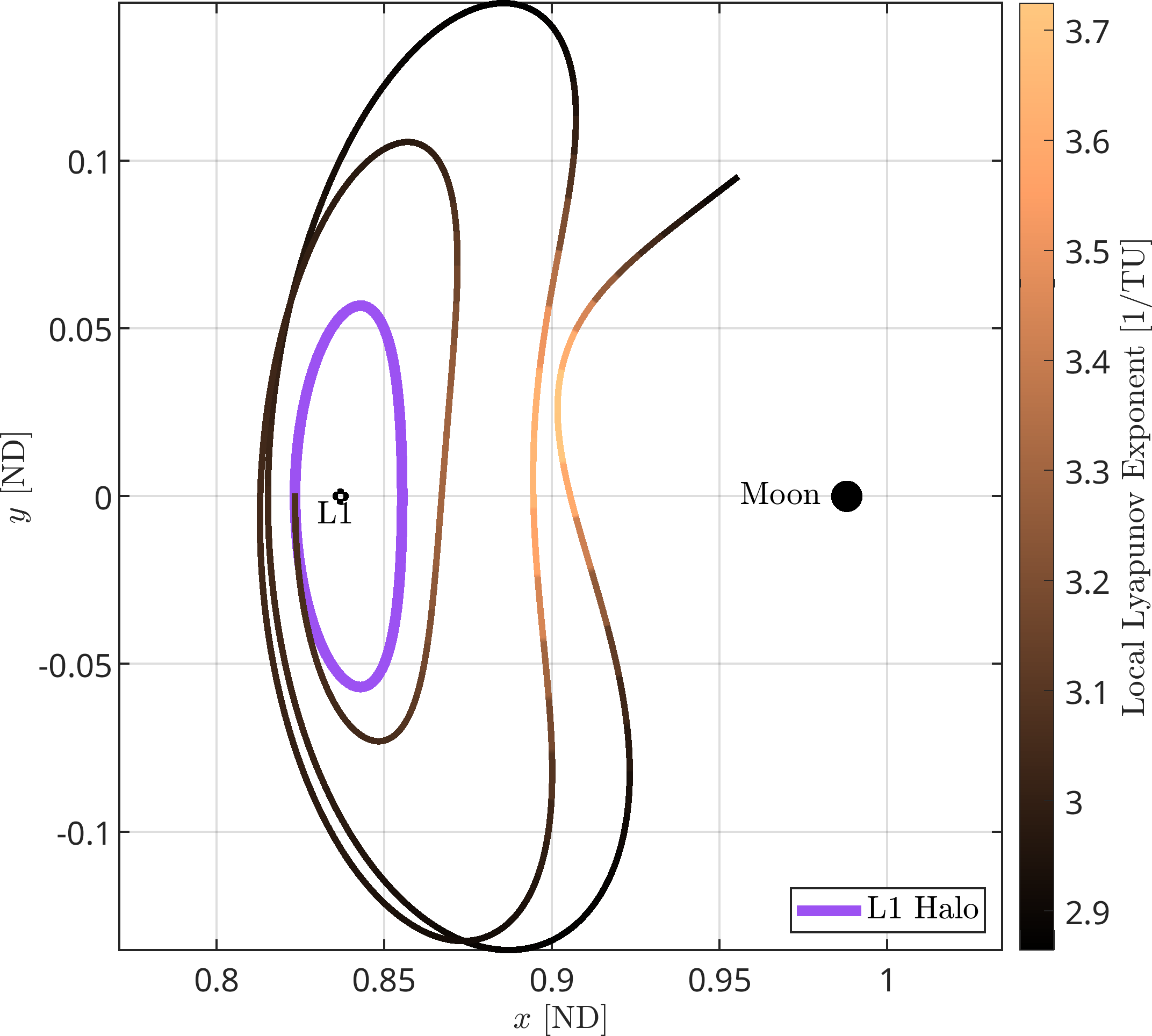}
        \caption{Second half of the trajectory, $xy$ projection}
    \label{fig:LLE_secondhalf}
    \end{subfigure}
    \caption{Finite-time local Lyapunov exponent along the robust trajectory}
    \label{fig:LLE}
\end{figure}
\begin{figure}[htbp]
    \centering
    \includegraphics[width=0.45\linewidth]{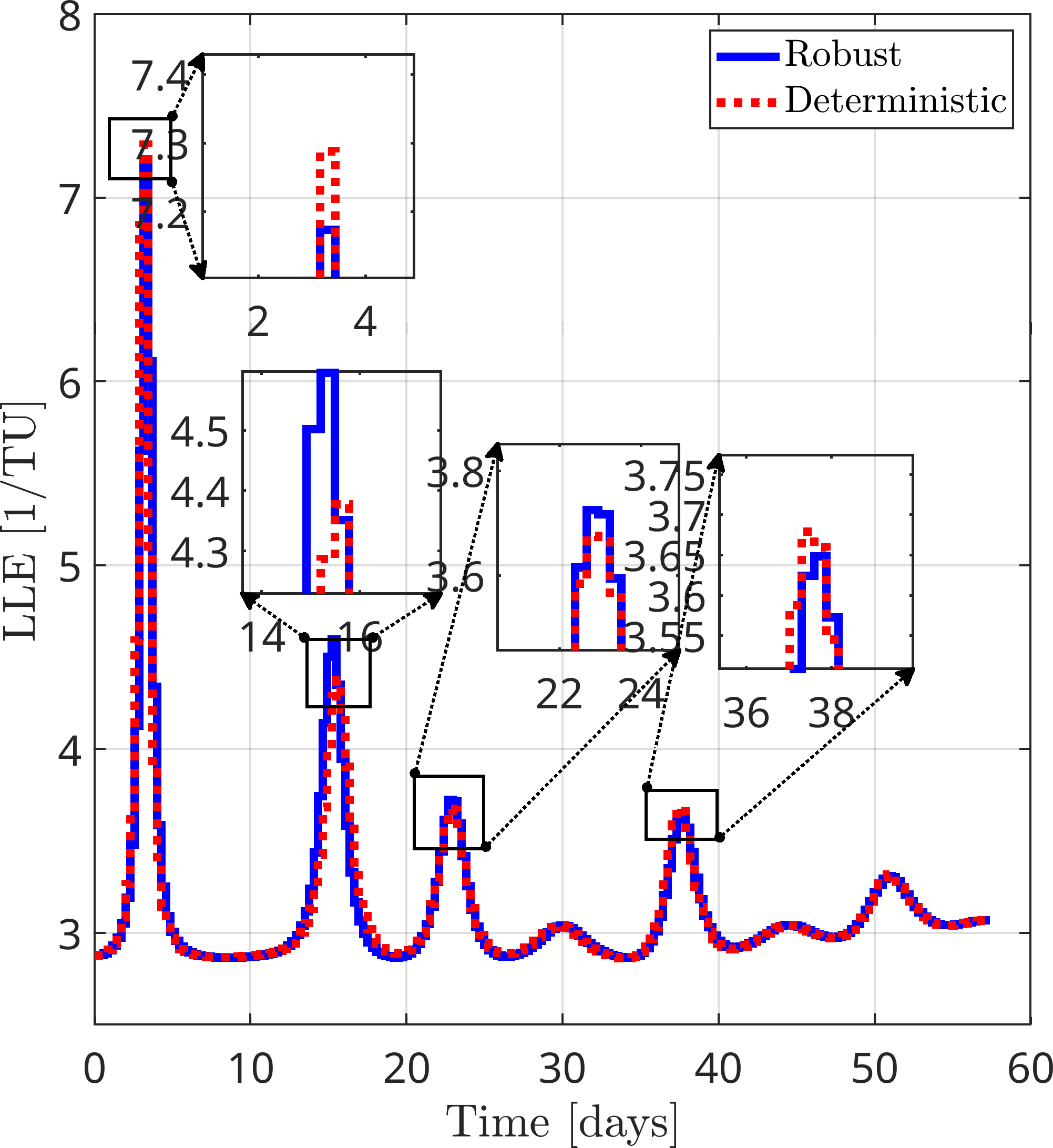}
    \caption{Comparison of the local Lyapunov Exponent along the trajectory for the deterministic and robust trajectories}
    \label{fig:LLE-comparison}    
\end{figure}

Figure \ref{fig:LLE} shows the finite-time LLE plotted along the robust trajectory. 
Due to the large LLE value around the flyby, the trajectory is split into two parts
 for visualization. 
 An immediate observation is the rise in LLE as the trajectory nears the $xz$-plane; 
 when the spacecraft is closer to the moon, the gravitational influence of the moon 
 is larger, leading to a higher LLE. 
 Compare the plots with \cref{fig:NRHOtoHalo-feedback_along_trajectory}; 
 the trends in \cref{fig:LLE_firsthalf} matches that of 
 \cref{fig:NRHOtoHalo-feedback_along_trajectory}, with large TC efforts 
 occuring around the flyby, and the sensitivity of this portion of the trajectory 
 can be numerically verified through the LLE. 

Finally, \cref{fig:LLE-comparison} shows the comparison of the LLE history
between the deterministic and robust trajectories. 
\highlight{
We see that for the first LLE peak, the robust trajectory has a lower LLE than the deterministic trajectory, 
while the second peak is larger for the robust trajectory. 
Note that the definition of the LLE involves a logarithm expression; this means that 
a similar difference in the LLE values for a large LLE value (first peak) indicates 
a larger difference in the local sensitivity of the trajectory.
The LLE-based analysis is \textit{independent of the trajectory correction policy used},
since} the LLEs are computed based on the \highlight{nominal} state transition matrix along the trajectory.
\highlight{This automated adjustment of the nominal trajectory for favorable statistical properties
represents the uncertainty-aware trajectory design process depicted in \cref{fig:comparison}.}

\section{Conclusions}
This work presents sequential convex programming-based chance-constrained covariance steering for safe trajectory and controller design in cislunar space. First, we perform a thorough comparison of the recent advances in stochastic optimal control, highlighting each approach's unique capabilities. 
We show that various operational uncertainties and dynamics uncertainties can be incorporated into the optimization, which directly steers the underlying state distribution through computing optimal trajectory correction policies. Numerical examples of transfers between periodic orbits in the Earth-Moon Circular Restricted Problem demonstrate that this framework is applicable in the highly nonlinear cislunar region. 
Building on recent advances in the stochastic optimal control and robust space trajectory optimization literature, we formulate an algorithm that empirically achieves near-optimality for the nonlinear covariance steering problem, runs faster than previous methods by orders of magnitude, and can accommodate scaling issues unique to space trajectory optimization. 

\section{Acknowledgements}
This material is based upon work supported by the Air Force Office of Scientific Research under award number FA9550-23-1-0512.
N. Kumagai acknowledges support for his graduate studies from the Shigeta Education Fund.
N. Kumagai thanks members of the Oguri Research Group for helpful discussions and Ethan Foss for pointing out the inverse-free discretization method. 
\section*{Appendix} \label{sec:appendix}
\margincomment{The contents of the Appendix have been moved from the introduction.}
\cref{tab:comparison} provides a comparison on SOC approaches based on the aspects that influence the fidelity, optimality, and computational effort. 
`UQ' compares the uncertainty quantification fidelity; 
`Control Policy Optimization' compares whether the method can output an optimized control policy. 
Navigation uncertainty is also a dividing factor that is an unignorable aspect of space trajectory optimization and is included. 
`Gradient Information' refers to the method in which the algorithm accesses derivative information; 
DDP-based methods inherently utilize up to the second-order derivative of the dynamics (state transition matrix (STM) and tensor (STT)), along with sigma points through unscented transform. 
SCP-based methods, by linearizing the dynamics and formulating the subproblem as a convex optimization, take advantage of highly optimized interior-point solvers that essentially perform the task of gradient evaluation. 
Finally, we provide a comparison of the main computational efforts that limit the algorithm complexity; 
the computational time for \highlight{\cite{oguri_stochastic_2022,vargheseNonlinearProgrammingApproach2025,oguri_stochastic-primer_2022}} are based on the authors' previous work.
\begin{table}[t]
\centering
 \fontsize{9}{9}\selectfont
\caption{Comparison of Stochastic Optimal Control Approaches for Space Trajectory Optimization}
\label{tab:comparison}
\begin{tabular}{lllllllll}
\hline
\textbf{Ref.} &
  \multicolumn{3}{l}{\textbf{Method}} &
  \textbf{\begin{tabular}[c]{@{}l@{}}UQ\\ (Fidelity)\end{tabular}} &
  \textbf{\begin{tabular}[c]{@{}l@{}}Control Policy \\ Optimization\end{tabular}} &
  \textbf{\begin{tabular}[c]{@{}l@{}}Considers \\ Navigation \\ Uncertainty\end{tabular}} &
  \textbf{\begin{tabular}[c]{@{}l@{}}Gradient \\ Information\end{tabular}} &
  \textbf{\begin{tabular}[c]{@{}l@{}}Main Computational\\ Effort (Time)\end{tabular}} \\ \hline
\cite{ozaki_tube_2020} &
  \multicolumn{3}{l}{DDP} &
  \begin{tabular}[c]{@{}l@{}}UT \\ (Med.)\end{tabular} &
  \begin{tabular}[c]{@{}l@{}}yes \\ (fitting at \\ $\sigma$ points)\end{tabular} &
  no &
  \begin{tabular}[c]{@{}l@{}}STM, STT, \\ and $\sigma$ points\end{tabular} &
  \begin{tabular}[c]{@{}l@{}}requires nonlin.\\ $\sigma$ point \\ propagation \\ (unclear)\end{tabular} \\ \hline
\cite{yuan_uncertainty-resilient_2024} &
  \multicolumn{3}{l}{DDP} &
  \begin{tabular}[c]{@{}l@{}}UT \\ (Med.)\end{tabular} &
  \begin{tabular}[c]{@{}l@{}}yes (linear \\ feedback \\ on state \\ deviation)\end{tabular} &
  yes &
  \begin{tabular}[c]{@{}l@{}}STM, STT, \\ and $\sigma$ points\end{tabular} &
  \begin{tabular}[c]{@{}l@{}}requires nonlin. \\ $\sigma$ point \\ propagation \\ (hours \cite{yuan_uncertainty-resilient_2024})\end{tabular} \\ \hline
\cite{ridderhof_chance-constrained_2020,oguri_stochastic_2022} &
  \multicolumn{3}{l}{\begin{tabular}[c]{@{}l@{}}SCP\\ (block Chol.)\end{tabular}} &
  \begin{tabular}[c]{@{}l@{}}LinCov \\ (Low)\end{tabular} &
  \begin{tabular}[c]{@{}l@{}}yes (linear \\ feedback\\  on state \\ deviation)\end{tabular} &
  \begin{tabular}[c]{@{}l@{}}no \cite{ridderhof_chance-constrained_2020}\\ yes \cite{oguri_stochastic_2022}\end{tabular} &
  \begin{tabular}[c]{@{}l@{}}linearization \& \\ discretization;\\ leverages \\ convex solver\end{tabular} &
  \begin{tabular}[c]{@{}l@{}}solves large LMI\\ (hours)\end{tabular} \\ \hline
\cite{benedikter_convex_2022} &
  \multicolumn{3}{l}{\begin{tabular}[c]{@{}l@{}}SCP\\ (full cov)\end{tabular}} &
  \begin{tabular}[c]{@{}l@{}}LinCov \\ (Low)\end{tabular} &
  \begin{tabular}[c]{@{}l@{}}yes (linear \\ feedback on \\ state history \\ deviation\\ /stochastic \\ process) \\ See \cref{remark:differences}.\end{tabular} &
  no &
  \begin{tabular}[c]{@{}l@{}}linearization \& \\ discretization;\\ leverages \\ convex solver\end{tabular} &
  \begin{tabular}[c]{@{}l@{}}solves multiple\\ small LMIs\\ (seconds to minutes)\end{tabular} \\ \hline
\cite{greco_robust_2022} &
  \multicolumn{3}{l}{\begin{tabular}[c]{@{}l@{}}NLP \\ (finite diff.)\end{tabular}} &
  \begin{tabular}[c]{@{}l@{}}GMM \\ (High)\end{tabular} &
  \begin{tabular}[c]{@{}l@{}}no (B-plane \\ targeter\\ in-the-loop)\end{tabular} &
  yes &
  \begin{tabular}[c]{@{}l@{}}Ref. uses \\ finite diff.;\\ availability of \\ analytical \\ derivative \\ unclear\end{tabular} &
  unclear \\ \hline
\highlight{\cite{vargheseNonlinearProgrammingApproach2025}} &
  \multicolumn{3}{l}{\begin{tabular}[c]{@{}l@{}}\highlight{NLP} \\ \highlight{(auto diff.)} \end{tabular}} &
  \begin{tabular}[c]{@{}l@{}} \highlight{LinCov} \\ \highlight{(Low)}\end{tabular} &
  \begin{tabular}[c]{@{}l@{}}\highlight{yes}\\ \highlight{(parameterization} \\ \highlight{of linear}\\ \highlight{feedback)}\end{tabular} &
  \highlight{no} &
  \begin{tabular}[c]{@{}l@{}}\highlight{STM, STT;}\\ \highlight{calculated via} \\ \highlight{automatic} \\ \highlight{differentiation} \end{tabular} &
  \highlight{tens of minutes} \\ \hline
\cite{oguri_stochastic-primer_2022} &
  \multicolumn{3}{l}{Indirect} &
  \begin{tabular}[c]{@{}l@{}}LinCov \\ (Low)\end{tabular} &
  no (open-loop) &
  no &
  \begin{tabular}[c]{@{}l@{}}Ref. uses \\ finite diff.\end{tabular} &
  \begin{tabular}[c]{@{}l@{}}solves 2PBVP \\ for mean and \\ covariance\\ (minutes)\end{tabular} \\ \hline
\end{tabular}
\end{table}

\bibliography{references}

\end{document}